\documentclass[12pt]{iopart}
\expandafter\let\csname equation*\endcsname\relax
\expandafter\let\csname endequation*\endcsname\relax 
\usepackage{amsthm,enumitem,color,tikz,iopams,datetime,dsfont,siunitx}
\usepackage{mathtools,amsmath,graphicx,subcaption,multirow,cite}
\DeclareGraphicsExtensions{.png,.pdf}
\usetikzlibrary{positioning}

\usepackage[titletoc]{appendix}
\definecolor{linkColor}{rgb}{0.4, 0.4, 0.4}
\usepackage[colorlinks, allcolors=linkColor, bookmarks=True]{hyperref}
\makeatletter
\newcommand{\customlabel}[2]{%
   \protected@write \@auxout {}{\string \newlabel {#1}{{#2}{\thepage}{#2}{#1}{}} }%
   \hypertarget{#1}{}
}
\makeatother

\newtheorem{definition}{Definition}[section]\newtheorem{theorem}{Theorem}[section]
\newtheorem{lemma}{Lemma}[section]
\makeatletter\let\c@definition\c@theorem\let\c@definition\c@theorem\makeatother
\makeatletter\let\c@lemma\c@theorem\let\c@lemma\c@theorem\makeatother
\makeatletter\let\c@example\c@theorem\let\c@example\c@theorem\makeatother
\newtheorem*{theorem*}{Theorem}
\newtheorem{approximation}{Approximation}

\newcommand{\F}[1]{\mathds{ #1 }}\newcommand{\C}[1]{\mathcal{ #1 }}
\newcommand{\norm}[1]{\left\| #1 \right\|}\newcommand{\sfrac}[2]{\,^{#1}\!\!/\!_{#2}}
\renewcommand{\IP}[2]{#1 \vcenter{\hbox{\resizebox{6pt}{!}{\ensuremath\cdot}}} #2}
\newcommand{\1}[0]{\F{1}}
\newcommand{\op}[1]{\operatorname{#1}}
\renewcommand{\bar}{\overline}\renewcommand{\hat}{\widehat}
\let\vecf\vec
\renewcommand{\vec}{\boldsymbol}\renewcommand{\epsilon}{\varepsilon}
\renewcommand{\i}{\mathrm{i}}
\DeclareMathOperator*{\E}{\mathop{\raisebox{-1pt}{\normalfont\large \ensuremath{\F{E}}}}}
\newcommand{\fitwidth}[1]{\noindent\resizebox{.99\textwidth}{!}{#1}\newline }
\newcommand{\fitwidtheq}[1]{\fitwidth{\ensuremath{#1}} }

\DeclareMathOperator*{\argmin}{\mathop{argmin}}
\newcommand{\cross}[1]{\left[#1\right]_\times}
\newcommand{\D}[1][\,]{#1{\rm d}}

\usepackage[normalem]{ulem}
\definecolor{darkgreen}{rgb}{0.0, 0.5, 0.0}
\newcommand{\edit}[2]{{\ifmmode\text{\color{red}\sout{\ensuremath{#1}}}\else {\color{red} \sout{#1}}\fi} {\color{darkgreen} #2}}
\renewcommand{\edit}[2]{#2}

\newcommand{\SUPP}[2]{#1} 

\begin{document}
\title[Scanning electron diffraction tomography of strain]{Scanning electron diffraction tomography of strain}

\author{Robert Tovey$^1$, Duncan N. Johnstone$^2$, Sean M. Collins$^{2,3}$, William R. B. Lionheart$^4$, Paul A. Midgley$^2$, Martin Benning$^5$, Carola-Bibiane Sch\"onlieb$^1$}
\address{1 Centre for Mathematical Sciences, University of Cambridge, Cambridge CB3 0WA, UK}
\address{2 Department of Materials Science and Metallurgy, University of Cambridge, 27 Charles Babbage Road, Cambridge CB3 0FS, UK}
\address{3 School of Chemical and Process Engineering and School of Chemistry, University of Leeds, Leeds LS2 9JT, UK}
\address{4 Department of Mathematics, University of Manchester, Manchester M13 9PL, UK}
\address{5 School of Mathematical Sciences, Queen Mary University of London, London E1 4NS, UK}
\ead{RobTovey@maths.cam.ac.uk}

\begin{abstract}
Strain engineering is used to obtain desirable materials properties in a range of modern technologies. Direct nanoscale measurement of the three-dimensional strain tensor field within these materials has however been limited by a lack of suitable experimental techniques and data analysis tools. Scanning electron diffraction has emerged as a powerful tool for obtaining two-dimensional maps of strain components perpendicular to the incident electron beam direction. Extension of this method to recover the full three-dimensional strain tensor field has been restricted though by the absence of a formal framework for tensor tomography using such data. Here, we show that it is possible to reconstruct the full non-symmetric strain tensor field as the solution to an ill-posed tensor tomography inverse problem. We then demonstrate the properties of this tomography problem both analytically and computationally, highlighting why incorporating precession to perform scanning precession electron diffraction may be important. We  establish a general framework for non-symmetric tensor tomography and demonstrate computationally its applicability for achieving strain tomography with scanning precession electron diffraction data.
\end{abstract}
\vspace{2pc}
\noindent{\it \edit{k}{K}eywords\/}: transverse ray transform, computed tomography, tensor tomography, scanning precession electron diffraction, 4D-STEM, strain mapping, strain tomography

\maketitle
	
\section{Introduction}

In this paper, we examine whether it is theoretically possible to recover the three-dimensional strain tensor field within a material using scanning electron diffraction (SED) data and tensor tomography methods.

Nanoscale strain is widely used to engineer desirable materials properties, for example, improving field effect transistor (FET) performance \cite{chu_strain_2009}, opening a bulk bandgap in topological insulator systems \cite{hsieh_tunable_2009} and enhancing ferroelectric properties \cite{schlom_strain_2007}. Strain also arises around crystal defects, which further affect materials properties. The strain tensor field is a rank 2 symmetric tensor field in three-dimensional space and is therefore fully described by 6 components at every 3D coordinate. 3D strain reconstruction of one or more strain components has been achieved using X-ray diffraction techniques, including: coherent Bragg diffractive imaging \cite{Pfeifer2006,Robinson2009,Newton2010}, micro-Laue diffraction using a differential aperture \cite{Larson2002}, and diffraction from polycrystalline specimens combined with back-projection methods \cite{Korsunsky2006, Korsunsky2011}. The spatial resolution of these X-ray techniques is however limited to ca. \SIrange{20}{100}{\nano\meter} and sub-\SI{10}{\nano\meter} resolution strain mapping is therefore dominated by (scanning) transmission electron microscopy ((S)TEM) techniques \cite{hytch_observing_2014}. These techniques include: imaging at atomic resolution \cite{Hytch2003, Galindo2007}, electron holography \cite{Hytch2011} and SED \cite{usuda_strain_2005, beche_strain_2013}. Amongst these, 3D strain has been assessed by atomic resolution tomography \cite{Goris2015} and in a proof-of-principle reconstruction of a single strain component using SED \cite{johnstone_nanoscale_2017}. In 2D strain mapping, SED has emerged as a particularly versatile and precise approach to strain mapping with few nanometre resolution \cite{cooper_combining_2015, cooper_high-precision_2017}.

\begin{figure}[ht]
\centering
	\begin{subfigure}{.35\textwidth}
		\includegraphics[width=\textwidth, trim={0 0 1423 0},clip]{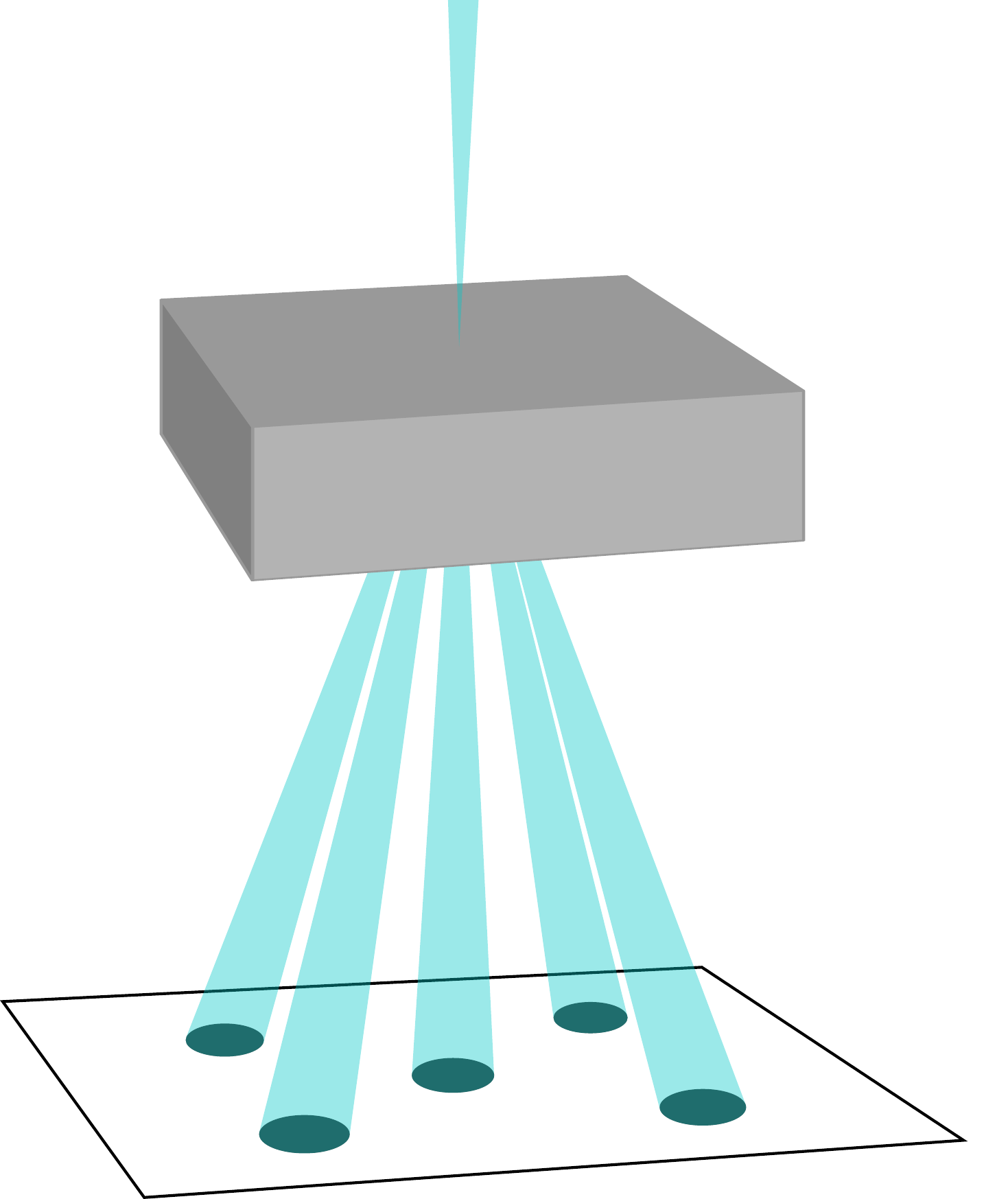}
		\caption{Electron diffraction}
	\end{subfigure}\hspace{20pt}
	\begin{subfigure}{.35\textwidth}
		\includegraphics[width=\textwidth, trim={467 0 956 0},clip]{multi-scale_EM}
		\caption{Scanning electron diffraction}
	\end{subfigure}
\caption{Schematic illustration of scanning electron diffraction. (a) A focused electron probe passes through a sample and a 2D electron diffraction pattern is recorded. If the sample is crystalline then the diffraction pattern comprises intense Bragg disks at positions related to the spacing of atomic planes in the crystal, as shown. (b) The probe is scanned over the sample and a diffraction pattern recorded at each position.}
\label{fig: sped}
\end{figure}

SED is a 4D-STEM technique \cite{ophus_four-dimensional_2019} based on the acquisition of a 2D transmission electron diffraction pattern at every probe position as a focused electron probe is scanned across the specimen in a 2D scan, as shown in Figure~\ref{fig: sped}. Each 2D electron diffraction pattern comprises intense Bragg disks (see Figure~\ref{fig: sped}) at scattering angles, $\theta_{B}$, related to the spacing of atomic planes, $d_{hkl}$, by Bragg's law,
\begin{equation*}
    \lambda = 2 d_{hkl} \sin \theta_{B}
\end{equation*}
where, $\lambda$ is the radiation wavelength \edit{}{and $hkl$ are the \emph{Miller indices} of the atomic plane defined as reciprocal of the intercepts of that plane with the crystal basis vectors}. For high-energy (ca. \SIrange{60}{300}{\keV}) incident electrons, the corresponding de Broglie wavelength is ca. \SIrange{0.0487}{0.0197}{\angstrom}. 

The largest $d_{hkl}$ values are typically \textless \SIrange{3}{5}{\angstrom} for most metals and ceramics\edit{and t}{. T}herefore the smallest scattering angles are typically ca. 10 mrad and the corresponding Bragg disks are associated with atomic planes near parallel to the incident electron beam direction. Strain alters the spacing of atomic planes and consequently causes the position of the Bragg disks to change and to be blurred if a varying strain field is sampled along the beam path. Components of strain perpendicular to the incident beam direction may therefore be mapped by tracking the position of Bragg disks as a function of probe position \cite{usuda_strain_2005}. 

Strain maps of three path averaged components of the strain tensor in 2D have been reported from a wide range of materials via SED \cite{cooper_high-precision_2017, haas_high_2017, pekin_situ_2018, bonef_interfacial_2016}. Determining the position of the Bragg disks in each diffraction pattern is a critical step, which has been explored in recent literature with cross-correlation based disk finding approaches achieving the best accuracy and precision \cite{beche_improved_2009, rouviere_improved_2013, pekin_optimizing_2017, zeltmann_patterned_2020}. The incorporation of double-conical electron beam rocking \cite{Vincent1994}, to record scanning precession electron diffraction (SPED) data, has further been demonstrated to improve precision both numerically \cite{Mahr2015} and experimentally \cite{cooper_combining_2015}. However, progress towards 3D strain mapping using S(P)ED data has so far been limited to a proof-of-principle reconstruction of one strain component in 3D \cite{johnstone_nanoscale_2017} by the lack of a framework for three-dimensional strain tensor field reconstruction using S(P)ED data. In this work, we establish such a framework for three-dimensional reconstruction of the full strain tensor field from S(P)ED data via consideration of an analytical forward model (see Sections~\ref{sec: background} to \ref{sec: approximations}). We show that a linearised approximation coincides with a non-symmetric tensor tomography problem and use this to demonstrate recovery of \edit{the full strain field}{both the symmetric infinitesimal strain tensor field and the non-symmetric rotation field to recover the full non-symmetric tensor field, i.e. the deformation gradient tensor field}. 

In Section~\ref{sec: Tensor Tomography} we show that the linearised model for diffraction data coincides with the transverse ray transform on tensor fields with non-symmetric tensors. Tensor tomography is a well-established area of mathematics with many physical applications \cite{Sharafutdinov1994}, including \edit{}{in the application of strain reconstruction in photoelastic tomography. Experimental and numerical details can be found in \cite{Tomlinson2006} and \cite{Lionheart2009} respectively. In two dimensions, the work of \cite{Puro2007} considers symmetric strain tomography, \cite{Derevtsov2015} considers symmetric tensors of arbitrary rank, and \cite{Gullberg1999} also considers the non-symmetric strain case.} 
\edit{(symmetric) strain tomography for polycrystalline materials [34].}{Symmetric strain tomography of polycrystalline materials in three dimensions is explored by \cite{Lionheart2015}.} 
Results for tensor fields of non-symmetric tensors are well explored in dimensions greater than three \cite{Novikov2007, Abhishek2020} and also in three dimensions for a general Riemannian geometry \cite{Holman2013}. We combine and simplify these results into a single cohesive argument to clarify properties of the inverse problem in the physical scenario. Finally, we validate our \edit{reasoning}{approach} computationally in Section~\ref{sec: numerics} with a range of  complex forward models for scanning electron diffraction. Our analytical and numerical results establish a robust framework for three-dimensional strain tensor field reconstruction using S(P)ED data.

\section{Fundamentals of electron diffraction for strain mapping} \label{sec: background}
\subsection{Notation}
We define that:
\begin{itemize}
	\item $\Psi_p\colon\F R^2\to \F C$ is the \emph{probe function} representing the incident electron wave with wavelength $\lambda>0$. This is extended to 3D by $\Psi_p(x,y,z) = \Psi_p(x,y)$, assuming \edit{the in free space }{}the wave \edit{would}{does} not change significantly over the specimen thickness.
	\item $u\colon\F R^3\to\F R$ and any subscripted derivation denotes the electrostatic potential of the specimen under investigation.
	\item $\vec x = (x,y,z)$ gives the standard coordinates on $\F R^3$ (direct space).
	\item $\vec K = (\vec k, k_z) = (k_x,k_y,k_z)$ gives the standard coordinates on $\F R^3$, in Fourier space (reciprocal space), which we treat as $\F R^3 = \F R^{2+1}$ interchangeably, where $\vec k$ is the 2D coordinate on the detector.
	\item $\C F[u](\vec K) = \int_{\F R^3} u(\vec x)\exp(-\i \IP{\vec x}{\vec K})\D\vec x$ denotes the Fourier transform.
	\item $\gamma = \left(\begin{smallmatrix}1&0\\0&1\\0&0\end{smallmatrix}\right)$ is the natural lifting from $\F R^2$ to the plane $k_z=0$ in $\F R^3$. i.e. $\gamma\vec K = \vec k, \gamma^\top \vec k = (\vec k,0)$.
	\item $D\colon\F R^2\to\F R$ and subscripted derivatives denote a two-dimensional diffraction pattern recorded on a flat detector.
	\item $\1_A(\vec x) = \begin{cases}
	1 & \vec x \in A \\ 0 &\text{else} \end{cases}$ denotes the indicator on the set $A$.
	\item $\delta_{\vec p}$ is the Dirac delta centred at $\vec p$. Formally, this is considered as a tempered distribution on the space of Schwartz functions. Informally, we write it as a function $\delta_{\vec p}(\vec x)$ where it aids clarity.
\end{itemize}

\subsection{General electron diffraction}

Forward models of electron diffraction typically assume coherent elastic scattering of fast incident electrons by the static Coulomb (electrostatic) potential, $u$, of a specimen in the far-field (Fraunhofer) diffraction limit \cite{Peng1996}. This neglects inelastic and mixed elastic/inelastic scattering contributions \cite{Wang2013}. Various further assumptions lead to a range of analytical and numerical forward models \cite{Kirkland2010, Paganin2006}. A key distinction is between \emph{kinematical} models in which a single electron undergoes at most one scattering event, and \emph{dynamical} models in which multiple scattering events can occur. While the latter is more physically accurate, the former captures many important aspects. We proceed here with a kinematical model for analysis and validate our results computationally using both kinematical and dynamical models in Section \ref{sec: numerics}.

The kinematical diffraction pattern, $D$, may be written as
\begin{equation}\label{eq: sphere model}
D(\vec k) = \left|\C F[\Psi_pu]\right|^2(\vec k, k_z(\vec k)), \qquad k_z(\vec k) = 2\pi\lambda^{-1}-\sqrt{4\pi^2\lambda^{-2}-|\vec k|^2}
\end{equation}
where $\vec k$ is the 2D coordinate on the detector, $\Psi_pu$ is a pointwise multiplication between the probe function and specimen potential, and the incident electron beam is chosen to be travelling parallel to the $z$-axis. The Fourier transform is only observed on the so\edit{}{-}called Ewald sphere,
$$\{(\vec k,k_z(\vec k)) : |\vec k|\leq 2\pi\lambda^{-1}\} = \left\{\vec K : K_z\leq 2\pi\lambda^{-1}, \left|\vec K - \begin{pmatrix} 0&0&2\pi\lambda^{-1}\end{pmatrix}\right| \leq 2\pi\lambda^{-1}\right\}.$$
on which Bragg's law is satisfied geometrically.

\subsection{Electron diffraction from a nanocrystal}

A crystal is defined as a material that has an essentially sharp diffraction pattern \cite{AC1992}. This means that most of the diffracted intensity is concentrated in Bragg peaks at particular positions. Diffraction patterns recorded from crystals are therefore highly structured and essentially sparse, though in practice diffuse scattering is always observed in addition to the Bragg peaks. This leads to the definition of an ideal defect-free crystal inferred from \eqref{eq: sphere model} to be a material with a sparse Fourier transform.
\begin{definition} The electrostatic potential of an \emph{ideal crystal} $u_0\colon\F R^3\to\F R$ satisfies 
    $$\C F[u_0] = \sum_{i=1}^\infty w_i\delta_{\vec p_i}$$
    for some weightings $w_i\in\F C$ and \emph{Bragg peaks} $\vec p_i\in\F R^3$.    
\end{definition}
For \textit{conventional crystals} the set of Bragg peaks are structured such that there exist some \textit{reciprocal lattice vectors} $\vec a^*$, $\vec b^*$ and $\vec c^*\in\F R^3$ such that 
$$\{\vec p_i: i\in\F N\} = \{h\vec a^*+k\vec b^*+l\vec c^* : h,k,l\in\F Z\}.$$
\edit{}{The classical description indexes over three integer indices, as on the right-hand side. In this work we only use a single positive index to compress notation, it will not be necessary to record the precise Miller indices for each $\vec p_i$.} The reciprocal lattice vectors link the sparsity of $\C F[u_0]$ and the periodicity of $u_0$ by the following Lemma.
\begin{lemma}\label{thm: spike locations}
    If $u_0$ has reciprocal lattice vectors $\vec a^*$, $\vec b^*$ and $\vec c^*$ then, with $V \coloneqq \op{det}((\begin{smallmatrix}\vec a^*&\vec b^*& \vec c^*\end{smallmatrix}))$, $u_0$ is periodic in the vectors
    $$ \vec a = \frac{2\pi}{V}\vec b^*\times\vec c^*, \qquad \vec b = \frac{2\pi}{V}\vec c^*\times\vec a^*, \quad \text{and}\quad \vec c = \frac{2\pi}{V}\vec a^*\times\vec b^*.$$
    In other words, there exists a repeating unit $v\colon[0,\max(|\vec a|,|\vec b|,|\vec c|)]^3\to\F R$, referred to as the \emph{unit cell} such that 
    $$u_0(\vec x) = \sum_{j_1,j_2,j_3=-\infty}^{\infty} v(\vec x + j_1\vec a+j_2\vec b+j_3\vec c).$$	
\end{lemma}

This is a standard result in crystallography and a short proof is provided in \ref{app: spike locations}. The key point is that both direct space and reciprocal space representations of the idealised $u_0$ are equivalent however \edit{}{each }have distinct advantages. The direct space parametrisation shows that ideal conventional crystals are built of repeating unit cells which will be the basis of a discretisation in Section~\ref{sec: physical model with strain}, however we generally adopt the reciprocal space (Fourier) parametrisation, which is natural for analysing diffraction patterns.

The diffraction pattern from a crystal truncated to a cubic volume can be computed directly using \eqref{eq: sphere model} as,
\begin{lemma}\label{thm: exact single crystal diffraction}
	\edit{If}{Let} $u(\vec x) = u_0(\vec x)\F V(\frac{\vec x}{\rho})$ where 
	$$\F V(\vec x) \coloneqq \1_{[-\sfrac12,\sfrac12]^3}(\vec x) = \begin{cases}1 & \vec x\in [-\sfrac12,\sfrac12]^3 \\ 0 &\text{ else}\end{cases}$$
	and $\rho>0$ is the width/depth of the crystal $u$\edit{ then}{. Define the function $\op{sinc}\colon\F R^d\to\F R$ by $$\op{sinc}(k) \coloneqq \frac{\sin(k)}{k}\quad\text{when } d=1, \qquad \text{and}\qquad \op{sinc}(\vec k) \coloneqq \op{sinc}(k_x)\op{sinc}(k_y)\ldots$$ when $d>1$. Then, } 
	$$ D(\vec k) = \left|\sum_{i=1}^\infty w_i\rho\op{sinc}(\rho[k_z(\vec k)-p_{i,z}]) f(\vec k-\gamma^\top \vec p_i)\right|^2 $$
	where
	\begin{equation}\label{eq: def f}
	f(\vec k) \coloneqq \rho^2 \C F[\Psi_p]\star [\op{sinc}(\rho\cdot)](\vec k) = \rho^2\int_{\F R^2} \C F[\Psi_p](\vec k-\vec k')\op{sinc}(\rho\vec k') \D\vec k'.
	\end{equation}
	
\end{lemma}

This is also a special case of Lemma~\ref{thm: full strain model} for which a proof is provided in \ref{app: approximations}. Single crystal diffraction patterns are therefore essentially a sum of spikes centred at points $\gamma^\top \vec p_i = (p_{i,x},p_{i,y})$ and with a diffracted intensity proportional to $\left|w_i\right|^2$, which is therefore related to the \textit{structure factor} in crystallography. The shape of \edit{}{the }spikes is dictated by $\Psi_p$, encoded in $f$. 

\subsection{Precession electron diffraction}\label{sec: precession intro}
Precession electron diffraction (PED) patterns are recorded using a dynamic ``double conical beam-rocking geometry'' \cite{Vincent1994}, as illustrated in Figure~\ref{fig: ped}. The incident electron beam is tilted away from the optic axis ($z$-axis) of the microscope by a precession angle, $\alpha$, and rotated about the $z$-axis above the specimen, with a counter rotation performed below the specimen to integrate diffracted intensities on the detector. Typical values for the precession angle in SPED experiments are ca. $\alpha = \SIrange{0.5}{2}{\degree}$.

Algebraically, PED can be equivalently described by switching to the frame of reference such that the electron beam is stationary and the sample is rotating around the $z$-axis. We can then write an expression for a PED pattern, $D_\alpha(\vec k)$, as,
\begin{multline}\label{eq: precessed sphere model}
D_\alpha(\vec k) = \E_t\left\{\left|\C F[\Psi_pu(R_t\vec x)]\right|^2(\vec k,k_z(\vec k))\text{ such that } t\in[0,2\pi) \text{ and }
\phantom{\left(\begin{smallmatrix} \cos(t)\\-\sin(t)\\1\end{smallmatrix}\right)}\right.
\\ \left.R_t = \left(\begin{smallmatrix} \cos(t)&\sin(t)&0\\-\sin(t)&\cos(t)&0\\0&0&1\end{smallmatrix}\right)\left(\begin{smallmatrix} 1&0&0\\0&\cos(\alpha)&\sin(\alpha)\\0&-\sin(\alpha)&\cos(\alpha)\end{smallmatrix}\right)\left(\begin{smallmatrix} \cos(t)&-\sin(t)&0\\\sin(t)&\cos(t)&0\\0&0&1\end{smallmatrix}\right) \right\}.
\end{multline}
\edit{}{The expectation $\E_t$ is defined formally in \ref{app: probability background}.} Experimentally, precession is a `summing' of diffraction patterns although here we write it equivalently as an average or expectation.
There is no intrinsic randomness but considering $t\sim \op{Uniform}[0,2\pi)$ ensures that the expectation is the desired value. Heuristically, standard diffraction patterns sample the Fourier transform pointwise, which is very sensitive to variations in intensity both due to the intersection of the Ewald sphere with the Fourier space and dynamical scattering, as can be seen in Figure~\ref{fig: first examples}.a. PED patterns integrate over these intensity variations (Figure~\ref{fig: first examples}.b) leading to improved agreement with a simpler and `more kinematical' model (Figure~\ref{fig: first examples} parts c and d). This idea will form a key approximation described in Section~\ref{sec: approximations}. Overall, Figure~\ref{fig: first examples} shows that the kinematical model has a much faster decay due to the Ewald sphere, however, it is also much less sensitive to oscillations in intensity. After precession, the two models agree well.

\begin{figure}[h]
	\hspace{-25pt}\begin{tikzpicture}
	\node at (0,0) {\includegraphics[width=.9\textwidth]{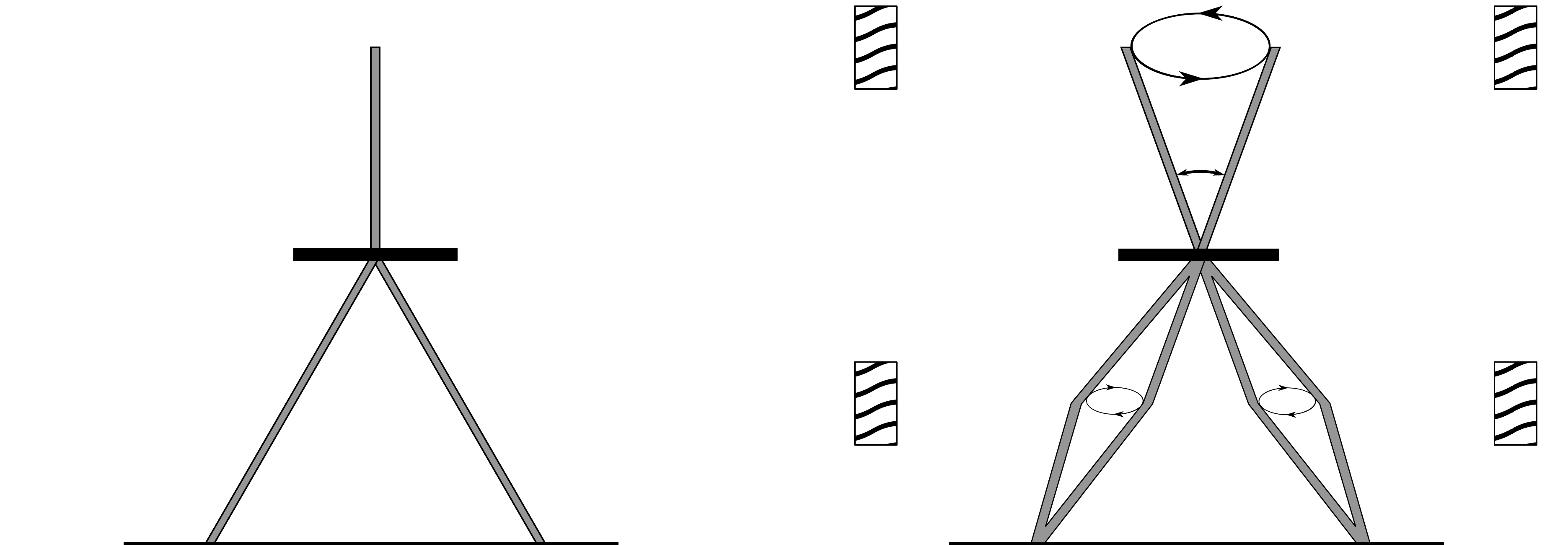}};
	\node at (-3.7,2.8) {standard beam};
	\node at (3.7,2.8) {precessing beam};
	\node at (-3.7,-2.7) {diffracted beams};
	\node at (+3.7,-2.7) {diffracted beams};
	\node[text width=4cm] at (9,2.1) {pre-specimen deflection coils};
	\node[text width=4cm] at (9,-1.2) {post-specimen deflection coils};
	\node at (3.8,1.2) {$2\alpha$};
	\draw[-] (0,3) -- (0,-3);
	\end{tikzpicture}
	\caption{Schematic of double conical beam-rocking geometry used to record precession electron diffraction (PED) patterns. The electron beam is rocked around the optic axis above the specimen and counter-rocked below the specimen to record electron diffraction patterns containing Bragg disks integrated over the rocking condition.}
	\label{fig: ped}
\end{figure}
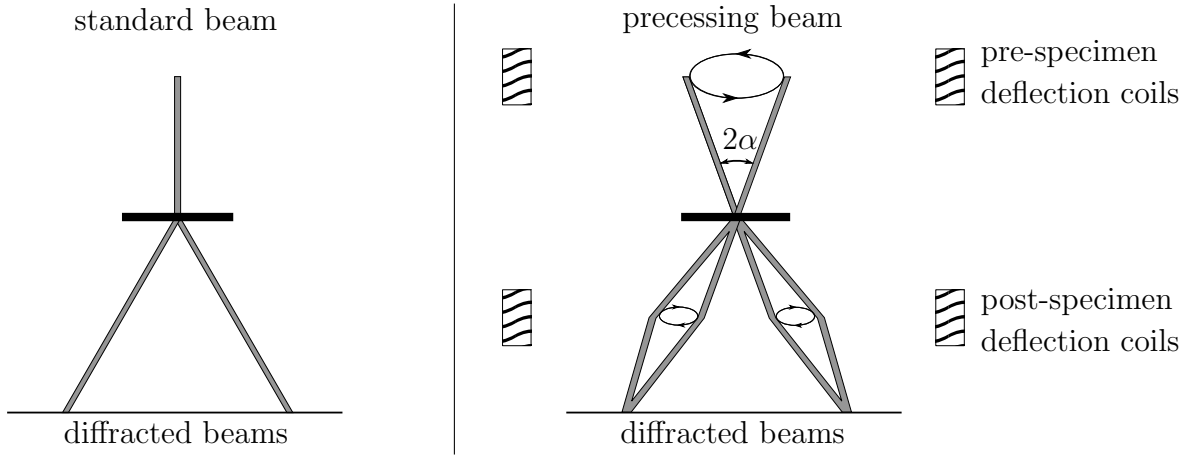
\begin{figure}\centering
	\begin{subfigure}{.4\textwidth}
		\fbox{\includegraphics[width=.9\textwidth, trim={140 50 540 55},clip]{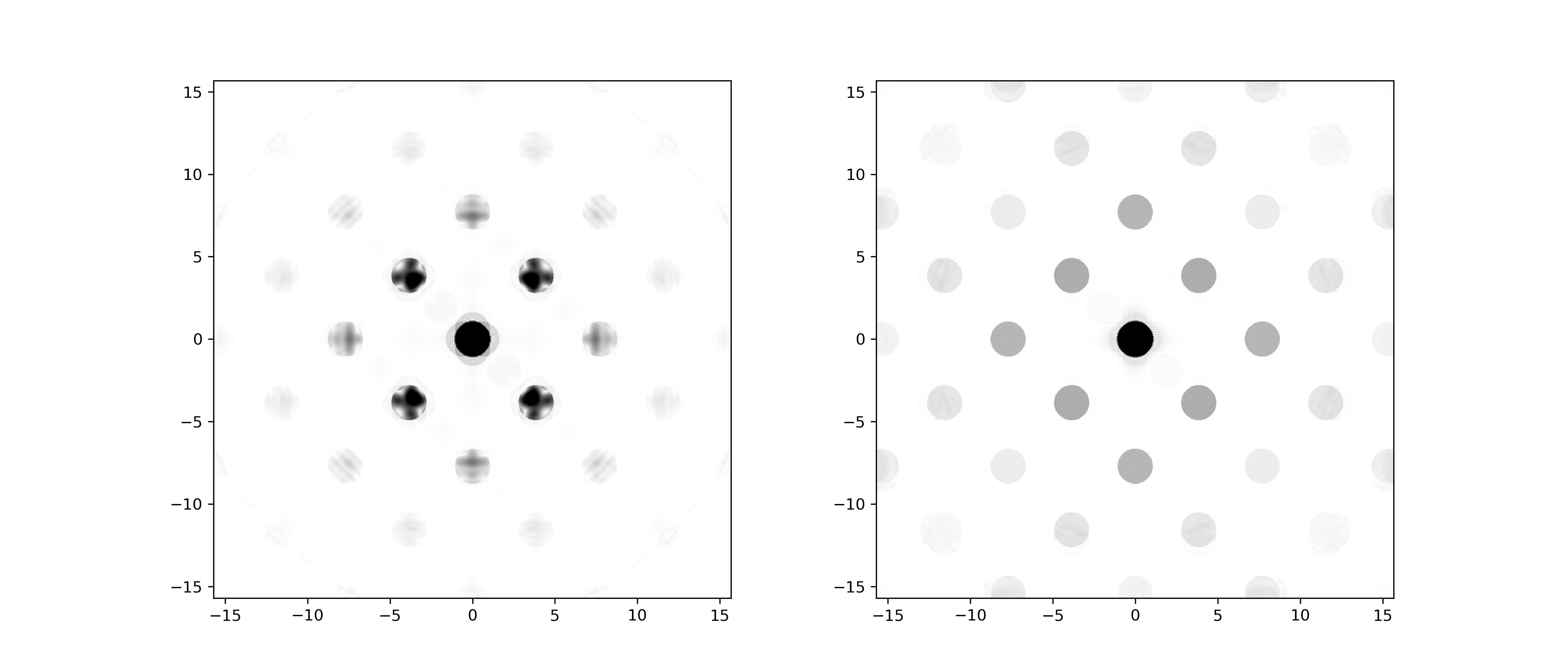}}
		\caption{Dynamical simulation}
	\end{subfigure}
	\begin{subfigure}{.4\textwidth}
		\fbox{\includegraphics[width=.9\textwidth, trim={565 50 115 55},clip]{multislice_demos}}
		\caption{Dynamical simulation with precession}
	\end{subfigure}
	\begin{subfigure}{.4\textwidth}
		\fbox{\includegraphics[width=.9\textwidth, trim={565 50 115 55},clip]{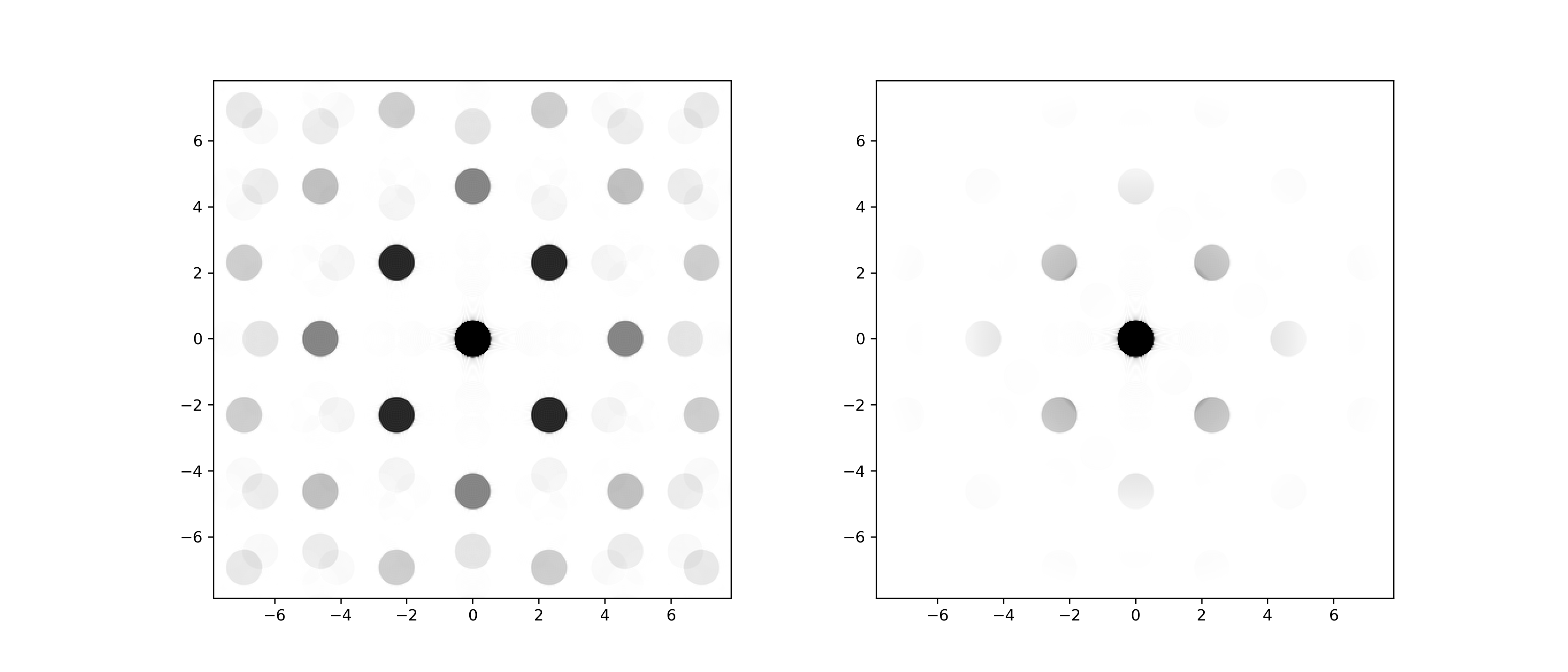}}
		\caption{Kinematical simulation}
	\end{subfigure}
	\begin{subfigure}{.4\textwidth}
		\fbox{\includegraphics[width=.9\textwidth, trim={565 50 115 55},clip]{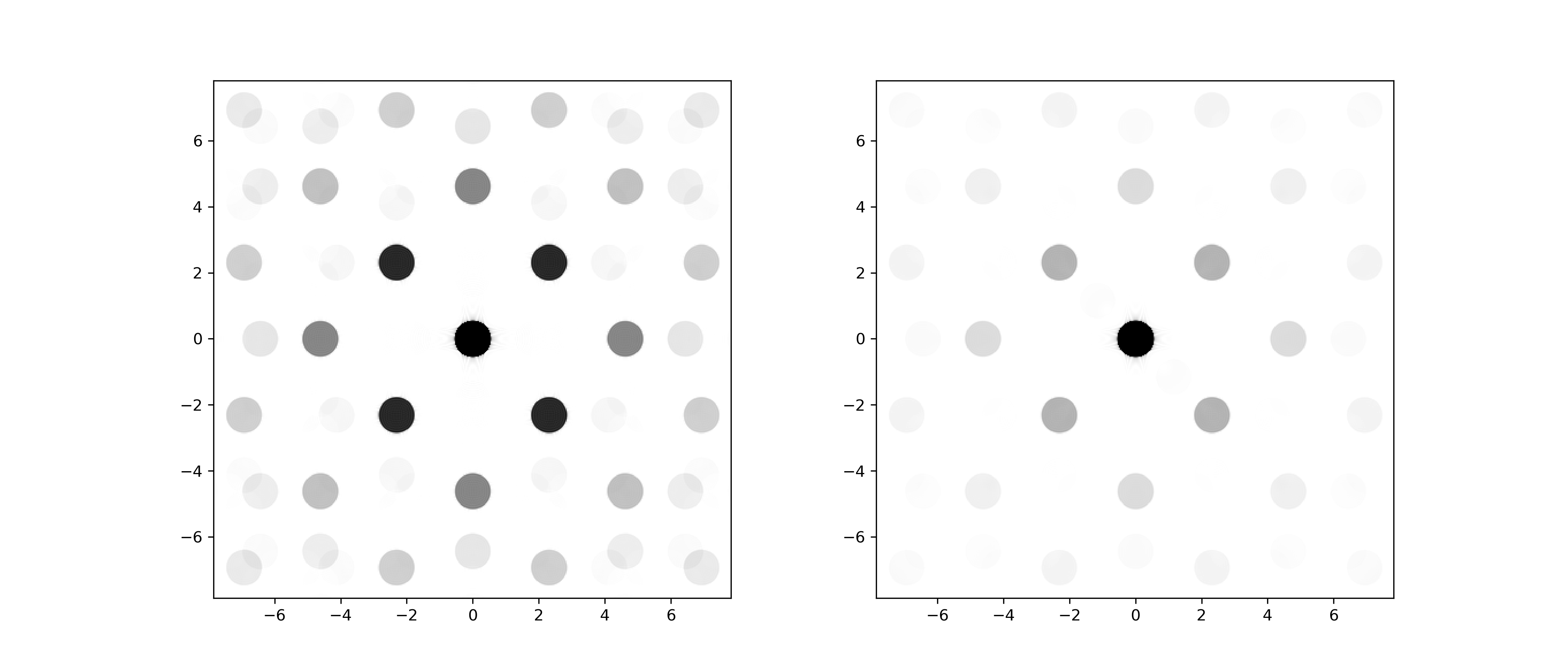}}
		\caption{Kinematical simulation with precession}
	\end{subfigure}
	\caption{Simulations of diffraction patterns from an unstrained Silicon crystal \edit{in the [001] orientation}{viewed along the [001] crystal direction i.e. along the crystallographic c-axis} (details in Section~\ref{sec: numerical simulation results}). (a) shows a dynamical simulation where complex spot inhomogeneities can be seen. (b) shows that with precession, the intensities in the dynamical simulation become much more homogeneous. (c)/(d) show kinematical simulations without/with precession. Note that precessed images qualitatively agree very closely with each other.}
	\label{fig: first examples}
\end{figure}

\subsection{Technical assumptions for diffraction imaging}\label{sec: assumptions}

We require five technical assumptions on the experimental setup which are used to justify the images seen in Figure~\ref{fig: first examples}.

\begin{enumerate}[leftmargin=\widthof{Asumption 1:+\quad},label=\textbf{Assumption \arabic*:},ref=Assumption \arabic*]
	\item \emph{small wavelength}, i.e. $\lambda\sim \SI{0.01}{\angstrom}$. \label{ass: high energy}
	\item \emph{thick crystals}, i.e. $\rho\sim\SI{1000}{\angstrom}$. \label{ass: thick crystals}
	\item \emph{non-overlapping Bragg disks}, i.e. there exists $\bar r>0$ such that for all  $i_1\neq i_2$: $|\vec p_{i_1}-\vec p_{i_2}|>2\bar r$, and $|\vec k|>\bar r\implies f(\vec k)= 0$. \label{ass: separated spots}
	\item \emph{symmetrical spots}, i.e. $f(-\vec k) = f(\vec k) \text{ for all } \vec k$. \label{ass: symmetric}
	\item \emph{narrow beam}, i.e. $|\vec x| > r \sim \SI{30}{\angstrom} \implies \Psi_p(\vec x) = 0$. \label{ass: local signal}
\end{enumerate}
These assumptions can all be met readily in typical SPED experiments by configuring the electron optics while considering the crystal lattice parameters. \ref{ass: high energy} and \ref{ass: thick crystals} justify that Figure~\ref{fig: first examples}.d is a sum of spots of shape $\C F[\Psi_p]$ \edit{centered}{centred} at $(p_{i,x},p_{i,y})$ with $p_{i,z}=0$. \ref{ass: separated spots} and \ref{ass: symmetric} state that the spots are symmetrical and well separated. For the kinematical model, \ref{ass: symmetric} is equivalent to the symmetry of $\Psi_p$ and precession helps to mitigate dynamical effects breaking this condition. Finally, \ref{ass: local signal} dictates the final spatial resolution of the strain mapping. The width of the support of $\Psi_p$ acts like a point-spread function on the model. Without treating this explicitly, we need to assume that the beam is contained within a single column of voxels in the final strain mapped volume.

\section{Electron diffraction from a strained crystal}\label{sec: physical model with strain}

There are two standard models of deformed crystals \cite{Howie1968}, the continuum (deformable ion) model, which we use here for our analysis, and the discrete or atomic (rigid ion) model, which we use in Section \ref{sec: numerics} for computation. In the continuum model, a deformed crystal is defined as follows.

\begin{definition}\label{def: strain}
	\edit{Let $u\colon\F R^3\to\F R$ be an electrostatic potential and $\vecf R\colon\F R^3\to\F R^3$ a deformation defined by}{If $u'\colon\F R^3\to\F R$ is defined pointwise by}
	$$u'(\vec x) = u(\vec x+\vecf R(\vec x))$$
	\edit{where $\vecf R$ the \emph{displacement map}}{for some crystal $u\colon\F R^3\to\F R$ and \emph{displacement map} $\vecf R\colon\F R^3\to\F R^3$, then we say that $u'$ is a \emph{deformed crystal}} and $\vecf F = \nabla\vecf R$ is the \emph{displacement gradient tensor field}.
\end{definition}

The \textit{strain tensor field}, $\vecf \epsilon$ is the symmetric part of the displacement gradient tensor, which we obtain as $\vecf \epsilon = \tfrac12(\vecf F + \vecf F^\top )$, while the skew part describes the rotation field. Since the diffraction pattern is sensitive to both strain and rotation, we consider general displacement gradient tensor fields, which can then be decomposed into strain and rotation parts following tomographic reconstruction.

\subsection{Strain in Fourier space}

Understanding the Fourier transform of a deformed crystal is the key to understanding the diffraction pattern produced by that deformed crystal within the kinematical model defined in \eqref{eq: sphere model} and \eqref{eq: precessed sphere model}. The following Theorem considers the case of a crystal subject to uniform affine deformation.

\begin{theorem}\label{thm: single strain crystal}
	If $$u'(\vec x) = u(\vec x + \vecf R(\vec x)) = u(A\vec x+\vec b) \text{ for some } A\in\F R^{3\times3},\ \vec b\in\F R^3,\ u\in L^2(\F R^{\edit{n}{3}};\F C),$$
	where $A$ is invertible, then we can express its Fourier transform as:
	$$\C F[u'](\vec K) = \op{det}(A)^{-1}e^{\i\IP{\vec b}{A^{-\top}\vec K}}\C F[u](A^{-\top}\vec K).$$
\end{theorem}
	
This is a standard result in Fourier analysis and a proof is given in \ref{app: single strain crystal}. The appearance of the determinant emphasises that larger potentials diffract more strongly. The only changes to the Fourier transform are a corresponding linear deformation and a change of phase depending on the translation. To translate this result to the ideal crystal $u_0$, where the Fourier transform is a distribution rather than a smooth function, we need a further technical lemma.

\begin{lemma}\label{thm: single strain infinite crystal}
    If $A\in\F R^{3\times3}$ is an invertible matrix and $\vec b\in\F R^3$\edit{}{,} then
	$$\C F[u_0(A\cdot + \vec b)](\vec K) = e^{\i\IP{\vec b}{A^{-\top}\vec K}}\sum_{i=1}^\infty w_i\delta_{A^\top\vec p_i}(\vec K).$$
\end{lemma}

Due to its importance of this standard result in this work we include a proof in \ref{app: single strain crystal}. This Lemma demonstrates a key feature for diffraction from deformed single crystals. Under linear deformation, the location of spikes in Fourier space is equivalently linearly deformed and therefore the locations of diffracted peaks are also linearly deformed.

\subsection{Strained diffraction patterns}

To model the diffraction patterns produced by crystals subject to more general deformation we consider a crystal subject to piece-wise affine deformation. This is an implicit assumption that deformed crystals are also piece-wise smooth and so can be well approximated in this framework. We thus assume the material, $u$, may be expressed as:
\begin{equation}
    u(\vec x) = \sum_{j=1}^N u_0(A_j\vec x+\vec b_j)\F V\left(\frac{\vec x-\vec \beta_j}{\rho}\right) \label{eq: strained crystal def}
\end{equation}
where:
\begin{align*}
j\in[N]=\{1,\ldots, N\} & \qquad\text{ indexes over all, finitely many, voxels},
\\ \vec \beta_j\in\F R^{3} &\qquad \text{ is the location of voxel }  j,
\\ A_j\in\F R^{3\times 3} &\qquad \text{ is $A_j = \op{id}+\nabla\vecf F$ within voxel } j,
\\ \vec b_j\in\F R^{3} &\qquad \text{ represents the  shift to align }u_0\text{ with the $u$ in voxel } j.
\end{align*}
In this work, each voxel is a volume element that is sufficiently large that the discrete atoms blur into a continuum. We note that the results which follow can be extended to more generic tensor fields by limiting the voxel size to zero and replacing the finite sum with a Riemann integral if additional smoothness assumptions are made on $u_0$, $\vec b_j,$ and $A_j$ to guarantee any necessary exchanges of limits.
\\We can now express the full diffraction pattern of a deformed crystal.
\begin{lemma}\label{thm: full strain model}
	If the probe is narrow (\ref{ass: local signal}) then
	\begin{align*}
	D(\vec k) &= |\C F[\Psi_p]\star \C F[u]|^2(\vec k, k_z(\vec k))
	\\&= \left|\sum_{\substack{i\in\F N,j\in[N]\\\gamma^\top\vec\beta_j=\vec 0}}\hat w_{i}(k_z(\vec k) - (A_j^\top\vec p_i)_z, \beta_{j,z}) e^{\i\IP{\vec b_j}{\vec p_i}} f(\vec k - \gamma^\top A_j^\top\vec p_i)\right|^2
	\end{align*}
	where $$\hat w_{i}(k, \beta) = w_i\rho\op{sinc}(\rho k) e^{-\i \beta k}.$$
\end{lemma}
The proof of this is in \ref{app: approximations}. A key point here is that \ref{ass: local signal} is used to guarantee that the beam is completely contained within a single column of voxels. Without loss of generality, this column is centred at $x=y=0$ which reduces the sum to indices $j$ such that $\gamma^\top\vec\beta_j=\vec 0$.

In the limit $\rho\to\infty, \lambda\to0$ (large voxels, high energy incident electrons) Lemma~\ref{thm: full strain model} simplifies to
\begin{equation}\label{eq: high energy equation}
\lim_{\substack{\lambda\to0, \\\rho\to\infty}}\frac{D(\vec k)}{\rho} = \left|\sum_{i,j\in I}w_i e^{\i\IP{\vec b_j}{\vec p_i}}\C F[\Psi_p](\vec k-\gamma^\top A_j^\top\vec p_i)\right|^2
\end{equation}
where $I = \{(i,j)\in\F N^2 : \gamma^\top\vec\beta_j=\vec0,\ (A_j^\top\vec p_i)_z = 0\}$. This helps to highlight two key properties of diffraction imaging under deformation:
\begin{itemize}
	\item `in-plane' deformation moves the centre of the spot linearly from $(p_{i,x},p_{i,y})$ to $\gamma^\top A_j^\top\vec p_i = ((A_j^\top\vec p_i)_x, (A_j^\top\vec p_i)_y)$. 
	\item `out-of-plane' deformation changes the intensity of each spot depending on $(A_j^\top\vec p_i)_z$. In the high-energy limit this is absolute, however, in practice the intensity of spot $i$ is dampened by a factor of $\op{sinc}(\rho (k_z(\vec k)-A_j^\top\vec p_i)_z)\leq 1$. In either case, dependence of the intensities on the out-of-plane strain $A_j$ is highly non-linear.
\end{itemize}

\bigbreak Lemma~\ref{thm: full strain model} provides a very explicit model for computing diffraction patterns from deformed crystals, yet it is still too complex to directly derive a linear correspondence for the strain mapping inverse problem. To do this we will use precession and also one final technical assumption, that we are in a `small strain' scenario. Formally, we state this as:
\begin{enumerate}[leftmargin=\widthof{Asumption 1:+\quad},label=\textbf{Assumption \arabic*:},ref=Assumption \arabic*]
	 \setcounter{enumi}{6}
	\item \emph{small strain}, i.e. $|A_j-\op{id}| <\sigma \text{ for each } j$. \label{ass: small strain}
\end{enumerate}

Informally, we need to ensure that diffraction patterns of strained crystals still look like blurred diffraction patterns of single crystals. That is to say that the diffraction patterns should still be essentially sharp with isolated, if blurred, Bragg disks.

\section{Linearised model of electron diffraction from deformed crystals}\label{sec: approximations}

A linear tomography model is developed from the kinematical diffraction model in this section following a physically motivated argument based on precession electron diffraction, which is supported by a parallel mathematically rigorous argument in \ref{app: approximations}. It is not clear how the necessary assumptions may be justified physically, and we therefore provide computational results in Section~\ref{sec: numerics} which demonstrate that our model can be quantitatively accurate.

\subsection{Strained diffraction patterns with precession}

In Section~\ref{sec: precession intro}, precession was motivated as a technique for simplifying the computation of the average deformation. Now, we make this statement more precise by using tools from probability theory. \edit{}{For brevity of notation, we will use the abbreviation $$ \E_{\gamma^\top\vec \beta_j=0}X(j) \coloneqq \frac{\sum\left\{X(j) \text{ s.t. } j\in[N],\ \gamma^\top\vec \beta_j=0\right\}}{|\{j\text{ s.t. } j\in[N],\ \gamma^\top\vec \beta_j=0\}|}$$ for all random variables $X\colon[N]\to\F C$. In probability theory, this is the standard conditional expectation with respect to the condition $\gamma^\top\vec \beta_j=0$.}

\begin{approximation}\label{thm: precession simplified model}
	If the beam energy is large and the strain is small (\ref{ass: high energy} and \ref{ass: small strain}), then
	\begin{align}
	D_\alpha(\vec k) &= \bar D_\alpha(\vec k) + \text{error} \label{eq: simple strain}
	\\\text{where }\quad \bar D_\alpha(\vec k) &\coloneqq \left|\sum_{i=1}^\infty\bar w_i\E_{\gamma^\top\vec \beta_j=0}f(\vec k - \gamma^\top A_j^\top\vec p_i)\right|^2 \label{eq: dbar def}
	\end{align}
	for some new weights $\bar w_i\in\F C$.
\end{approximation}

$\bar D_\alpha$ is exactly the simple model one would hope for in strain mapping. Ignoring the squared norm, such a diffraction pattern is the average of each idealised diffraction pattern with an average spot-shape and weighted with an average structure-factor $\bar w_i$ independent of the strain parameter $A_j$. The cost for assuming that the raw data $D_\alpha$ obeys such a simple model is determined by the error term. With too little precession $D_\alpha$ will look close to Figure~\ref{fig: first examples}.a with inhomogeneous spot intensities and so the error will be large, however, too much precession introduces a new blurring of the data which also contributes to the error. The important point is that the error should not bias the computations that we go on to make in Approximation~\ref{thm: final model}, in particular it should not bias the centres of spots away from $\bar D_\alpha$.
	
Algebraically, it is hard to quantify the precise sweet-spot but the proof motivates a rule-of-thumb for the choice of precession angle. It should be just sufficiently large to ensure that 
$$k_z(\vec k) - (R_t^\top A_j^\top\vec p_i)_z = 0 \qquad \text{ for some } t\in[0,2\pi), |\vec k| = |\vec p_i| $$
for all \edit{$j$ and $i$ such that $p_{i,z}=0$}{relevant Bragg peaks $\vec p_i$ and strains $A_j$}. \edit{This is sufficient to guarantee that for any spot visible in Figure~\ref{fig: first examples} and all small strains ($|A_j-\op{id}|<\sigma$), the precession angle $\alpha$ is large enough such that $R_t^\top A_j^\top \vec p_i$ lies on the Ewald sphere for some $t$.}{i.e. we require $\alpha$ to be sufficiently large such that the strained Bragg peak $A_j^\top\vec p_i$ exactly intersects the Ewald sphere for some $t$. For example, all visible peaks in Figure~\ref{fig: first examples} satisfy $p_{i,z}=0$, $|\vec p_i|\leq P$ for some $P>0$, and we assume small strain ($|A_j-\op{id}|<\sigma$).} \edit{After a geometrical argument detailed in \ref{app: precession rule-of-thumb}, to analyse the deformation of spots $\vec p_i$ with $|\vec p_i|<P$, we suggest the relation}{A short geometrical argument given in \ref{app: precession rule-of-thumb} suggests the relation}
\begin{equation}\label{eq: precession angle}
    \alpha \approx \cos^{-1}\left(1-\frac{\sigma^2}{2}\right) + \sin^{-1}\left(\frac{\lambda P}{4\pi}\right)\sim \sigma + \frac{\lambda P}{4\pi}.
\end{equation}
In words, the precession angle should be larger than the maximum rotation due to the deformation plus the distance between the flat hyperplane and the Ewald sphere.
	
One final observation from this approximation is that the coordinate $(A_j^\top \vec p_i)_z$ has disappeared completely, other than in the choice of $\alpha$ above, and the remaining expression only depends on $\gamma^\top A_j^\top \vec p_i$. This indicates that diffraction imaging is insensitive to deformations parallel to the beam direction, which will dictate our choice of tensor tomography model in Section~\ref{sec: Tensor Tomography}.

\subsection{Linearised diffraction model}

The final approximation is to linearise the forward model. $\bar D_\alpha$ is now a sufficiently simple, but still non-linear, model to go from a deformed crystal to a diffraction image, however, what we want is a simple linear model to map from a deformed crystal to its average deformation tensor. To do this we propose a simple pre-processing procedure to apply to the raw data $D_\alpha$ which corresponds to a linear forward model with respect to the deformation parameters. In particular, we propose computing centres of mass for each of the observed diffraction spots.
	
\begin{approximation}\label{thm: final model}
	If the conditions of Approximation~\ref{thm: precession simplified model} hold and the diffracted spots are symmetric and non-overlapping (\ref{ass: separated spots} and \ref{ass: symmetric}) then
	\begin{equation}\label{eq: com approx} 
	\E_{\gamma^\top \beta_j = 0}\gamma^\top A_j^\top \vec p_i = \frac{\int_{|\vec k-\vec q_i|<\bar r} \vec k D_\alpha(\vec k)\D\vec k}{\int_{|\vec k-\vec q_i|<\bar r} D_\alpha(\vec k)\D\vec k} + \text{error}
	\end{equation}
	where $\bar r>0$ is the separation of spots given by \ref{ass: separated spots}.
\end{approximation}

This \edit{theorem}{Approximation} directly indicates that centres of mass are a good linear model for deformed diffraction patterns and this is confirmed by the numerical results in Section~\ref{sec: numerics}. More generally, this provides a motivation that the centre of deformed spots are equal to average deformation tensors. Other centre detection methods are common in the literature, and we also give numerical comparison of the accuracy and robustness of each method.

\section{Non-symmetric Tensor Tomography}\label{sec: Tensor Tomography}
Up to this point we have considered how to compute a single average deformation tensor from a single diffraction pattern. The role of this section is to identify this process with an inverse problem capable of reconstructing a 3D strain map from many average deformation tensors, and then highlight the relevant properties of this inverse problem. 

The fact that strain maps are rank 2 tensor fields immediately puts us in the domain of tensor tomography, and the physical characteristics of diffraction imaging seen in Approximation~\ref{thm: precession simplified model} highlights the \emph{transverse ray transform} (TRT) as the natural parallel. In particular, Approximation~\ref{thm: precession simplified model} shows that precessed diffraction patterns are insensitive to out-of-plane strain. This aligns with equivalent reasoning in \cite{Lionheart2015} for polycrystalline materials where it was shown that this corresponds to the TRT. The only difference for polycrystalline materials is that the forward model is insensitive to the skew component of deformations, and so we will extend the analysis of the TRT to account for general tensor fields of non-symmetric tensors.

\subsection{Notation and definitions}
We define that:
\begin{itemize}
	\item $\vecf f\colon \F R^3\to \F R^{3}$ denotes a \edit{}{vector field, i.e. a} rank 1 tensor field.
	\item $\vecf F\colon \F R^3\to \F R^{3\times 3}$ denotes a rank 2 tensor field.
	\item $\F S^2 = \{\vec \xi\in\F R^3 \text{ such that } |\vec \xi| = 1\} $ is the unit sphere in $\F R^3$.
	\item $\C S$ denotes a Schwartz space.
	\item The set of all lines is identified with the set
		$$ T\F S^2 = \{(\vec \xi, \vec x)\in \F R^3\times\F R^3 \text{ such that } |\vec \xi| = 1, \IP{\vec x}{\vec \xi} = 0\} \subset \F S^2\times\F R^3$$
		under the understanding that $(\vec\xi, \vec x)\mapsto \{\vec x + t\vec \xi \text{ s.t. }t\in\F R\}$ parametrises all lines in $\F R^3$ (uniquely with orientation).
	\item $\op{Sym}(\F R^{3\times 3}) = \{ A\in\F R^{3\times 3} \text{ such that } A^\top  = A\}$ is the set of symmetric matrices.
	\item $\op{Sym}(\vecf F)(\vec x) = \frac12(\vecf F(\vec x) + (\vecf F(\vec x))^\top )$ is the pointwise symmetrisation of $\vecf F$.
	\item $\cross{\vec x}$ is the \edit{}{skew-symmetric} matrix such that $\cross{\vec x}\vec y = \vec x\times\vec y$ for all $\vec y\in\F R^3$. Similarly, $\cross{\vecf f}(\vec x) = \cross{\vecf f(\vec x)}$ is the pointwise operation.
	\item $\vec \xi^\perp = \{ \vec x\in\F R^3 \text{ such that } \IP{\vec x}{\vec \xi} = \vec 0\}$ is the orthogonal complement.
	\item $\op{id}= \left(\begin{smallmatrix} 1&0&0\\0&1&0\\0&0&1\end{smallmatrix}\right)$ is the identity matrix.
\end{itemize}

We follow the definitions of tensor tomography operators given in \cite{Sharafutdinov1994}, defining the \emph{longitudinal ray transform} (LRT), $I\colon \C S(\F R^3;\F R^3)\to \C S(T\F S^2;\F R)$
$$ I\vecf f(\vec \xi,\vec x) = \int_{\F R} \IP{\vecf f(\vec x+t\vec \xi)}{\vec \xi} \D t$$
and the \emph{transverse ray transform} (TRT), $J\colon \C S(\F R^3;\F R^{3\times3})\to \C S(T\F S^2;\F R^{3\times 3})$
\begin{equation}\label{eq: TRT def}
J\vecf F(\vec \xi,\vec x) = \int_{\F R} \Pi_{\vec \xi}^\top  \vecf F(\vec x+t\vec \xi)\Pi_{\vec \xi} \D t
\end{equation}
where $\Pi_{\vec \xi} = \op{id} - \frac{\vec \xi\vec \xi^\top }{|\xi|^2}$ is the orthogonal projection map such that $\IP{(\Pi_{\vec \xi} \vec x)}{\vec \xi} = 0$ for all $\vec x\in\F R^3$.  

\subsection{Electron Diffraction and the Transverse Ray Transform}

The first step is to generalise the notation from Section~\ref{thm: final model} to consider diffraction patterns where the electron beam is not parallel to the $z$-axis or through the point $x=y=0$. This is mainly an algebraic exercise, first to upgrade from average spot centres to average strain tensors, and then to realise the generalisation.

\begin{lemma}\label{thm: noncolinear spots}\hfill
	\begin{enumerate}
		\item \emph{Tensors from vectors:} $\E_{\gamma^\top \vec \beta_j=0} \gamma^\top A_j^\top \gamma$ can be computed from the values of $\E_{\gamma^\top \vec \beta_j=0} \gamma^\top A_j^\top \vec p_i$ for any two non-colinear points, say $\vec p_1, \vec p_2$, such that $p_{i,z}=0$.
		\item \emph{Generalising notation: } If we choose $\vec x = \left(\begin{smallmatrix}
		0\\0\\0\end{smallmatrix}\right)$ and $\Pi_\xi = \op{id}-\frac{\vec\xi\vec\xi^\top }{|\vec\xi|^2}$ with $\vec\xi=\left(\begin{smallmatrix}
		0\\0\\1\end{smallmatrix}\right)$, then
		$$ {\operatorname*{\E}_{\gamma^\top \vec \beta_j=\vec x}} \Pi_\xi^\top A_j^\top \Pi_\xi = \begin{pmatrix}
		\operatorname*{\E}\limits_{\gamma^\top \vec \beta_j=\vec0} \gamma^\top A_j^\top \gamma & \begin{matrix}0\\0\end{matrix}
		\\ \begin{matrix}\quad0&\quad0\end{matrix} & 0
		\end{pmatrix}. $$
	\end{enumerate}
\end{lemma}
\begin{proof}
	Let square brackets temporarily denote the horizontal concatenation of vectors into matrices. Note that 
	$$[\vec p_1,\vec p_2] = \left(\begin{smallmatrix} p_{1,x}&p_{2,x}\\p_{1,y}&p_{2,y}\\0&0 \end{smallmatrix}\right) = \left(\begin{smallmatrix} 1&0\\0&1\\0&0 \end{smallmatrix}\right)\left(\begin{smallmatrix} p_{1,x}&p_{2,x}\\p_{1,y}&p_{2,y}\end{smallmatrix}\right) = \gamma\left(\begin{smallmatrix}	p_{1,x}&p_{2,x}\\p_{1,y}&p_{2,y}\end{smallmatrix}\right).$$
	$\vec p_i$ are non-colinear and so this final matrix is invertible. Thus, we have
	\begin{align*}
	\E_{\gamma^\top \vec \beta_j=0} \gamma^\top A_j^\top \gamma &= \E_{\gamma^\top \vec \beta_j=0} \left\{\gamma^\top A_j^\top \gamma \left(\begin{smallmatrix}		p_{1,x}&p_{2,x}\\p_{1,y}&p_{2,y}\end{smallmatrix}\right)\left(\begin{smallmatrix}		p_{1,x}&p_{2,x}\\p_{1,y}&p_{2,y}\end{smallmatrix}\right)^{-1}\right\}
	\\&= \E_{\gamma^\top \vec \beta_j=0} \left\{\gamma^\top A_j^\top [\vec p_1,\vec p_2]\right\}\left(\begin{smallmatrix}	p_{1,x}&p_{2,x}\\p_{1,y}&p_{2,y}\end{smallmatrix}\right)^{-1}
	\\&= \left[\E_{\gamma^\top \vec \beta_j=0} \gamma^\top A_j^\top \vec p_1, \E_{\gamma^\top \vec \beta_j=0} \gamma^\top A_j^\top \vec p_2\right]\left(\begin{smallmatrix}		p_{1,x}&p_{2,x}\\p_{1,y}&p_{2,y}\end{smallmatrix}\right)^{-1}
	\end{align*}
	which verifies the first part. The final part is a simple algebraic argument:
	$$\Pi_\xi = \op{id} - \left(\begin{smallmatrix} 0\\0\\1\end{smallmatrix}\right)\left(\begin{smallmatrix} 0&0&1\end{smallmatrix}\right) = \left(\begin{smallmatrix} 1&0&0\\0&1&0\\0&0&0\end{smallmatrix}\right) = \left(\begin{smallmatrix} 1&0\\0&1\\0&0\end{smallmatrix}\right)\left(\begin{smallmatrix} 1&0&0\\0&1&0\end{smallmatrix}\right) = \gamma\gamma^\top .$$
	From this we see $$\E_{\gamma^\top \vec \beta_j=\vec x} \Pi_\xi^\top A_j^\top \Pi_\xi = \gamma\left\{\E_{\gamma^\top \vec \beta_j=0}\gamma^\top A_j\gamma\right\}\gamma^\top  = \begin{pmatrix}
	\operatorname*{\E}\limits_{\gamma^\top \vec \beta_j=0} \gamma^\top A_j^\top \gamma & \begin{matrix}0\\0\end{matrix}
	\\ \begin{matrix}\quad0&\quad0\end{matrix} & 0
	\end{pmatrix}$$
	as required.
\end{proof}

This Lemma is where \ref{ass: thick crystals} becomes relevant. If the crystal is of finite thickness then all $\vec p_i$ are always present in the diffraction pattern (not just when $p_{i,z}=0$), although the intensity may be small. In practice, all visible spots will lie on this hyperplane which is justified by \ref{ass: thick crystals}.

The final step to aligning with the TRT is to replace discrete sums with line integrals. We compute
\begin{align}
\edit{}{|\{j\;\op{s.t.}\;\gamma^\top \vec\beta_j=\vec 0\}|}\E_{\gamma^\top \vec\beta_j=\vec x}\Pi_{\vec \xi}^\top A_j^\top \Pi_{\vec \xi} &= \edit{\frac{\int_{\gamma^\top \vec\beta_j=0} \Pi_{\vec \xi}^\top A(\vec \beta_j)^\top \Pi_{\vec \xi} \D j}{|\{j\;\op{s.t.}\;\gamma^\top \vec\beta_j=\vec 0\}|}}{\int_{\gamma^\top \vec\beta_j=0} \Pi_{\vec \xi}^\top A(\vec \beta_j)^\top \Pi_{\vec \xi} \D j} &
\\&\edit{\propto}{=} \int_{-\infty}^\infty \Pi_{\vec \xi}^\top A(\vec x+t\vec \xi)^\top\Pi_{\vec \xi} \D t. \label{eq: com to TRT}
\end{align}
This is now in the classical TRT format as in \eqref{eq: TRT def}. \edit{}{The additional scaling on the left-hand side is an important difference between diffraction imaging and the TRT although we postpone discussion until Section~\ref{sec: physics vs mathematical TRT}.}

\subsection{Invertibility of tensor ray transforms}

The results that we show in this section are necessary for the analysis of the strain tomography inverse problem but are not novel in the field of tensor tomography. Instead, the intention is to bring together all \edit{of }{}the relevant results to clearly explain the relatively special case of three dimensions. In dimension strictly greater than three the TRT is invertible and its properties are already well explored \cite{Novikov2007,Abhishek2020}. In a general three dimensional Riemannian manifold the analysis is given by \cite{Holman2013} however this is much more complicated than needed in the Euclidean case. Our argument closely follows that of \cite[Section 4]{Novikov2007}, however we also retain the symmetric component of the tensor field.

We note that we have formulated these ray transforms in a Schwartz space setting in order to make use of existing results regarding invertibility of both the LRT and the TRT\edit{ in the next section}{}. Recent work has extended the stability results of the LRT to Sobolev spaces on compact domains \cite{Boman2018} and the same techniques appear valid for the TRT but this has not yet been made rigorous.

The invertibility of the LRT and the TRT for symmetric rank 2 tensor fields is established in \cite[Chapter 2]{Sharafutdinov1994} \edit{}{(LRT), }and \cite[Chapter 5, page 155]{Sharafutdinov1994} \edit{respectively}{(TRT)}. We state the key results:

\begin{theorem}\label{thm: old invertibility}
For all $\vecf f,\vecf g \in \C S(\F R^3;\F R^{3})$, $\vecf F,\vecf G \in \C S(\F R^3;\op{Sym}(\F R^{3\times 3}))$, orthogonal bases $\{\vec e_1,\vec e_2,\vec e_3\}$
\begin{align}
I\vecf f(\vec\xi,\vec x) = I\vecf g(\vec \xi,\vec x) \text{ for all } (\vec \xi,\vec x)\in T\F S^2 \hspace{7pt}\not\hspace{-7pt}\implies &\vecf f=\vecf g
\\ J\vecf F(\vec \xi,\vec x) = J\vecf G(\vec \xi,\vec x) \text{ for all } \vec x, \;\vec \xi\in \cup_i \vec e_i^\perp \iff &\vecf F=\vecf G
\end{align}
In particular,
$$ I\vecf f = I\vecf g \iff \exists \varphi\in\C S(\F R^3;\F R) \text{ such that } \vecf f = \vecf g + \nabla\varphi.$$
\end{theorem}

In other words, vector fields $\vecf f$ and $\vecf g$ can never be completely distinguished\edit{ experimentally whereas}{, only the solenoidal component can be reconstructed from LRT-type data. On the other hand,} any two symmetric tensor fields, $\vecf F$ and $\vecf G$, can be distinguished from experimental data \edit{}{of TRT-type }from three well-chosen tilt-series.
\\We can now state the invertibility of the TRT for non-symmetric rank 2 tensor fields. 
\begin{theorem}\label{thm: new invertibility}
	For all $\vecf F,\vecf G \in \C S(\F R^3;\F R^{3\times 3})$, orthogonal bases $\{\vec e_1,\vec e_2,\vec e_3\}$
	\begin{align}
	J\vecf F(\vec \xi,\vec x) = J\vecf G(\vec \xi,\vec x) \text{ for all } (\vec \xi,\vec x)\in T\F S^2 \iff &\vecf F=\vecf G + \cross{\nabla\varphi}, \text{ some }\varphi\in \C S(\F R^3;\F R)
	\\ J\vecf F(\vec \xi,\vec x) = J\vecf G(\vec \xi,\vec x) \text{ for all } \vec x, \;\vec \xi\in \cup_i \vec e_i^\perp \iff &\op{Sym}(\vecf F)=\op{Sym}(\vecf G)
	\end{align}
\end{theorem}
The proof of Theorem~\ref{thm: new invertibility} relies on the following decomposition lemma.

\begin{lemma}\label{thm: decomp lemma}
\edit{The TRT of a general tensor field can be decomposed as:}{Decomposition of the TRT of non-symmetric tensors:}
	\begin{enumerate}
		\item For all $\vecf F \in \C S(\F R^3;\F R^{3\times 3})$ there exists $\vecf f\in \C S(\F R^3;\F R^{3})$ such that $$\vecf F = \op{Sym}(\vecf F) + \cross{\vecf f}$$
        \item \edit{}{For all $\vecf F \in \C S(\F R^3;\F R^{3\times 3})$,} $J\op{Sym}(\vecf F)(\vec \xi,\vec x) = \op{Sym}(J\vecf F(\vec \xi,\vec x))$
		\item \edit{}{For all $\vecf f\in \C S(\F R^3;\F R^{3})$,} $J\cross{\vecf f}(\vec \xi,\vec x) = I\vecf f(\vec \xi,\vec x) \cross{\vec \xi}$
\end{enumerate}
\end{lemma}

\begin{proof}[Proof of Lemma~\ref{thm: decomp lemma}.]
	Part \textit{(i)} is a simple algebraic equivalence. For any $A\in\F R^{3\times3}$ 
	\fitwidtheq{$$A - \op{Sym}(A) = \frac12\begin{pmatrix}
	0 & A_{1,2}-A_{2,1} & A_{1,3}-A_{3,1}
	\\ -(A_{1,2}-A_{2,1}) & 0 & A_{2,3}-A_{3,2}
	\\ -(A_{1,3}-A_{3,1}) & -(A_{2,3}-A_{3,2}) & 0
	\end{pmatrix} = \cross{\frac12\begin{pmatrix}A_{3,2}-A_{2,3} \\ A_{1,3}-A_{3,1} \\ A_{2,1}-A_{1,2} \end{pmatrix}}$$}
	The expression $\vecf F = \op{Sym}(\vecf F) + \cross{\vecf f}$ is simply the pointwise extension of this equality. Confirming $\vecf f\in \C S$ is also clear as $\cross{\vecf f} = \frac12\vecf F-\frac12\vecf F^\top \in \C S$ and so we must have $\vecf f_i\in\C S$ for each $i=1,\ldots,3$. By the linearity of $J$
	$$J\op{Sym}(\vecf F) = J\left(\tfrac12(\vecf F + \vecf F^\top )\right) = \tfrac12\left(J\vecf F + J(\vecf F^\top )\right).$$
	Hence, to prove part \textit{(ii)}, it suffices to show:
	$$J(\vecf F^\top )(\vec \xi,\vec x) = \int_{\F R} \Pi_{\vec \xi}\vecf F(\vec x+t\vec \xi)^\top \Pi_{\vec \xi} \D t = \int_{\F R} (\Pi_{\vec \xi}\vecf F(\vec x+t\vec \xi)\Pi_{\vec \xi})^\top  \D t = (J\vecf F(\vec \xi,\vec x))^\top .$$
    \bigbreak\noindent The proof of part \textit{(iii)} is also by direct evaluation, fix $\vec a\in\F R^3$. We claim $\Pi_{\vec\xi}\cross{\vec a}\Pi_{\vec\xi}=\IP{\vec a}{\vec \xi}\cross{\vec\xi}$. First note the trivial case, when $\vec\xi = \vec e_3=(\begin{smallmatrix}0&0&1\end{smallmatrix})^\top $
    $$ \Pi_{\vec e_3}\cross{\vec a}\Pi_{\vec e_3} = \left(\begin{smallmatrix}1&0&0\\0&1&0\\0&0&0\end{smallmatrix}\right) \left(\begin{smallmatrix}0&-a_3&a_2\\a_3&0&-a_1\\-a_2&a_1&0\end{smallmatrix}\right) \left(\begin{smallmatrix}1&0&0\\0&1&0\\0&0&0\end{smallmatrix}\right) = \left(\begin{smallmatrix}0&-a_3&0\\a_3&0&0\\0&0&0\end{smallmatrix}\right) = a_3\cross{\vec e_3}.$$
    To generalise this, recall the property of cross products $ (R\vec a)\times (R\vec b) = R(\vec a\times \vec b)$ for all rotation matrices $R$ and vectors $\vec b$. Because of this, we know 
    $$R\cross{\vec a}R^\top \vec b = R(\vec a\times R^\top \vec b) = (R\vec a)\times \vec b = \cross{R\vec a}\vec b, \qquad R\Pi_{\vec\xi}R^\top  = R\left(\op{id}-\frac{\vec\xi\vec\xi^\top }{|\xi|^2}\right)R^\top  = \Pi_{R\vec\xi}.$$
    If we choose $R$ such that $\vec\xi=R\vec e_3$, it follows
    \begin{multline*}
        \Pi_{\vec\xi}\cross{\vec a}\Pi_{\vec \xi} = \Pi_{R\vec e_3}\cross{\vec a}\Pi_{R\vec e_3} = R\Pi_{\vec e_3} R^\top \cross{\vec a}R\Pi_{\vec e_3}R^\top  = R(\Pi_{\vec e_3} \cross{R^\top \vec a}\Pi_{\vec e_3}) R^\top 
        \\= R(\IP{R^\top \vec a}{\vec e_3} \cross{\vec e_3})R^\top  = \IP{\vec a}{\vec\xi}\cross{\vec \xi}
    \end{multline*}
    as required. Finally, we can extend this equality pointwise to $\vecf f$
	$$ J\cross{\vecf f}(\vec \xi, \vec x) = \int_{\F R} \Pi_{\vec \xi}\cross{\vecf f(\vec x+t\vec \xi)}\Pi_{\vec \xi} \D t = \int_{\F R} \IP{\vecf f(\vec x+t\vec \xi)}{\vec \xi}\cross{\vec \xi} \D t = I\vecf f(\vec \xi,\vec x)\cross{\vec \xi}.$$
\end{proof}
\begin{proof}[Proof of Theorem~\ref{thm: new invertibility}.]
	As $J$ is a linear map, it suffices to prove the theorem in the case $\vecf G=0$. Applying the decomposition established in Lemma~\ref{thm: decomp lemma}:
	\begin{align*}
	J\vecf F(\vec \xi,\vec x) = 0 &\iff J\op{Sym}(\vecf F)(\vec \xi,\vec x) + J\cross{\vecf f}(\vec \xi,\vec x) = 0 &\text{Lemma~\ref{thm: decomp lemma}(i)}
	\\&\iff J\op{Sym}(\vecf F)(\vec \xi,\vec x) + I\vecf f(\vec \xi,\vec x)\cross{\vec \xi} = 0 &\text{Lemma~\ref{thm: decomp lemma}(iii)}
	\\&\iff J\op{Sym}(\vecf F)(\vec \xi,\vec x) = \op{Sym}(0) & \text{Lemma~\ref{thm: decomp lemma}(ii)+symmetric/}
	\\ &\phantom{\iff}\;\text{ and }I\vecf f(\vec \xi,\vec x)\cross{\vec \xi} = 0-\op{Sym}(0)& \text{skew decomposition}
	\\&\iff J\op{Sym}(\vecf F)(\vec \xi,\vec x) = I\vecf f(\vec \xi,\vec x) = 0 & \vec \xi\neq\vec 0
	\end{align*}
	Hence, by Theorem~\ref{thm: old invertibility} we know that $\vecf f$, and thus the skew component of $\vecf F$, is never uniquely determined but the symmetric component can be recovered through the equality of Lemma~\ref{thm: decomp lemma}(ii). 
\end{proof}
	
This decomposition, Lemma~\ref{thm: decomp lemma}, is very powerful for understanding the tensor tomography problem from an analytical stand-point. Theorem~\ref{thm: new invertibility} lifts invertibility results from \cite{Sharafutdinov1994, Desai2016} to demonstrates when (or how much) the non-symmetric TRT is invertible. On the other hand, results from \cite{Desai2016} also compute the (pseudo) inverses of the LRT and symmetric TRT which could be lifted to a filtered-back projection-like pseudo-inverse of the non-symmetric TRT. This in turn allows us to characterise the sensitivity to noise, the inverse problem of non-symmetric tomography is only mildly ill-posed as the singular values decay at only a polynomial rate.

\subsection{Physical setting of the transverse ray transform}\label{sec: physics vs mathematical TRT}

There are two key differences before applying the standard theory. Firstly, the scaling constant switching between average to integral was ignored in \eqref{eq: com to TRT}, however, practically this represents a necessary re-scaling by the thickness of the specimen to go from the raw data (average deformation) to the linear model (integral of deformation). This also appears in the analysis as a violation of the small strain assumption, \ref{ass: small strain}, because the vacuum outside \edit{of }{}the crystal has a deformation of order one. This scaling allows us to account for the violation by incorporating \edit{outside}{alternative} knowledge of the specimen. Experimentally, a high-angle annular dark-field (HAADF) STEM image can be recorded at each specimen orientation to record the object thickness\cite{Midgley2003}.

Secondly, we draw attention to the different acquisition geometries expected from the theoretical and physical settings. In Theorem~\ref{thm: new invertibility} the TRT requires three (well\edit{}{-}chosen) continuous tilt-series to be suitably invertible, however in Lemma~\ref{thm: noncolinear spots} we can only compute the signal if two non-colinear spots are visible in the diffraction pattern. This places a strict constraint on the choice of $\vec\xi$ and, in particular, the set of physically feasible beam orientations is a discrete set. Because of this, datasets in this application will always be in a limited data scenario.

\section{Computational Validation}\label{sec: numerics}

We perform a computational analysis to validate the approximation of the forward model  \eqref{eq: com approx} used to relate electron diffraction to the TRT and to confirm that the TRT inverse problem, which is under-determined and has a non-trivial null-space, can be solved accurately using a realistic amount of data. These tasks are separated for reasons of efficiency by first rigorously testing the forward model to provide a worse case estimate for the error level, and then simulate error at and above this level for numerical tomographic reconstructions.

\subsection{Forward model validation}

Physical validation requires quantifying the level of error in \eqref{eq: com approx} knowing the exact deformation and comparing against the computed deformation determined based on measuring the centres of the Bragg diffraction disks. This will be assessed in three ways. First we use the MULTEM \cite{Lobato15} package as a dynamical simulator, in particular modelling multiple scattering effects. We suggest this provides a worst-case analysis as the simulation model is more accurate than the analytical model used, and also we use thick crystals where the extra complexity should be most apparent. Second we compare against a kinematical simulation, the analytical model of this study, and observe that the accuracy is equivalent to that of the dynamical simulation. The previous two comparisons use crystals formed as in \eqref{eq: strained crystal def} which ignores strain on the boundaries of each block. The final comparison uses crystals with dislocations where the deformation is continuously defined to verify robustness to the specific definition of strain.

\subsubsection{Piece-wise affine deformation phantom}\label{sec: numerics modelling}
Phantom crystals are defined with piece-wise linear deformations as in \eqref{eq: strained crystal def} but using the rigid-ion model of deformation, i.e. where the \edit{array}{lattice} of atoms is deformed but atoms remain spherical and the same size. To do this, we create a three parameter collection of phantoms.

\begin{definition}
	Given an initial crystal structure, randomly sampled phantoms are built up atomistically using three parameters:
	\begin{enumerate}
		\item $L\in\F N$ is the number of layers. Phantoms consist of layers stacked orthogonal to the $z$-axis. Each layer is under constant affine transformation and is (approximately) the same thickness.
		\item $d\in\{1,2,3\}$ is the rank of the displacement gradient tensor. $d=1$ corresponds to a simple isotropic scaling of the initial crystal structure. $d=2$ also allows for rotation and shearing within the layer. $d=3$ allows any generic $3\times 3$ tensor.
		\item $\sigma\geq0$ is the average magnitude of atomic displacement. In particular, over each randomly generated displacement gradient tensor, the average of the spectral norm of the perturbation from identity will be equal to $\sigma$. Note that this average is not enforced in each simulated phantom.
	\end{enumerate}
\end{definition}

\subsubsection{Continuous deformation phantom}\label{sec: numerics modelling analytical}

\begin{figure}[ht]
	\centering
	\begin{subfigure}{.45\textwidth}\centering
		\includegraphics[width=.8\textwidth,trim={55 50 325 60},clip]{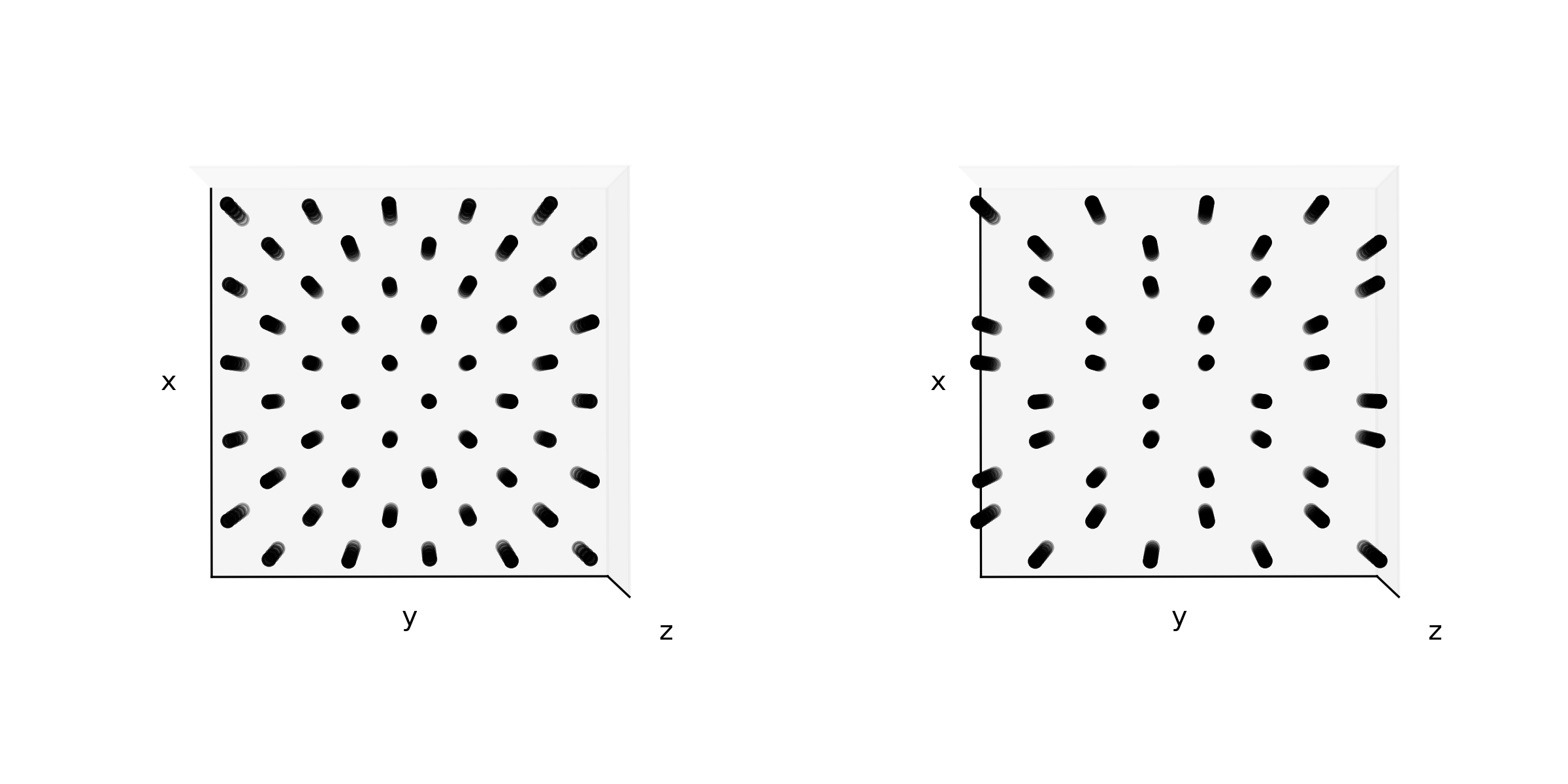}
		\caption{[001] crystal}
	\end{subfigure}
	\begin{subfigure}{.45\textwidth}\centering
		\includegraphics[width=.8\textwidth,trim={340 50 40 60},clip]{phantoms}
		\caption{[011] crystal}
	\end{subfigure}
	\begin{subfigure}{.45\textwidth}\centering
		\includegraphics[width=.8\textwidth,trim={35 45 300 70},clip]{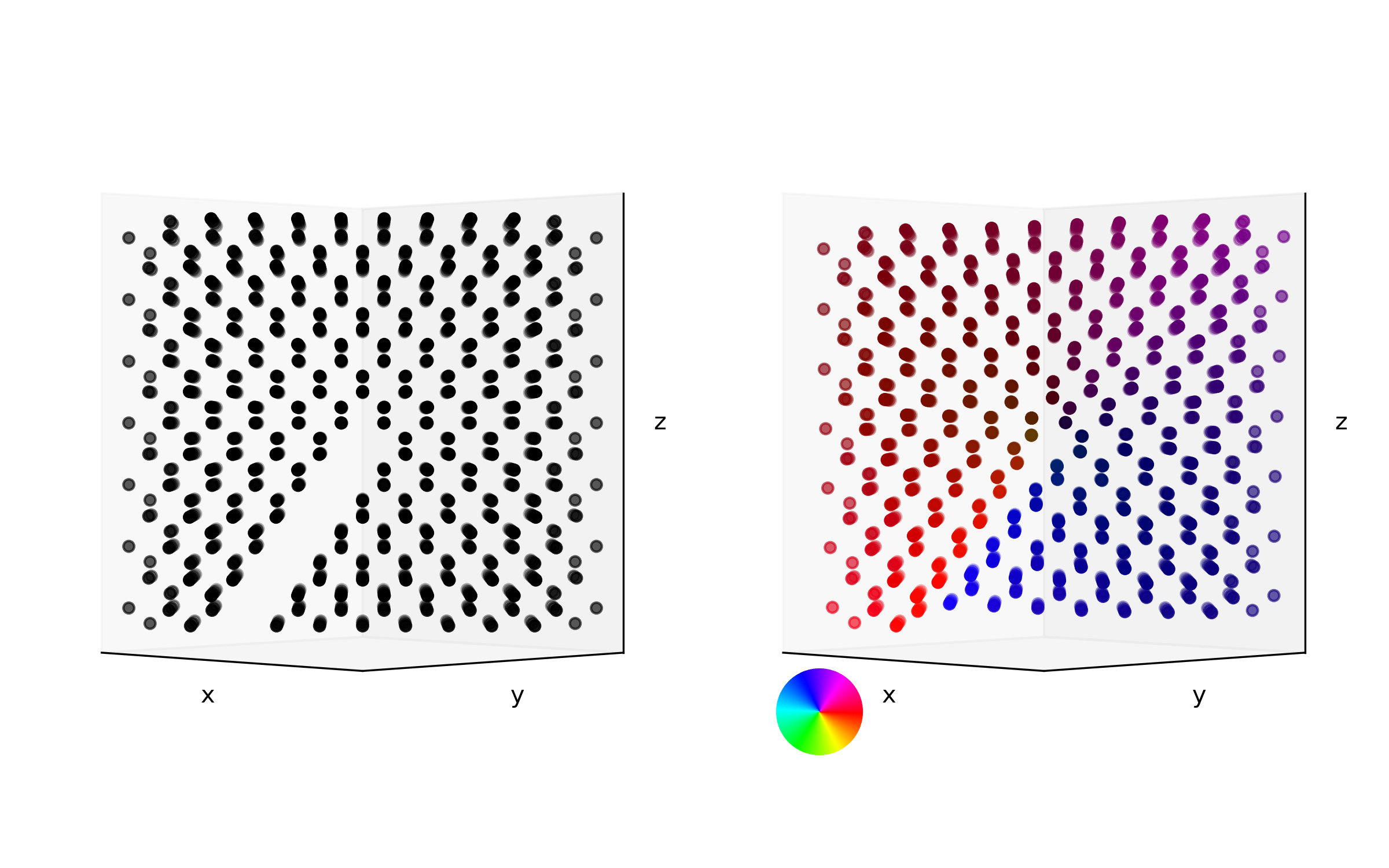}
		\caption{Unstrained dislocation}
	\end{subfigure}
	\begin{subfigure}{.45\textwidth}\centering
		\includegraphics[width=.8\textwidth,trim={320 45 15 70},clip]{dislocation_phantom}
		\caption{Strained dislocation}
	\end{subfigure}
	\caption{Parts (a) and (b) show columns of atoms parallel to the beam for two Silicon crystals. To form a dislocation, we remove a half-plane of atoms (c) atoms were displaced according to \eqref{eq: displacement field} (d). Deformation direction/magnitude is visualised by colour/brightness respectively.}
	\label{fig: phantom crystals}
\end{figure}

When defined atomistically, the deformation gradient tensor requires a choice of interpolation to be related our discussion in Section~\ref{sec: physical model with strain}. The piece-wise affine deformation phantoms test this in a discrete sense where deformations were defined to be piece-wise constant with discontinuities at the interfaces between layers. To test accuracy under continuous deformation we define use the continuum deformation associated with a dislocation with Burger's vector $\vec b \propto \left(1,1,1\right)^\top $ and line vector $\vec u \propto \left(1,-1,0\right)^\top $ \cite{hirsch1966}. In the notation of Definition~\ref{def: strain}
\begin{equation}
    \vecf R(\vec x) \coloneqq \frac{4(1-\nu)\phi+\sin(2\phi)}{8\pi(1-\nu)}\vec b + \frac{2(1-2\nu)\log(r) + \cos(2\phi)}{8\pi(1-\nu)}(\vec u\times\vec b) - \vec x \label{eq: displacement field}
\end{equation}
where $(r,\phi,Z)$ are the cylindrical coordinates of $\vec x$ on the basis $(\vec u\times\vec b, \vec u\times\vec b\times \vec u, \vec u)$. A dislocation can be modelled naively by removing all atoms along the half-plane defined by $\phi=\pi$ from a crystal and displacing the atoms according to \eqref{eq: displacement field}. This is visualised in Figure~\ref{fig: phantom crystals} and 100 diffraction patterns were simulated from this phantom with beam parallel to the $z$-axis and $\alpha=2$.
	
\subsubsection{Electron diffraction simulations}

Electron diffraction patters were simulated using a wavelength $\lambda = \SI{0.02}{\angstrom}$ (energy of ca. \SI{300}{\keV}) and an average deformation magnitude of $\sigma=0.01(=1\%)$ assuming an incident electron probe function of the form,
$$ \Psi_p(x,y) = \frac{J_1(r\sqrt{x^2+y^2})}{r\sqrt{x^2+y^2}}, \qquad \C F[\Psi_p](\vec k) = \begin{cases}\tfrac{1}{r^2} & |\vec k|<r\\ 0&\text{else}\end{cases}$$
where $J_1$ is a Bessel function of the first kind. This probe function corresponds to the probe produced using a circular disk aperture in ideal probe forming optics, and the choice of $r$ corresponds to an aperture of 2 mrad.

Kinematical simulations are computed using \eqref{eq: precessed sphere model} while dynamical simulations were performed using the multislice approach as implemented in the MULTEM package \cite{Lobato15}. The number of points chosen to simulate precession was chosen sufficiently large such that the errors reported in Table~\ref{tbl: quantitative results} had converged, details are given in Section~\SUPP{\ref{supp: precession discretisation}}{S.2}. Precession angles in the range \SIrange{0}{2}{\degree} were used.

\subsubsection{Results: Linearised model accuracy} \label{sec: numerical simulation results}

Computations presented in this section were performed to determine the accuracy with which Bragg disk centres can be measured and related to the average deformation tensor. We also considered how this accuracy may be affected by both the choice of precession angle and the choice of disk centre detection algorithm. \edit{The theory suggested in}{Equation} \eqref{eq: com approx} \edit{}{suggests }that the centre of mass is an accurate measure of average deformation, \edit{and }{}we compare this with disk-detection by cross-correlation, which is common in strain mapping literature, to evaluate which method is best at linearising the forward model.

Precession angle is an experimental parameter\edit{}{,} which \eqref{eq: precession angle} suggests should be chosen as 
$$\alpha \approx \cos^{-1}\left(1-\frac{\sigma^2}{2}\right) + \sin^{-1}\left(\frac{\lambda P}{4\pi}\right)\approx 0.01 + \frac{0.02\times 5}{4\pi} \approx \SI{0.6}{\degree} + \SI{0.5}{\degree} = \SI{1.1}{\degree}.$$
It is more common in practice to use $\alpha< \SI{1}{\degree}$, however in this case we wish to quantify any expected advantage of using larger values.
The following simulations were performed:
\begin{itemize}
	\item 72 MULTEM simulated diffraction patterns with one random phantom for each combination of $\alpha\in\{\SI{0}{\degree},\SI{0.5}{\degree},\SI{1}{\degree},\SI{2}{\degree}\}$, $L\in\{1,3,15\}$, $d\in\{1,2,3\}$ and crystal orientations [001] and [011] (see Figure~\ref{fig: phantom crystals}). The full phantom objects were \SI{1000}{\angstrom} thick.
	\item 2160 kinematical simulated diffraction patterns with 30 random phantoms for each combination of $\alpha\in\{\SI{0}{\degree},\SI{0.5}{\degree},\SI{1}{\degree},\SI{2}{\degree}\}$, $L\in\{1,3,15\}$, $d\in\{1,2,3\}$ and crystal directions [001] and [011] parallel to the optic axis ($z$-axis). The full crystals were \SI{250}{\angstrom} thick.
	\item 540 high-energy kinematical simulated diffraction patterns without precession and with 30 random phantoms for each combination of $L\in\{1,3,15\}$, $d\in\{1,2,3\}$ and crystal directions [001] and [011] parallel to the optic axis ($z$-axis). The full crystals were \SI{250}{\angstrom} thick.
	\item 100 kinematical simulated diffraction patterns of the dislocation phantom with $\alpha=\SI{2}{\degree}$ at different beam locations. The full crystal was \SI{250}{\angstrom} thick.
\end{itemize}

\begin{figure}[ht]
	\begin{subfigure}{.49\textwidth}
		\includegraphics[width=\textwidth,trim={132 30 580 48},clip]{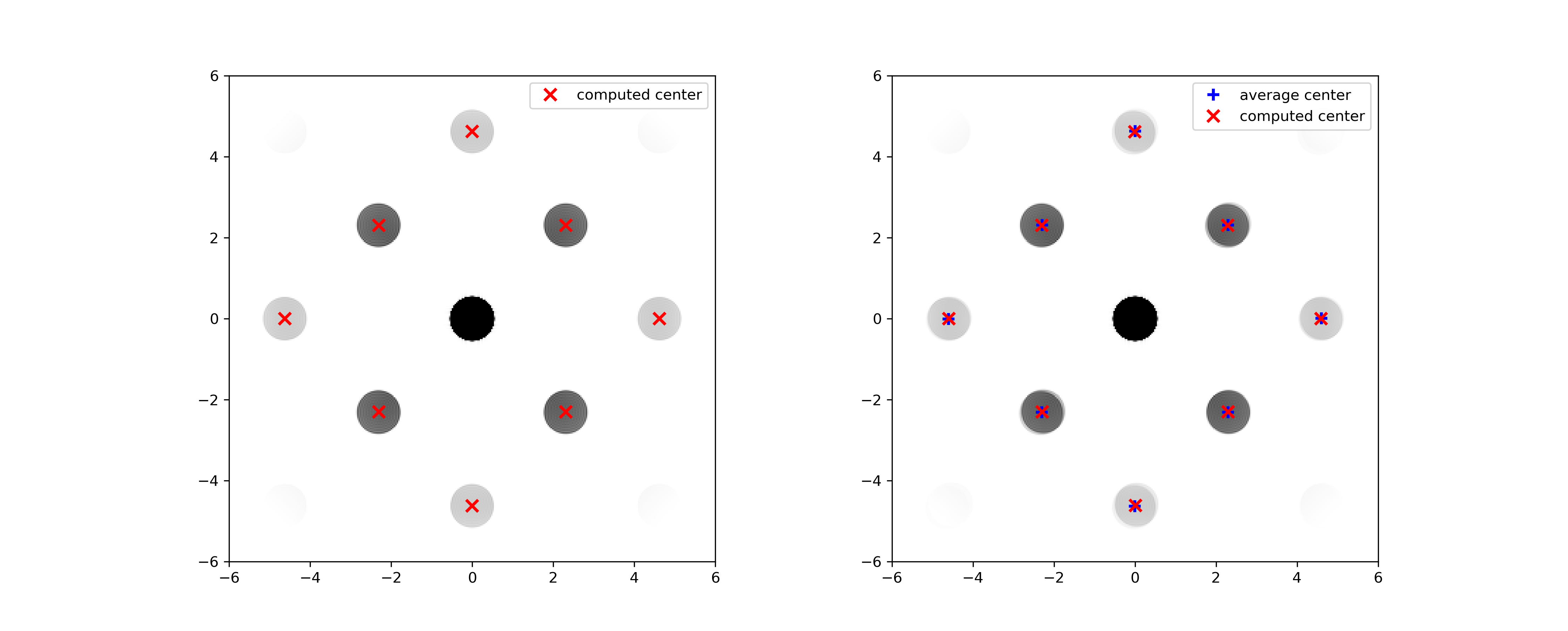}
		\caption{Undeformed diffraction pattern}
	\end{subfigure}
	\begin{subfigure}{.49\textwidth}
		\includegraphics[width=\textwidth,trim={590 30 122 48},clip]{bad_prediction}
		\caption{Deformed diffraction pattern}
	\end{subfigure}
	\begin{subfigure}{.49\textwidth}
		\includegraphics[width=\textwidth,trim={580 30 115 48},clip]{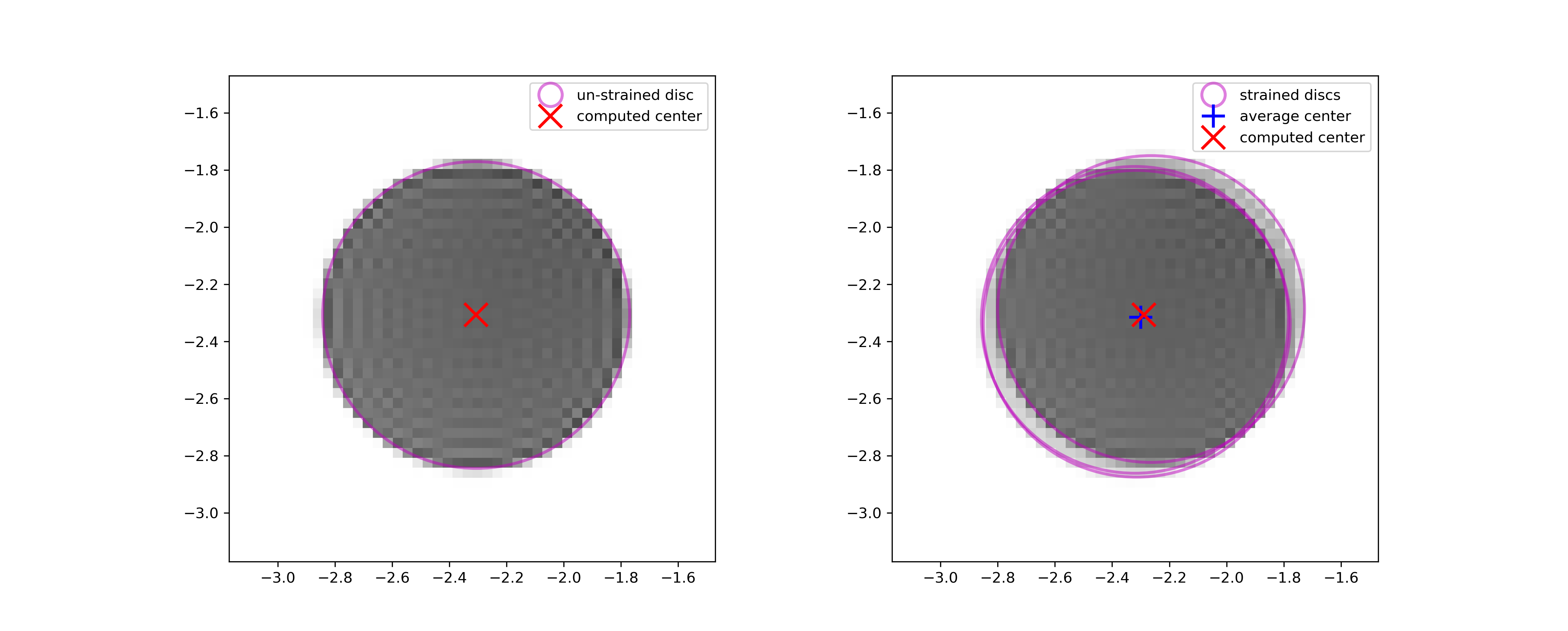}
		\caption{Superimposed deformed discs}
	\end{subfigure}
	\begin{subfigure}{.49\textwidth}
		\includegraphics[width=\textwidth,trim={580 30 115 48},clip]{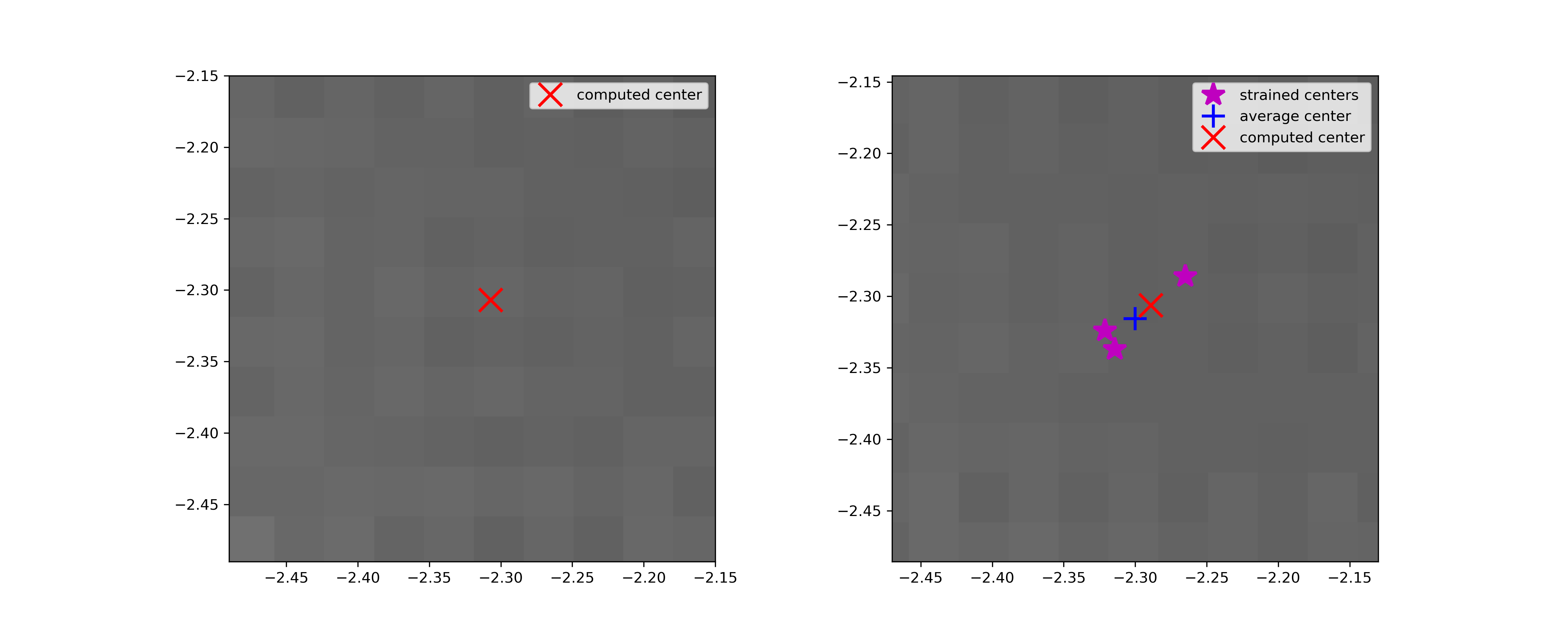}
		\caption{Distribution of deformed centres}
	\end{subfigure}
	\caption{Evaluation pipeline: we compute centres on undeformed (a) and deformed crystals (b). The TRT uses (a) to predict an average centre (blue cross) at the centre of mass of deformed centres (purple stars). (c) and (d) give sub-plots of (b) to visualise the effect of deformation.}\label{fig: examples}
\end{figure}

In particular, we used a disk-detection method involving patch-wise (least squares) registration between each Bragg disk in the strained diffraction pattern and the corresponding Bragg disk in an unstrained diffraction pattern, the exact form of this is in \eqref{eq: centre defs}. A sketch of this pipeline is provided in Figure~\ref{fig: examples}.

For each precessed diffraction pattern we compute centres for each spot in the inner-most ring on the pattern. Predicted and computed centres are only ever compared like-for-like relative to centres computed with the same algorithm using an undeformed sample as reference.

In the notation of \eqref{eq: com approx}, we define the true centres of each spot as $$\vec c_{true} = \E_j \gamma^\top A_j\vec p_i$$
\edit{in the case of}{for} discrete piece-wise affine deformation, \edit{}{and }for the continuous deformation map $\vecf R$
$$\vec c_{true} = \frac{1}{\SI[parse-numbers=false]{30^2}{\angstrom\squared} \ T}\int_{|x-p_{i,x}|\leq \SI{15}{\angstrom}}\int_{|y-p_{i,y}|\leq \SI{15}{\angstrom}}\int_0^T\gamma^\top\vec p_i + \gamma^\top\nabla\vecf R(x,y,z)\vec p_i \D z\D[]y\D[]x$$
computes an average deformation where $T$ is the known thickness of the phantom.

\begin{equation}\label{eq: centre defs}
    \vec c_{com} = \frac{\int_{|\vec k-\gamma^\top \vec p_i|<\bar r} \vec k D_\alpha(\vec k)\D\vec k}{\int_{|\vec k-\gamma^\top \vec p_i|<\bar r} D_\alpha(\vec k)\D\vec k}, \qquad\vec c_{reg} = \argmin_{\vec p} \int_{|\vec k-\gamma^\top \vec p_i|<\bar r} |D_\alpha(\vec k) - D_\alpha^0(\vec k+\vec p-\gamma^\top\vec p_i)|^2\D\vec k
\end{equation}
where $\gamma^\top\vec p_i$ is computed from an undeformed diffraction pattern, $D^0_\alpha$. We fit a zero-mean Gaussian to the TRT modelling errors and so the important error measure is the Euclidean distance between detected and expected centres, i.e. the error variance. We report the values of 
\begin{equation}\label{eq: error metric}
    \text{error } = 100\cdot \frac{|\vec c_{true}-\vec c|}{|\vec c_{true}|}
\end{equation}
which is scaled to percentage error. This is convenient because with $\sigma=1\%$, a naive detection algorithm of $\vec c = \gamma^\top\vec p_i$ (i.e. zero strain) corresponds to an average error of 1. 

\begin{table}
	\centering{%
	\begin{tabular}{|c|c|c|c|c|c|}
		\hline\multirow{2}{*}{Method} & \multicolumn{4}{|c|}{dynamical simulation} & \multirow{2}{*}{high-energy} \\
		               & $\alpha=\SI{0}{\degree}$ &$\alpha=\SI{0.5}{\degree}$& $\alpha=\SI{1}{\degree}$ & $\alpha=\SI{2}{\degree}$ &      \\ \hline
		Centre of Mass &   1.95           &   0.40           &    0.20          & 0.08 &   0.03 \\ 
		Registered     &   0.46           &   0.10           &    0.05          & 0.04 &   0.04 \\\hline
		\multicolumn{6}{c}{}\\
		\hline\multirow{2}{*}{Method} & \multicolumn{4}{|c|}{kinematical simulation} & dislocation \\
		               & $\alpha=\SI{0}{\degree}$ &$\alpha=\SI{0.5}{\degree}$& $\alpha=\SI{1}{\degree}$ & $\alpha=\SI{2}{\degree}$ & phantom \\ \hline
		Centre of Mass &   0.44           &   0.27           &    0.08          & 0.04 &   0.12 \\ 
		Registered     &   0.11           &   0.06           &    0.05          & 0.04 &   0.10 \\\hline
	\end{tabular}%
	}
	\caption{All values given are mean relative Euclidean error \eqref{eq: error metric}. The dislocation phantom is described in Section~\ref{sec: numerics modelling analytical}, the remainder use layered phantoms. The dynamical simulation is computed by MULTEM, high-energy/kinematical from \eqref{eq: high energy equation}/\eqref{eq: precessed sphere model}.} \label{tbl: quantitative results}
\end{table}

Table~\ref{tbl: quantitative results} summarises the results of this comparison. The key observations are:\nopagebreak
\begin{itemize}
    \item An average error of one pixel width would correspond to an error of 0.7\%, all results below this are super-resolved.
	\item Increasing the precession angle in this range reduces the errors for all models and centre detection methods.
	\item Comparing centre detection algorithms, both have comparable maximum accuracy yet the registration method appears much more robust to changes in the simulation model and phantom. It also converges much faster with respect to precession angle with little gain between $\SI{1}{\degree}$ and $\SI{2}{\degree}$.
	\item Errors for the continuously deformed phantom are noticeably worse than with the piece-wise constant phantoms. An example disc is shown in Figure~\SUPP{\ref{fig: continuous deformation disc}}{S.1}, and we see it is qualitatively very different from those shown previously in Figure~\ref{fig: phantom crystals}. The more smoothly varying deformation causes an elliptical blurring of the disc which may make it harder to consistently detect the centre.
	\item The high-energy model achieved optimal accuracy without precession. This suggests the dominant benefit of precession in these examples is smoothing the non-linearities of the model rather than accounting for rotation out-of-plane.
\end{itemize}

In a worst case, with $\SI{2}{\degree}$ of precession, we observe errors between peak-finding and the TRT model approximation of $0.12\%$, corresponding to a signal-to-noise ratio (SNR) over 8. On the other hand, the registration peak finding algorithm consistently achieves an average accuracy of $0.04\%$ corresponding to an SNR of $25$. 

\subsection{Tomographic Reconstruction Validation}\label{sec: synth recon}

Section~\ref{sec: Tensor Tomography} considers the analytical properties of the continuous inverse problem we wish to solve, but in practice we have corrupted and limited data. Here, we perform a reconstruction from a realistic dataset and analyse the accuracy both quantitatively and qualitatively. For this purpose we generate a phantom, $\vecf F^\dagger\colon [-1,1]^3\to \F R^{3\times3}$ and simulate data using the model
$$\vecf d(\vec \xi,\vec x) = J\vecf F^\dagger(\vec \xi,\vec x) + \Pi_{\vec \xi}\vecf \eta(\vec \xi,\vec x)\Pi_{\vec \xi}$$
where $\vecf \eta$ is a 0-mean isotropic white noise tensor field. 

\begin{figure}
    \centering
	\begin{subfigure}[b]{.49\textwidth}
		\includegraphics[width=\textwidth, trim={890 0 467 0},clip]{multi-scale_EM}
		\caption{Tilt-rotate SPED}
	\end{subfigure}
	\begin{subfigure}[b]{.49\textwidth}
		\includegraphics[width=\textwidth,trim={0 0 0 0}, clip]{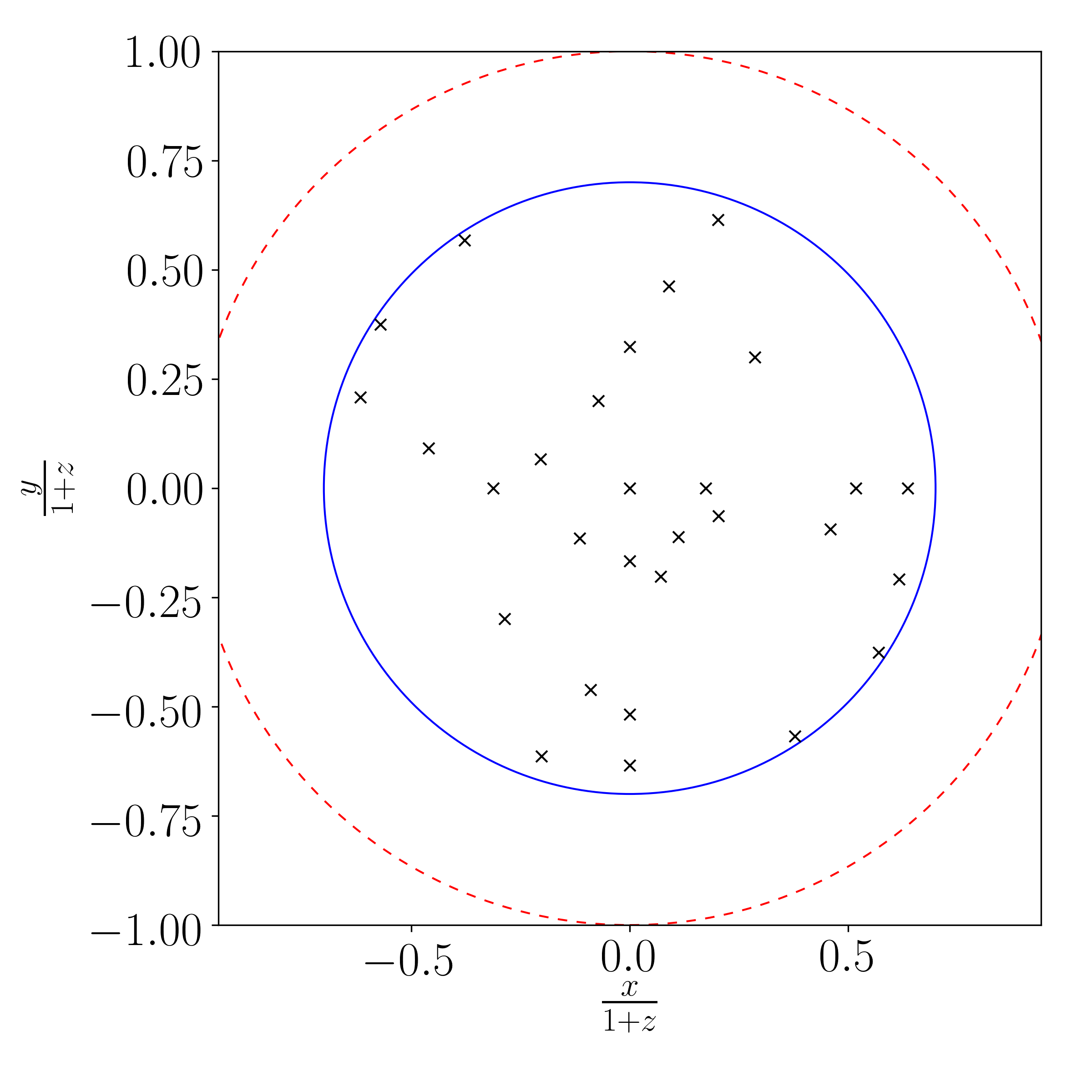}
		\caption{[001] pole figure}
	\end{subfigure}
	\caption{Visualisation of acquisition geometry. (a) Tilt-rotate specimen holders allow two degrees of rotation in the form of Euler angles. The `rotate' angle has full range and the `tilt' is limited to \SI{70}{\degree}. (b) Stereographic projection of chosen tilt directions. Crosses indicate 42 zone axes within the limited tilt range (\SI{70}{\degree} indicated by solid blue line). The full \SI{90}{\degree} tilt range is indicated by the dashed red line. } \label{fig: stereographic tilts}
\end{figure}

The choice of noise level and phantom were chosen to be slightly more challenging than suggested by the results of Section~\ref{sec: numerical simulation results}. In particular, a phantom is chosen with $L=3$ layers, $d=3$ dimensional deformation, and an average deformation magnitude of $\sigma =2\%$. We perform reconstructions with two levels of noise. Firstly at $0.1\%$, in line with Table~\ref{tbl: quantitative results}, and then at $1\%$ (SNR $=2$) to account for any noise and physical modelling errors not considered previously. The final experimental choice is to specify a scan geometry which respects practical time and hardware constraints. 42 tilt directions, shown in Figure~\ref{fig: stereographic tilts}, were selected and for each tilt we scan over a $50\times50$ grid of beam positions. Directions are chosen down zone axes to guarantee at least two non-colinear Bragg discs in the diffraction patterns, as required by Lemma~\ref{thm: noncolinear spots}. For any rotation $R\in\F R^{3\times3}$ of the sample, the corresponding direction is the third column of $R$, \edit{}{say }$(x,y,z)^\top$. This direction is then plotted with a stereographic projection at point $\left(\frac{x}{1-z},\frac{y}{1-z}\right)$. The plot in Figure~\ref{fig: stereographic tilts}.b assumes that the crystal is initially perfectly aligned with the $z$-axis. If this is not true then an offset must be computed and passed on \edit{}{for }the computation of Euler angles for the specimen holder.

For a reconstruction method we choose to perform a standard Total Variation reconstruction \cite{Goris2012,Leary2013,Collins2017}:
$$\vecf F_{recon} = \argmin_{\vecf F} \frac12\norm{J\vecf F-\vecf d}_2^2 + \beta\int_{[-1,1]^3}|\nabla \vecf F(\vec x)|_{Frobenius}\D\vec x$$
where $\beta$ is a manual tuning parameter. With perfect data there is little need for regularisation ($\beta=0$) however, in this example the Total Variation functional is compensating for:
\begin{itemize}
    \item Measurement noise, representing modelling errors at either at $0.1\%$ or $1\%$ magnitude
    \item Limited angular range, within $\SI{70}{\degree}$ of the initial orientation
    \item Limited projections, 42 projections with a $50\times50$ grid of beam positions
    \item Analytical null space, part of the skew component of the tensor field is unobserved by the TRT ($J\vecf F = J(\vecf F + \cross{\nabla \varphi})$ for all $\varphi\in C_0^1$) as discussed in Section~\ref{sec: Tensor Tomography}.
\end{itemize}
While each of these factors has a different physical or analytical origin, numerically they are all incorporated into a choice of $\beta>0$. In both reconstructions the parameter was coarsely tuned to $\beta=\frac{10^{-4}}{2}$.

Many other choices of variational methods exist, for instance those compared in \cite{Leary2013}, which each account for noise and `fill in' missing data in their own characteristic fashion. Total Variation is commonly chosen because it promotes sparse jumps in the reconstruction \cite{Leary2013,Ehrhardt2015}. This specifically reflects the structure of the phantom in this section but is also often accurate for other physical samples. \edit{}{For a direct FBP-like reconstruction algorithm, c.f. \cite{Sharafutdinov2007}.}

\begin{figure}
    \begin{subfigure}[t]{0.33\textwidth}
        \includegraphics[width=.99\textwidth, trim={5 115 950 159}, clip]{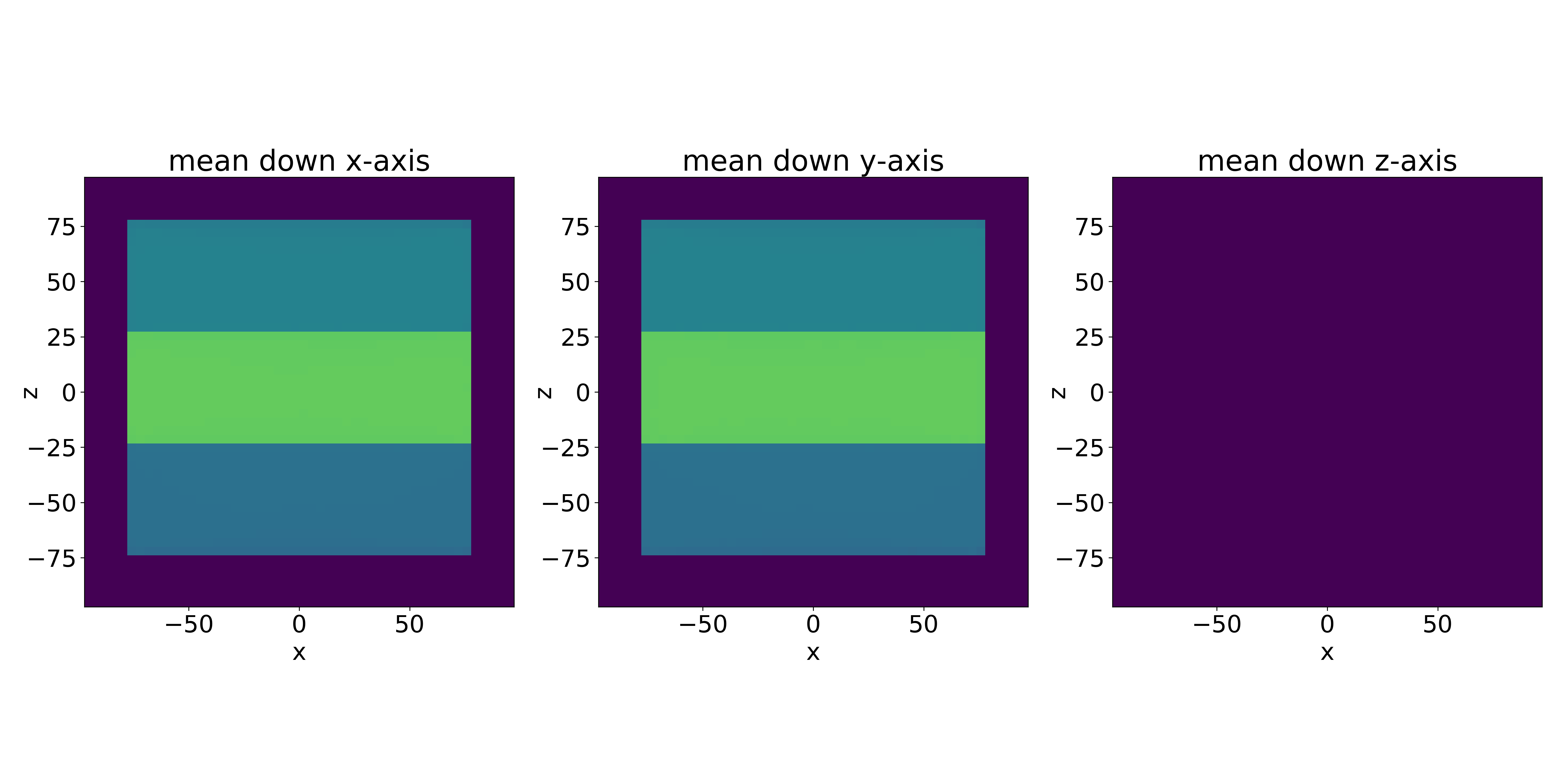}
        \caption{Low noise reconstruction}
    \end{subfigure}
    \begin{subfigure}[t]{0.33\textwidth}
        \includegraphics[width=.99\textwidth, trim={5 115 950 159}, clip]{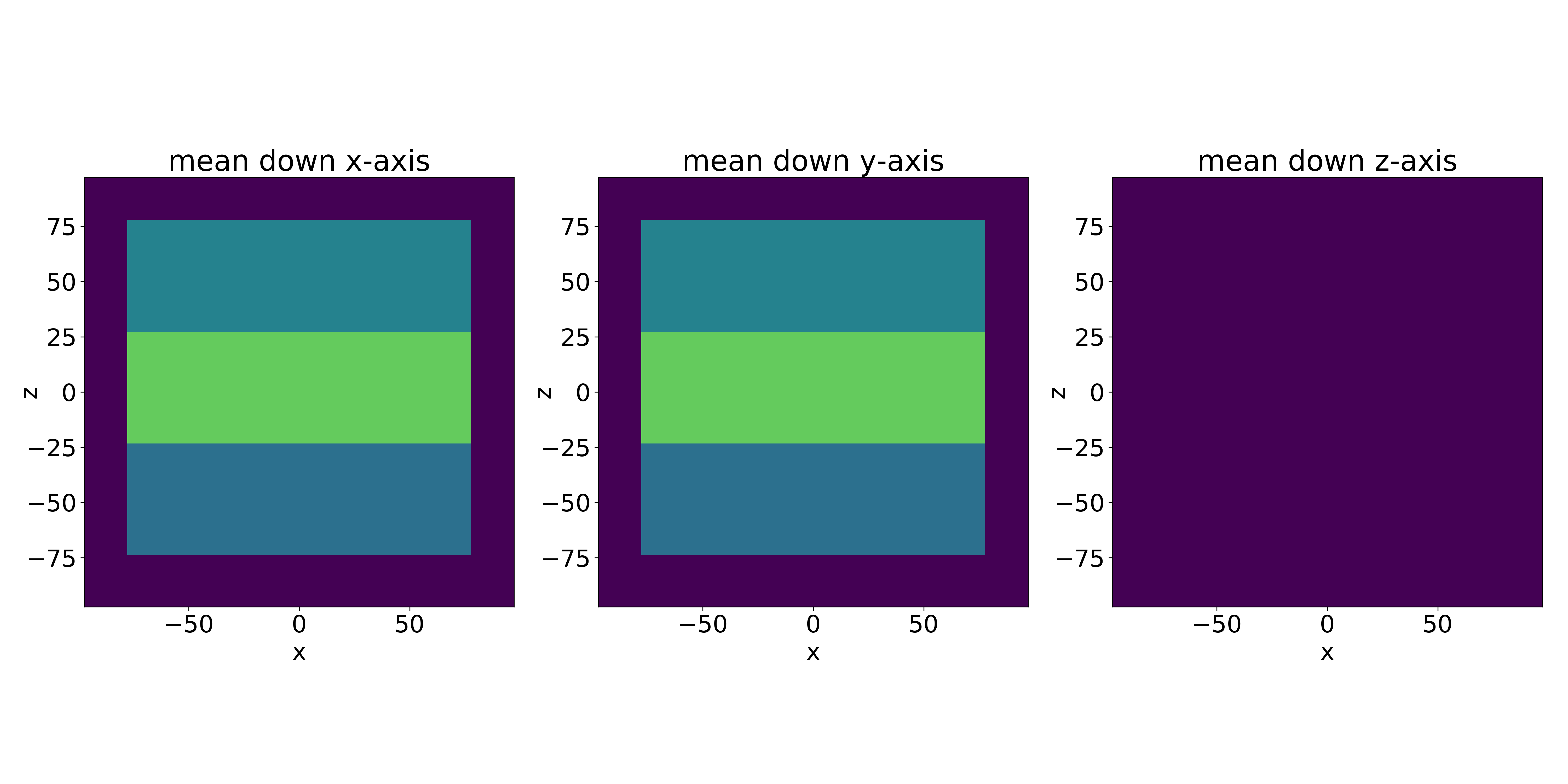}
        \caption{Ground truth}
    \end{subfigure}
    \begin{subfigure}[t]{0.33\textwidth}
        \includegraphics[width=.99\textwidth, trim={5 115 950 159}, clip]{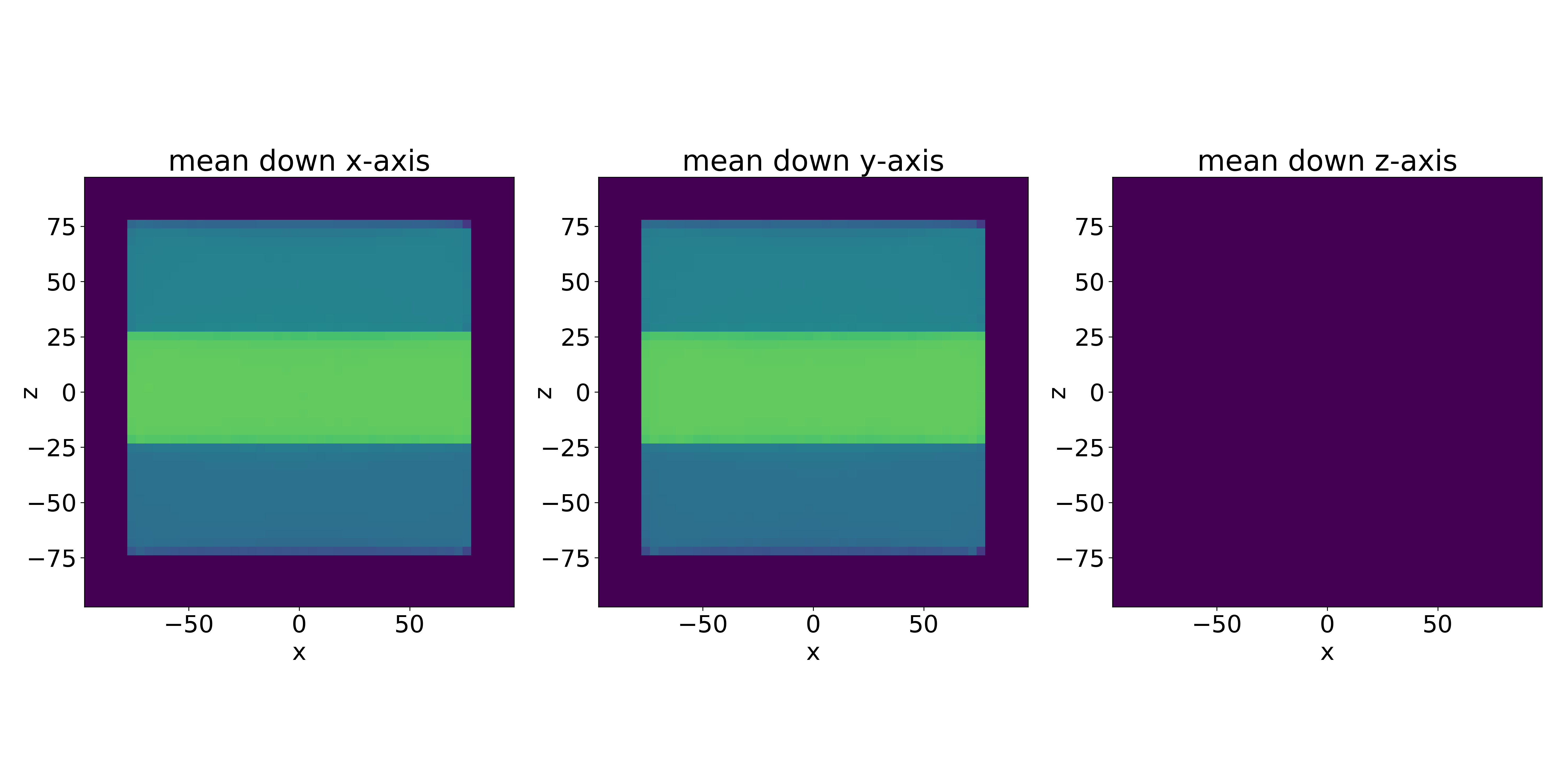}
        \caption{High noise reconstruction}
    \end{subfigure}
	\caption{2D renders of ground truth and reconstructions. The phantom consists of three deformed layers arranged along the $z$-axis. Deformation tensors are max-projected from 5D to 3D then averaged down the $y$-axis into 2D.}
	\label{fig: phantom recon}
\end{figure}

\begin{figure}\centering
	\begin{tikzpicture}
	\node at (0,0) {\includegraphics[height=.3\textheight, trim={20 10 10 5},clip]{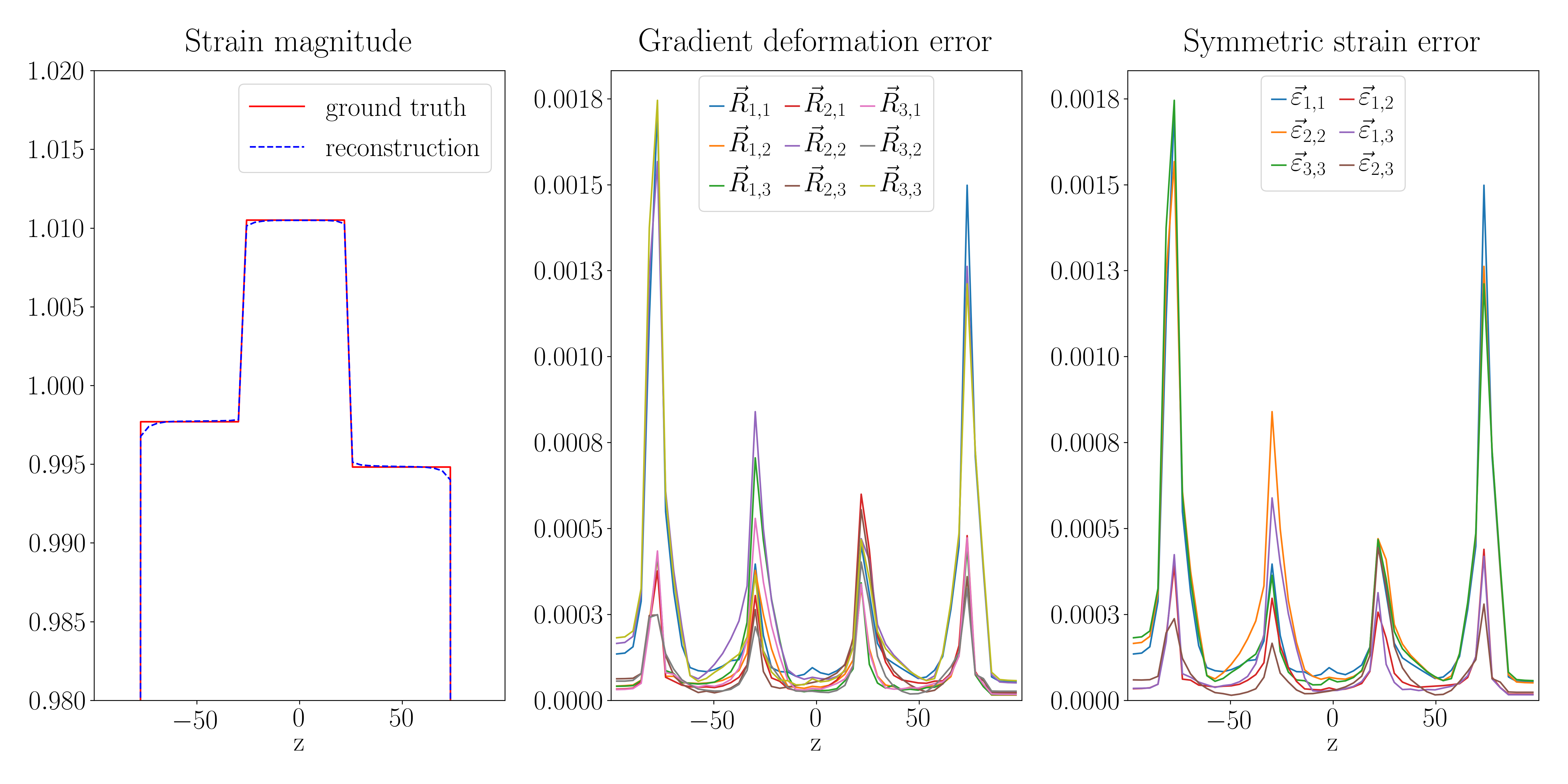}};
	\node at (-.5\textwidth,0) [rotate=90,text width=100pt, align=center] {Low noise};
	\node at (0,-.3\textheight) {\includegraphics[height=.3\textheight, trim={20 10 10 5},clip]{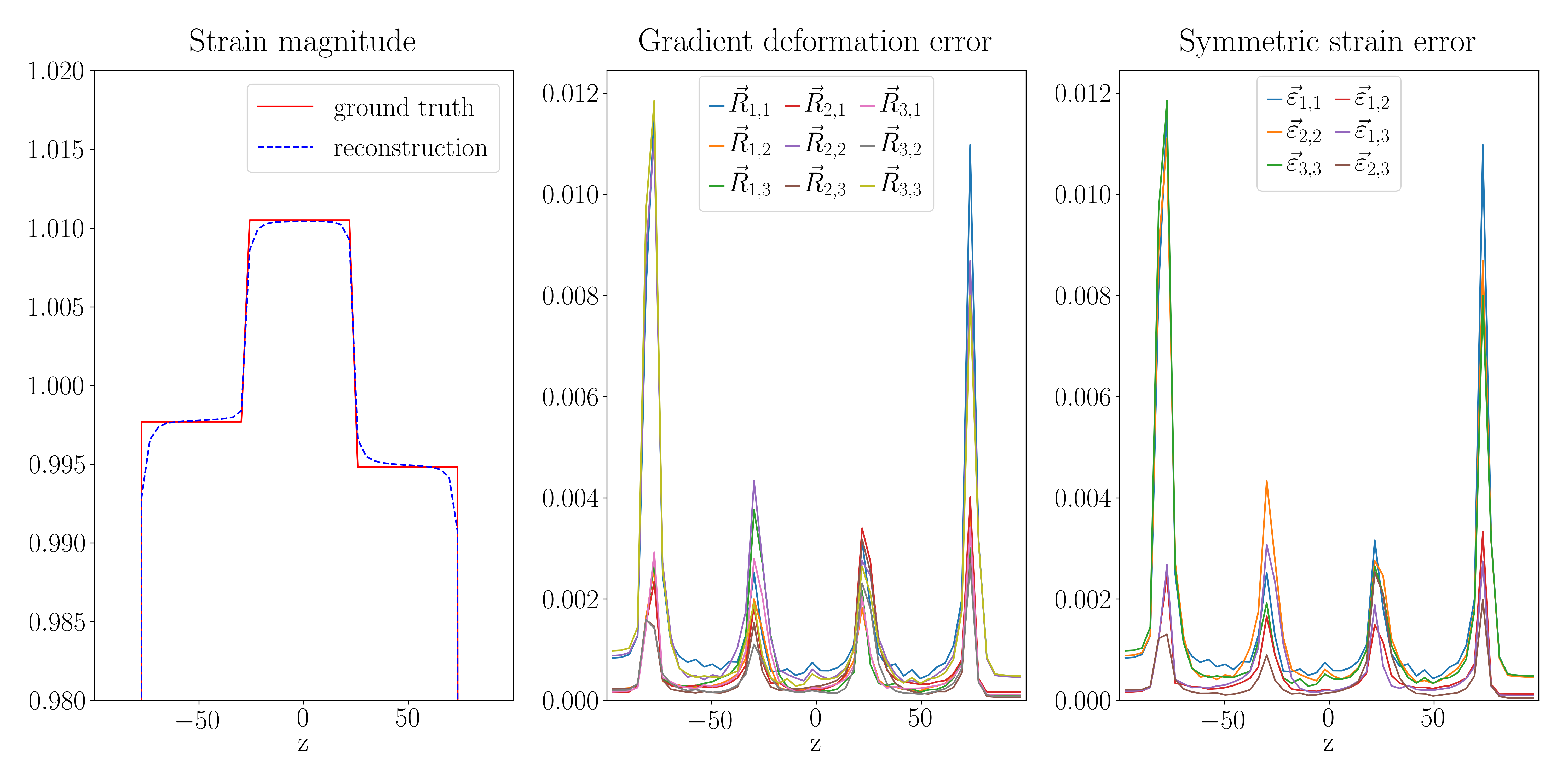}};
	\node at (-.5\textwidth,-.3\textheight) [rotate=90,text width=100pt, align=center] {High noise};
	\end{tikzpicture}
	\caption{1D projections of reconstruction volumes and error distributions. In the volume, tensors are projected to scalars by selecting the maximum. For errors, the $99^{\text{th}}$ percentile of each nine components of the gradient deformation tensor (or six components of symmetric strain) are plotted.}
	\label{fig: phantom recon slices}
\end{figure}

Figure~\ref{fig: phantom recon} visually compares reconstructions from low/high noise data against the original phantom. We see that the general symmetry of the phantom is preserved, reconstructions are uniform in the $x$- and $y$-axes and partition into clear slices along the $z$-axis. The errors in the reconstruction from low noise are imperceptible but at higher noise levels we see the layer interfaces have a slight blur. Figure~\ref{fig: phantom recon slices} allows us to quantify these errors more precisely. The cross-section at low noise shows that the structure of the deformation is well recovered, however, in the flat regions a small consistent error is made in the deformation tensor. Errors at high noise have the same structure but larger magnitude, approximately a factor of $6.5$ larger error for a factor of $10$ larger noise. The interfaces are still well identified however the jump is more visibly blurred. In all cases, approximately 3 pixels away from a discontinuity errors improve rapidly, by up to a factor of 10. More detailed error comparisons are given in Figure~\SUPP{\ref{fig: phantom recon slices comparison}}{S.3} show that the $99^{\text{th}}$ percentile is a representative reference for the general structure of errors. 

\subsubsection{Interpretation of errors}
The three main sources of error are from noise (or modelling error), acquisition geometry, and the null space in the skew component. Figure~\ref{fig: phantom recon slices} allows us to compare the impact of each of these factors. 

Errors from the noise should distribute uniformly over the reconstruction. This error should be constant in uniformly strained regions (near $z=0$ and $z=\pm 50$) and is visibly much lower than the error at interfaces. Comparing between low and high noise, we also see that errors scale approximately linearly with the noise level.

Tomography with a limited angular range, called limited angle tomography, is common in electron microscopy \cite{Quinto1993, Leary2013, Tovey2019}. From this literature, we know that all changes in the deformation in directions near-orthogonal (at angles greater than \SI{70}{\degree}) to the $z$-axis are all missing from the observed data. Uncertainty over where jumps occur lead to blurred edges which are seen as the spreading of errors in the $z$ direction. The jump from crystal to vacuum provides a worst-case scenario and, as previously commented, this blurring is approximately 3 to 4 pixels in radius. This radius is consistent between different noise levels and the full/symmetric strain components.

Comparing the second and third columns of Figure~\ref{fig: phantom recon slices}, it is clear that the error on the symmetric components is no smaller than the full error. This indicates that the contribution of error from the null space of the skew component is negligible.

Returning to the question of symmetric strain reconstruction, there are two conventions for the definition of strain, either
$$ \vecf\epsilon = \frac12(\vecf F + \vecf F^\top) \qquad \text{ or } \qquad\vecf\epsilon =\sqrt{\vecf F^\top\vecf F}.$$
In words, symmetric strain is either the arithmetic mean or the geometric mean. In both cases they are easily computed from the full displacement gradient tensor $\vecf F_{recon}$. The analysis of Section~\ref{sec: Tensor Tomography} nicely aligns with the first definition, which is the default in this work, but the error analysis above shows that this does not bias the quality of the reconstruction. Both definitions of symmetric strain are reconstructed with equivalent levels of accuracy in this example. The total accuracy is dictated much more by the acquisition geometry than noise (representing modelling error) or missing data in the skew component.

\section{Discussion and future steps}\label{sec: conclusion}

In this work we have proposed a tomographic model for the 3D strain mapping problem, analysed and extended the known analytical properties of the resulting inverse problem, and provided numerical results for both modelling and inversion steps. Table~\ref{tbl: quantitative results} shows that our forward model is accurate up to an SNR of 20 with respect to the TRT model. Beam precession is key to this accuracy and disk registration based methods are more stable to this uncertainty relative to centres of mass. We note that the 0.08\% found here agrees exactly with a comparable quantity of \cite{Mahr2015} despite very different deformations considered in each study. In Figure~\ref{fig: phantom recon slices}, we have shown that reconstructions can be performed with realistic experimental parameters and achieve accurate results, even with errors much larger than predicted. Because of this, it is our belief that there is no benefit to quantifying errors beyond this SNR of 20 within the scope of this work, namely the validity of the dynamical model tested in Section~\ref{sec: numerics}. Realistic reconstructions can be recovered well at this level of error, although this does not take into account many factors such as detector performance, electron optical aberrations, and inelastic scattering. These factors could potentially be the dominant causes of reconstruction error in practice and should be minimised experimentally and assessed further.

Our proposed framework requires diffraction patterns to be recorded near zone axes where diffraction patterns have straight lines of spots. It would be much faster experimentally to acquire many diffraction patterns with single arced lines of spots, called `off-axis' diffraction patterns, however this would take us away from the TRT model. On the theoretical side, there have been recent advances in histogram tomography which could assist the well-posedness of strain tomography \cite{Lionheart2018}. In particular, in this study we compute the centre of mass (first moment) of each spot in a diffraction pattern and the resulting inverse problem is the TRT which always has a large null-space. If we also extracted second order moments from each spot, then the model is no longer the TRT and the extra data may remove the issue of non-unique solutions. Finally, there is an interesting conflict in the desired scan geometry due to the physical model and the TRT. As commented in Section~\ref{sec: physics vs mathematical TRT}, theoretical results of the TRT rely on three orthogonal tilt series whereas we only have access to data at a discrete set of orientations which, dependent on the crystal structure, arise (approximately) uniformly over the sphere. It would be interesting to unify these two pressures and analyse characteristics of the tomography problem in such a constrained geometry.

\ack
R.T. acknowledges funding from EPSRC Grant No. EP/L016516/1 for the Cambridge Centre for Analysis. C.-B.S. acknowledges support from the Leverhulme Trust project `Breaking the non-convexity barrier', EPSRC Grant No. EP/M00483X/1, EPSRC centre Grant No. EP/N014588/1, and from CHiPS (Horizon 2020 RISE project grant). R.T. and C.-B.S. also acknowledge the Cantab Capital Institute for Mathematics of Information. M.B. acknowledges the Leverhulme Trust Early Career Fellowship Grant No. ECF-2016-611. S.M.C. acknowledges the Henslow Research Fellowship at Girton College, Cambridge and a University Academic Fellowship at the University of Leeds. W.R.B.L. acknowledges the Royal Society Wolfson Research Merit Award and the EPSRC for EP/P02226X/1, EP/M022498/1, EP/K00428X/1. P.A.M. was supported by the EPSRC under grant EP/R008779/1. We acknowledge the support of NVIDIA Corporation for the provision of the Quadro P6000 and TITAN Xp GPUs used for this research.

\section*{References}
\bibliography{ref}

\begin{appendix}
\section{Lemma~\ref{thm: spike locations}}\label{app: spike locations}
\begin{theorem*}
If $u_0$ has reciprocal lattice vectors $\vec a^*$, $\vec b^*$ and $\vec c^*$ then, with $V \coloneqq \op{det}((\begin{smallmatrix}\vec a^*&\vec b^*& \vec c^*\end{smallmatrix}))$, $u_0$ is periodic in the vectors
    $$ \vec a = \frac{2\pi}{V}\vec b^*\times\vec c^*, \qquad \vec b = \frac{2\pi}{V}\vec c^*\times\vec a^*, \quad \text{and}\quad \vec c = \frac{2\pi}{V}\vec a^*\times\vec b^*.$$
    In other words, there exists a repeating unit $v\colon[0,\max(|\vec a|,|\vec b|,|\vec c|)]^3\to\F R$ such that 
    $$u_0(\vec x) = \sum_{i_1,i_2,i_3=-\infty}^{\infty} v(\vec x + i_1\vec a+i_2\vec b+i_3\vec c).$$	
\end{theorem*}
\begin{proof}
	By direct computation,
	$$u_0(\vec x) = (2\pi)^{-3}\sum_{i=1}^\infty w_i\int_{\F R^3} \delta_{\vec p_i}(\vec y) e^{\i \IP{\vec x}{\vec y}}\edit{}{\D\vec y} = (2\pi)^{-3}\sum_{i=1}^\infty w_ie^{\i \IP{\vec x}{\vec p_i}}.$$
	Suppose $\vec a = \frac{2\pi}{V}\vec b^*\times \vec c^*$, for $\vec p_i = i_1\vec a^* + i_2\vec b^* + i_3\vec c^*$ we have
	$$\exp(\i\IP{\vec a}{\vec p_i}) = \exp(\i\tfrac{2\pi}{V} i_1\IP{(\vec b^*\times \vec c^*)}{\vec a^*}) = \exp(\i\tfrac{2\pi}{V} i_1V) = 1. $$
	Therefore $u_0(\vec x+\vec a) = u_0(\vec x)$ for all $\vec x$. By symmetry, the same holds for $\vec b$ and $\vec c$. Define
	$$v(\vec x) = \begin{cases}
	u_0(\vec x) & 0\leq \IP{\vec x}{\vec a}< |\vec a|^2, 0\leq \IP{\vec x}{\vec b}< |\vec b|^2 \text{ and } 0\leq \IP{\vec x}{\vec c}< |\vec c|^2
	\\0 & \text{else}
	\end{cases}.$$
	Note that $\vec a$, $\vec b$ and $\vec c$ span $\F R^3$, therefore for any $\vec x$ there exist $x_j\in\F R$ such that $\vec x = x_1\vec a + x_2\vec b + x_3\vec c$. We can conclude
	\begin{align*}
	\sum_{\F Z^3} v(\vec x + i_1\vec a+i_2\vec b+i_3\vec c) &= \sum_{\F Z^3} v((x_1+i_1)\vec a + (x_2+i_2)\vec b + (x_3+i_3)\vec c)
	\\&= v((x_1-\lfloor x_1\rfloor)\vec a + (x_2-\lfloor x_2\rfloor)\vec b + (x_3-\lfloor x_3\rfloor)\vec c)
	\\&= u_0((x_1-\lfloor x_1\rfloor)\vec a + (x_2-\lfloor x_2\rfloor)\vec b + (x_3-\lfloor x_3\rfloor)\vec c)
	\\&= u_0(\vec x)
	\end{align*}
	as required.
\end{proof}

\section{Theorem~\ref{thm: single strain crystal} and Lemma~\ref{thm: single strain infinite crystal}}\label{app: single strain crystal}
\begin{theorem*}[Theorem~\ref{thm: single strain crystal}]
	If $$u'(\vec x) = u(\vecf E(\vec x)) = u(A\vec x+\vec b) \text{ for some } A\in\F R^{3\times3},\ \vec b\in\F R^3,\ u\in L^2(\F R^{\edit{n}{3}};\F C)$$
	where $A$ is invertible, then we can express its Fourier transform as:
	$$\C F[u'](\vec K) = \op{det}(A)^{-1}e^{\i\IP{\vec b}{A^{-\top }\vec K}}\C F[u](A^{-\top }\vec K).$$
\end{theorem*}
\begin{proof}
	\begin{align*}
	\C F[u'](\vec K) &= \int_{\F R^3} u'(\vec x) \exp\left(-\i\IP{\vec x}{\vec K}\right)\D\vec x
	\\&= \int_{\F R^3} u(A\vec x+\vec b) \exp\left(-\i\IP{(A^{-1}A\vec x)}{\vec K}\right)\D\vec x
	\\&= \int_{\F R^3} u(A\vec x+\vec b) \exp\left(-\i\IP{(A\vec x)}{(A^{-\top }\vec K)}\right)\D\vec x
	\\&= \int_{\F R^3} u(\vec x'+\vec b)\exp\left(-\i\IP{(\vec x'+\vec b-\vec b)}{A^{-\top }\vec K}\right)\frac{\D\vec x'}{\op{det}(A)}
	\\&= \op{det}(A)^{-1}\int_{\F R^3} u(\vec x')\exp\left(-\i\IP{\vec x'}{A^{-\top }\vec K}+\i\IP{\vec b}{A^{-\top }\vec K}\right)\D\vec x'
	\\&= \op{det}(A)^{-1}e^{\i\IP{\vec b}{A^{-\top }\vec K}}\C F[u](A^{-\top }\vec K).
	\end{align*}
\end{proof}

\begin{theorem*}[Lemma~\ref{thm: single strain infinite crystal}]
	$$\C F[u_0(A\cdot + \vec b)](\vec K) = e^{\i\IP{\vec b}{A^{-\top }\vec K}}\sum_{i=1}^\infty w_i\delta_{A^\top \vec p_i}(\vec K)$$
	for all $A\in\F R^{3\times3}, \vec b\in\F R^3$.
\end{theorem*}
\begin{proof}
	From Theorem~\ref{thm: single strain crystal}, we know 
	\begin{align*}
	\C F[u_0(A\cdot+\vec b)](\vec K) &= \op{det}(A)^{-1}e^{\i\IP{\vec b}{A^{-\top }\vec K}}\C F[u_0](A^{-\top }\vec K)
	\\&= \sum_{i=1}^\infty w_i\op{det}(A)^{-1}e^{\i\IP{\vec b}{A^{-\top }\vec K}}\delta_{\vec p_i}(A^{-\top }\vec K)
	\end{align*}
	Thus, to complete the lemma, it suffices to show
	$$\op{det}(A)^{-1}\delta_{\vec p}(A^{-\top }\vec K) = \delta_{A^\top \vec p}(\vec K)$$
	This is verified by an arbitrary test function, $\varphi\in C_c^\infty$
	$$ \int_{\F R^3} \op{det}(A)^{-1}\delta_{\vec p}(A^{-\top }\vec K) \varphi(\vec K) \D\vec K = \int_{\F R^3} \delta_{\vec p}(\vec K) \varphi(A^\top \vec K) \D\vec K = \varphi(A^\top \vec p). $$
\end{proof}

\section{Probability background}\label{app: probability background}
We recap some concepts and technical results from probability theory which are needed in the proofs in \ref{app: approximations}.
\begin{definition}\hfill
	\begin{itemize}
		\item $X\colon t\mapsto X_t\in \F C$ can be called a \emph{random variable} where random (complex) values of $X$ can be sampled by sampling indices $t\in[t_0,t_1]$ uniformly at random, i.e. $t\sim \op{Uniform}[t_0,t_1)$.
		\item For a random variable, $X$, we define its \emph{expectation} to be
		$$ \E X = \E_tX_t = \frac{\int_{t_0}^{t_1} X_t \D t}{|t_1-t_0|}.$$
		If $t$ is a discrete index then integral can be replaced with summation and, for $t_0,t_1\in\F Z$, this becomes 
		$$\E X = \E_tX_t = \frac{\sum_{t_0}^{t_1-1} X_t}{|t_1-t_0|}.$$
		In quantum physics this would also be called the \emph{expectation value}, $\E X = \langle X\rangle$.
		\item For a random variable $X$ and value $x\in\F C$ we say denote the probability that $X_t=x$ to be $$\op{probability}(X = x) = \op{probability}(X_t = x) = \frac{|\{t\text{ s.t. } X_t=x\}|}{|t_1-t_0|}.$$
		\item Let $X$ and $Y$ be random variables. For $x,y\in\F C$, we say the events $X_t=x$ and $Y_t=y$ are \emph{independent} if 
		$$\op{probability}(X_t=x \text{ and } Y_t=y) = \op{probability}(X_t=x) \cdot \op{probability}(Y_t=y).$$
		\item We say that two random variables, $X$ and $Y$, are \emph{independent} if 
		$$\text{events } X_t=x \text{ and } Y_t=y\text{ are independent for all } x,y.$$
	\end{itemize}
\end{definition}
The following are two standard results from probability theory.
\begin{lemma}\label{thm: product of independent}
	If $X$ and $Y$ are two independent random variables then
	$$\E(XY) = \E(X)\E(Y).$$
	Also, for any continuous function, $\varphi$, the random variable $\varphi(X)_{i,j,t} = \varphi(X_{i,j,t})$ is also independent of $Y$.
\end{lemma}
\begin{lemma}\label{thm: average paths}
	Let $\varphi\colon\F R^n\to\F R$ be a twice differentiable function and let $\eta\colon[0,1]\to\F R^n$ be a random variable, then $$\E_t \varphi(\eta_t) = \varphi(\E_t\eta) + O\left(\norm{\eta-\edit{\E}{\E_t}\eta}_\infty^2\norm{\nabla^2\varphi}_\infty\right).$$
\end{lemma}
\begin{proof}
    \edit{}{Denote $\bar\eta \coloneqq \E_t\eta$.} Because the length of the path (size of domain of $\eta$) is one, we can replace expectation with integration and then replace $\varphi$ with its standard Taylor expansion:
	\begin{align*}
	\E_t\varphi(\eta_t) &= \int_0^1\varphi(\eta_t)\D t \qquad\qquad\leadsto \text{(Taylor expansion about } \edit{\E}{\bar\eta})
	\\&= \int_0^1 \left[\varphi(\edit{\E}{\bar\eta}) + \IP{(\eta_t-\edit{\E}{\bar\eta})}{\nabla \varphi(\edit{\E}{\bar\eta})} + O\left(|\eta_t-\edit{\E}{\bar\eta}|^2\norm{\nabla^2 \varphi}_\infty\right)\right] \D t
	\\&= \varphi(\edit{\E}{\bar\eta}) + \IP{\left(\int_0^1 \eta_t-\edit{\E}{\bar\eta}\D t\right)}{\nabla \varphi(\edit{\E}{\bar\eta})} + O\left(\norm{\eta-\edit{\E}{\bar\eta}}_\infty^2\norm{\nabla^2\varphi}_\infty\right)
	\\&= \varphi(\edit{\E}{\bar\eta}) + \IP{(\E_t[\eta]-1\edit{\E}{\bar\eta})}{\nabla \varphi(\edit{\E}{\bar\eta})} + O\left(\norm{\eta-\edit{\E}{\bar\eta}}_2^2\norm{\nabla^2\varphi}_\infty\right)
	\end{align*}
	as required.
\end{proof}

\section{Algebraic justification of Section~\ref{sec: approximations}} \label{app: approximations}
In Section~\ref{sec: approximations} we give a physical motivation for why centre of mass calculations should accurately predict average strain values. In this section we provide a formal mathematical link from the precessed Ewald sphere model to this same goal. Some of the statistical assumptions made here are highly technical and it would be hard to justify them from a purely physical perspective. However, as in the main text, the core justification which we rely upon is the simulation study of the whole pipeline.

In the remainder of this section, we first sketch the proof at a high level listing a sequence of results, and then the details of the longer proofs appear at the end of the section.

To recap, the starting point of this argument is the physical model which we assume to be exact. In particular,
$$D(\vec k) = |\C F[\Psi_p]\star\C F[u]|^2(\vec k,k_z(\vec k))$$
where $$ u(\vec x) = \sum_{j=1}^N u_0(A_j\vec x+\vec b_j)\F V\left(\frac{\vec x-\vec \beta_j}{\rho}\right).$$
The first step is to spatially localise the strain which generates the signal in a diffraction pattern. A key part of tensor tomography is that signal only depends on the strain contained in blocks (or voxels) which directly intersect the electron beam. This is formalised in the following result which also expands the notation to take advantage of the particular structure of $u$.
\begin{theorem*}[Theorem~\ref{thm: full strain model}]
	If \ref{ass: local signal} (narrow probe) holds then
	\begin{align*}
	D(\vec k) &= \left|\sum_{i\in\F N,j\in[N]}\hat w_{i}(k_z(\vec k) - (A_j^\top \vec p_i)_z, \beta_{j,z}) e^{\i\IP{\vec b_j}{\vec p_i}} f(\vec k - \gamma^\top A_j^\top \vec p_i)\right|^2
	\\\text{ where } & f(\vec k) = \left[\C F[\Psi_p]\star [\rho^2\op{sinc}(\rho\vec k')]\right](\vec k)
	\\ \text{ and }&\hat w_{i}(k,\beta) = w_i\rho\op{sinc}(\rho k) e^{-\i \beta k}.
	\end{align*}
\end{theorem*}
As can be seen, this expression is still highly non-linear in terms of the strain parameter $A_j$ and indeed too non-linear for us to develop a linear model directly. To make this possible we introduce the beam precession technique. The heuristic is that if we modify our data at acquisition, then we can analyse it as if it were generated by a simpler model. There is of course a trade-off here, a small amount precession will usefully smooth the problem, but too much will begin to blur out the desired structures. The following two lemmas make this statement precise and their combination is exactly the statement of Approximation~\ref{thm: precession simplified model}.
\begin{lemma}
	$$D_\alpha(\vec k) = \left|\E_t\C F[\Psi_p u(R_t\vec x)](\vec k, k_z(\vec k))\right|^2 + \op{Var}_t\left[\C F[\Psi_p u(R_t\vec x)](\vec k, k_z(\vec k))\right].$$
\end{lemma}
\begin{proof}
	This is just the definition of variance,
	$$ \E_t |Y_t|^2 = |\E_t Y_t|^2 + \op{Var}_t Y_t$$
	for the appropriate choice of random variable $Y_t$.
\end{proof}
This claim is the most physically ambiguous. It suggests that coherent precession is a good approximation of incoherent precession. The justification of this is that it only needs to be true when looking at the whole pipeline. In particular, the variance term may be large in magnitude but, so long as it does not strongly bias the locations of the centres of the diffracted peaks, then it will not affect the overall approximation.

\begin{lemma}\label{thm: large precession penalty}
	If $A_jR_t\gamma$, $A_jR_t\left(\begin{smallmatrix}
	0\\0\\1\end{smallmatrix}\right)$, and $\vec b_j$ are independent random variables over $(j,t)$  and the high-energy limit is valid (\ref{ass: high energy}), then
	$$\left|\E_t\C F[\Psi_p u(R_t\vec x)](\vec k, k_z(\vec k))\right|^2 = \underbrace{\left| \sum_{i=1}^\infty \bar w_i\E_{\gamma^\top \vec\beta_j=0} f(\vec k-\gamma^\top A_j^\top \vec p_i)\right|^2}_{\eqqcolon\bar D_\alpha} + O(\alpha^2)$$
	for some new weights $\bar w_i\in\F C$.
\end{lemma}
While the assumptions of this lemma appear highly technical, they are also not unreasonable from a physical standpoint. $\vec b_j$ and $A_j$ represent translations and strains respectively. These are physically distinct quantities and so $\vec b_j$ should also be statistically independent \edit{from}{of} the first two terms. Looking closer at these terms, $A_jR_t\gamma$ is the first two columns of $A_jR_t$ and $A_jR_t\left(\begin{smallmatrix} 0\\0\\1\end{smallmatrix}\right)$ is the third. If $\alpha$ is small, then this is just a statement that the out-of-plane strain is independent to the in-plane strain.	

Finally, we need to justify that the centre of mass is an accurate predictor of average strain. The difficulty here again is the squared modulus, if the exit waves were imaged directly then this step would be direct but the square has the potential to bias towards points of maximum intensity. The assumptions necessary at this stage are some form of symmetry on the strain field, the technical statement is as follows.
\begin{lemma}\label{thm: final model exact case}
	Suppose the conditions of Approximation~\ref{thm: precession simplified model} hold and the diffracted spots are symmetric and non-overlapping (\ref{ass: separated spots} and \ref{ass: symmetric}). If the random variables $A_{j_1}+A_{j_2}$ and $A_{j_1}-A_{j_2}$ are independent over random pairs of indices $\{(j_1,j_2) \text{ s.t. } \gamma^\top \vec\beta_{j_1} = \gamma^\top \vec\beta_{j_2} = \vec0\}$ then for each $i$:
	\begin{equation*}\E_{\gamma^\top \beta_j = 0}\gamma^\top A_j^\top \vec p_i = \frac{\int_{|\vec k-\gamma^\top \vec p_i|<\bar r} \vec k \bar D_\alpha(\vec k)\D\vec k}{\int_{|\vec k-\gamma^\top \vec p_i|<\bar r} \bar D_\alpha(\vec k)\D\vec k} 
	\end{equation*}
	whenever the denominator is not 0 and where $\bar r>0$ is the separation of spots given by \ref{ass: separated spots}.
\end{lemma}
Theorem~\ref{thm: final model} is the combination of the previous three lemmas where any violations of the exact conditions becomes absorbed into the error term. We now provide the proofs of these lemmas.

\begin{proof}[Proof of Theorem~\ref{thm: full strain model}.]
	By Theorem~\ref{thm: single strain crystal}, we have
	$$\C F\left[\F V\left(\frac{\vec x-\vec \beta_j}{\rho}\right)\right](\vec K) = \rho^3\op{sinc}(\rho\vec K)e^{-\i\IP{\vec \beta}{\vec K}} = \left[\rho^2\op{sinc}(\rho \vec k)e^{-\i\IP{(\gamma^\top \vec \beta)}{\vec k}}\right]\left[\rho\op{sinc}(\rho k_z)e^{-\i\beta_zk_z}\right]$$
	which can be split into its $(x,y)$ and $z$ components. Combining this again with Theorem~\ref{thm: single strain crystal} we can expand
	\begin{align*}
	\C F[\Psi_p]\star\C F[u] &= \sum_{i\in\F N,j\in [N]}\left( \C F[\Psi_p]\star \C F\left[\F V\left(\frac{\vec x-\vec \beta_j}{\rho}\right)\right]\right)\star \left[w_ie^{\i\IP{\vec b_j}{A_j^{-\top }\vec K}}\delta_{A_j^\top \vec p_i}(\vec K)\right].
	\end{align*}
	To simplify this, observe that for all smooth functions, $\phi,\psi$:
	\begin{align*}
	\phi\star (\psi\delta_{\vec p})(\vec K) &= \int_{\F R^3} \phi(\vec K-\vec K')\psi(\vec K')\delta_{\vec p}(\vec K')\D\vec K'
	\\&= \phi(\vec K-\vec p)\psi(\vec p).
	\end{align*}
	Thus we derive
	$$ \C F[\Psi_p]\star\C F[u](\vec K) = \sum_{i\in\F N,j\in[N]}\left( \C F[\Psi_p]\star \C F\left[\F V\left(\frac{\vec x-\vec \beta_j}{\rho}\right)\right]\right)(\vec K- A^\top _j\vec p_i) \left[w_ie^{\i\IP{\vec b_j}{\vec p_i}}\right].$$
	
	Finally, if the beam is smaller than the width of a single block (from \ref{ass: local signal},  $2r<\rho$) then only one column of blocks directly on the beam path contribute to the diffraction signal. Without loss of generality, this is the set of blocks $j$ such that $\beta_{j,x} = \beta_{j,y} = 0$, equivalently $\gamma^\top \vec\beta_j=\vec 0$. With this simplification, we can expand
	\begin{align*}
	\C F[\Psi_p]\star \C F\left[\F V\left(\frac{\vec x-\vec \beta_j}{\rho}\right)\right] (\vec K) &= \left[\C F[\Psi_p]\star[\rho^2\op{sinc}(\rho\vec\kappa)e^{-\i\IP{(\gamma^\top \vec \beta_j)}{\vec\kappa}}]\right](\vec k)\left[\rho\op{sinc}(\rho k_z)e^{-\i\beta_{j,z}k_z}\right]
	\\&= \left[\C F[\Psi_p]\star[\rho^2\op{sinc}(\rho\vec\kappa)]\right](\vec k)\left[\rho\op{sinc}(\rho k_z)e^{-\i\beta_{j,z}k_z}\right]
	\\&\equiv f(\vec k)\left[\rho\op{sinc}(\rho k_z)e^{-\i\beta_{j,z}k_z}\right].
	\end{align*}
	Substituting this above gives
	$$ \C F[\Psi_p]\star\C F[u](\vec K) = \sum_{i\in\F N,\gamma^\top \vec\beta_j=\vec 0} w_i[\rho\op{sinc}(\rho\cdot)e^{-\i\beta_{j,z}\cdot}](k_z-(A^\top _j\vec p_i)_z) \;e^{\i\IP{\vec b_j}{\vec p_i}} f(\vec k-\gamma^\top A_j^\top \vec p_i)$$
	as required.
\end{proof}

\begin{lemma}
	$$R_t\gamma-\E_tR_t\gamma = O(\alpha).$$
\end{lemma}
\begin{proof}
	This is a purely algebraic proof:
	\begin{align*}
	R_t\gamma &= \left(\begin{smallmatrix} \cos(t)&\sin(t)&0\\-\sin(t)&\cos(t)&0\\0&0&1\end{smallmatrix}\right)\left(\begin{smallmatrix} 1&0&0\\0&\cos(\alpha)&\sin(\alpha)\\0&-\sin(\alpha)&\cos(\alpha)\end{smallmatrix}\right)\left(\begin{smallmatrix} \cos(t)&-\sin(t)&0\\\sin(t)&\cos(t)&0\\0&0&1\end{smallmatrix}\right)\left(\begin{smallmatrix} 1&0\\0&1\\0&0\end{smallmatrix}\right)
	\\&= \left(\begin{smallmatrix} \cos^2(t)+\cos(\alpha)\sin^2(t)&(\cos(\alpha)-1)\cos(t)\sin(t)\\(\cos(\alpha)-1)\cos(t)\sin(t)&\sin^2(t) + \cos(\alpha)\cos^2(t)\\-\sin(\alpha)\sin(t)&-\sin(\alpha)\cos(t)\end{smallmatrix}\right)
	\\&= \gamma + \left(\begin{smallmatrix} (\cos(\alpha)-1)\sin^2(t)&(\cos(\alpha)-1)\cos(t)\sin(t) \\(\cos(\alpha)-1)\cos(t)\sin(t)&(\cos(\alpha)-1)\cos^2(t) \\-\sin(\alpha)\sin(t)&-\sin(\alpha)\cos(t)\end{smallmatrix}\right)
	\\&= \gamma - \left(\begin{smallmatrix} 0&0 \\0&0 \\\sin(t)&\cos(t)\end{smallmatrix}\right)\sin\alpha + O(\alpha^2).
	\end{align*}
	This gives
	\begin{align*}
	R_t\gamma - \E_t R_t\gamma &= -\left[\left(\begin{smallmatrix} 0&0 \\0&0 \\\sin(t)&\cos(t)\end{smallmatrix}\right) - \E_t\left(\begin{smallmatrix} 0&0 \\0&0 \\\sin(t)&\cos(t)\end{smallmatrix}\right)\right]\sin\alpha+ O(\alpha^2) 
	\\&= O(\alpha).
	\end{align*}
\end{proof}

\begin{proof}[Proof of Lemma~\ref{thm: large precession penalty}.]
	Note that we can re-write
	{\small\begin{align*}
		\bigg|\E_t\C F[\Psi_pu(R_t&\vec x)](\vec k,k_z(\vec k))\bigg|^2\notag
		\\ &= \left|\sum_{i\in\F N,\gamma^\top \vec \beta_j=0}\E_t\left[\hat w_i(k_z(\vec k) - ((A_jR_t)^\top \vec p_i)_z, \beta_{j,z}) e^{\i\IP{\vec b_j}{\vec p_i}} f(\vec k - \gamma^\top (A_jR_t)^\top \vec p_i)\right]\right|^2
		\\&= \left|\sum_{i\in\F N}\E_{j,t}\left[X_{i,j,t}(\vec k)Y_{i,j,t}(\vec k)Z_{i,j,t}(\vec k)\right]\right|^2\times|\{j \text{ s.t. } \gamma^\top \vec\beta_j=0\}|^2
		\end{align*}}
	where the constant comes from switching the sum to average over $j$, and we define
	\begin{align*}
	X_{i,j,t}(\vec k) &= \hat w_{i}(k_z(\vec k)-(R_t^\top A_j^\top \vec p_i)_z,\beta_{j,z})
	&&= X_i\left(\vec k|A_jR_t\left(\begin{smallmatrix} 0\\0\\1\end{smallmatrix}\right)\right)
	\\Y_{i,j,t}(\vec k) &= \exp(\i\IP{\vec b_j}{\vec p_i}) &&= Y_i(\vec k|\vec b_j) ,
	\\Z_{i,j,t}(\vec k) &= f(\vec k - \gamma^\top R_t^\top A_j^\top \vec p_i),
	&&= Z_i(\vec k|A_jR_t\gamma).
	\end{align*}
	These functions should be considered as \edit{a set of }{}random variables which, for each $\vec k$, are indexed over $i$ and the indices $(j,t)$ are considered the source of randomness. Restating the assumptions of the theorem onto these variables, we know $X_i$, $Y_i$, and $Z_i$ are independent in $(j,t)$. With this, we can apply Lemma~\ref{thm: product of independent} to simplify
	\begin{align*}
	\bigg|\E_t\C F[\Psi_pu(R_t\vec x)](\vec k,k_z(\vec k))\bigg|^2 &\propto \left|\sum_{i\in\F N}\E_{j,t} \left[X_{i,j,t}(\vec k)Y_{i,j,t}(\vec k)Z_{i,j,t}(\vec k)\right]\right|^2 
	\\&= \left|\sum_{i\in\F N}\E \left[X_{i}(\vec k)\right]\E \left[ Y_{i}(\vec k)\right]\E \left[ Z_{i}(\vec k)\right]\right|^2.
	\end{align*}
	Each of these factors now simplifies:
	\begin{itemize}
		\item Assuming the high-energy limit we have $k_z(\vec k)=0$ and so $\E X_i(\vec k)= \op{constant}_i$.
		\item $Y_i$ is not a function of $\vec k$, so we trivially have $\E Y_i(\vec k) = \op{constant}_i'$.
		\item By \ref{thm: average paths} we can approximate
		\begin{align*}
		\E Z_i(\vec k) = \E_j\E_t f(\vec k -\gamma^\top R_t^\top A_j^\top \vec p_i) 
		&= \E_j f(\vec k -\E_t\gamma^\top R_t^\top A_j^\top \vec p_i) + O(|R_t\gamma-\E_tR_t\gamma|^2)
		\\&= \E_j f(\vec k -\E_t\gamma^\top R_t^\top A_j^\top \vec p_i) + O(\alpha^2).
		\end{align*}
	\end{itemize}
	The theorem is concluded by gathering the new constants into $\bar w_i$.
\end{proof}

\begin{lemma}\label{thm: sym lemma}
	Under \ref{ass: symmetric},
	$$\int_{\F R^2} \vec kf(\vec k - \vec c)f(\vec k+\vec c)\D\vec k = 0$$
	for all $\vec c$.
\end{lemma}
\begin{proof}
	The proof is direct:
	\begin{align*}
	\int_{\F R^2} \vec k f(\vec k-\vec c)f(\vec k+\vec c)\D\vec k &= \int_{\F R^2} -\vec k f(-\vec k-\vec c)f(-\vec k+\vec c)|\op{det}(-\op{id})|\D\vec k & \edit{}{\vec k\mapsto -\vec k}
	\\&\edit{}{= -\int_{\F R^2} \vec k f(-(\vec k+\vec c))f(-(\vec k-\vec c))\D\vec k }
	\\&= -\int_{\F R^2} \vec k f(\vec k+\vec c)f(\vec k-\vec c)\D\vec k & \text{\ref{ass: symmetric}}
	\end{align*}
	\edit{}{therefore $2\int_{\F R^2} \vec k f(\vec k-\vec c)f(\vec k+\vec c)\D\vec k = 0$ as required.}
\end{proof}

\begin{proof}[Proof of Lemma~\ref{thm: final model exact case}.]
	To compute centres of mass, it shall be convenient to define some new functions and abbreviations. We define the functions:
	$$F_0(\vec c) = \int_{\F R^2} f(\vec k-\sfrac12\vec c)f(\vec k+\sfrac12\vec c)\D\vec k, \qquad F_1(\vec c) = \int_{\F R^2} \vec kf(\vec k-\sfrac12\vec c)f(\vec k+\sfrac12\vec c)\D\vec k$$
	and the points
	$$ \vec q_i = \gamma^\top \vec p_i, \qquad \vec q_{i,j} = \gamma^\top  A^\top _j\vec p_i$$
	to \edit{remove excessive use of the $\gamma^\top $ projection}{simplify notation}. Also note by Lemma~\ref{thm: sym lemma} that $F_1(\vec c) =0$ for all $\vec c$. We now compute the relevant integrals starting with the formula of \eqref{eq: dbar def}:
	$$\int_{|\vec k - \vec q_i|<\bar r} \bar D_\alpha(\vec k)\D\vec k = \int_{|\vec k - \vec q_i|<\bar r} \left|\sum_{i'\in\F N}\E_{\gamma^\top \vec \beta_j=\vec0}\left[\bar w_{i'}f(\vec k-\gamma^\top \vec q_{i',j})\right]\right|^2\D\vec k.$$
	Using \ref{ass: separated spots}, we know that the only $i'$ for which the summand is non-zero on this integral domain is $i'=i$. Also using \ref{ass: separated spots}, we know that $f(\vec k-\vec q_{i,j})=0$ outside \edit{of }{}the \edit{restricted}{integral} domain and so we drop this constraint to simplify notation.
	$$\int_{|\vec k - \vec q_i|<\bar r} \bar D_\alpha(\vec k)\D\vec k = \int_{|\vec k - \vec q_i|<\bar r} \left|\E_{\gamma^\top \vec \beta_j=\vec0}\left[\bar w_{i}f(\vec k-\vec q_{i,j})\right]\right|^2\D\vec k = \int_{\F R^2} \left|\E_{\gamma^\top \vec \beta_j=\vec0}\left[\bar w_{i}f(\vec k-\vec q_{i,j})\right]\right|^2\D\vec k.$$
	Next, we expand the brackets noting that $f$ is a real\edit{}{-}valued function:
	\begin{multline*}
	\int_{|\vec k - \vec q_i|<\bar r} \bar D_\alpha(\vec k)\D\vec k = \int_{\F R^2} |\bar w_i|^2 \E_{\substack{\gamma^\top \vec \beta_{j}=\vec0\\\gamma^\top \vec \beta_{J}=\vec0}} \left[f(\vec k-\vec q_{i,j})f(\vec k-\gamma^\top \vec q_{i,J}) \D\vec k\right] 
	\\= |\bar w_i|^2 \E_{\substack{\gamma^\top \vec \beta_{j}=\vec0\\\gamma^\top \vec \beta_{J}=\vec0}} \left[\int_{\F R^2} f(\vec k-\vec q_{i,j})f(\vec k-\gamma^\top \vec q_{i,J}) \D\vec k\right].
	\end{multline*}
	The coordinate translation $\vec k\mapsto \vec k - \frac12(\vec q_{i,j} + \vec q_{i,J})$ simplifies this to a special case of $F_0$:
	\begin{multline*}
	\int_{|\vec k - \vec q_i|<\bar r} \bar D_\alpha(\vec k)\D\vec k = |\bar w_i|^2 \E_{\substack{\gamma^\top \vec \beta_{j}=\vec0\\\gamma^\top \vec \beta_{J}=\vec0}} \left[\int_{\F R^2} f(\vec k-\sfrac12(\vec q_{i,j}-\vec q_{i,J}))f(\vec k+\sfrac12(\vec q_{i,j}-\vec q_{i,J})) \D\vec k\right]
	\\ = |\bar w_i|^2 \E_{\substack{\gamma^\top \vec \beta_{j}=\vec0\\\gamma^\top \vec \beta_{J}=\vec0}} F_0(\vec q_{i,j}-\vec q_{i,J}).
	\end{multline*}
	Similarly, 
	\begin{align*}
	\int_{|\vec k - \vec q_i|<\bar r} \vec k &\bar D_\alpha(\vec k)\D\vec k = \int_{|\vec k - \vec q_i|<\bar r} \vec k\left|\sum_{i'\in\F N}\E_{\gamma^\top \vec \beta_j=\vec0}\left[\bar w_{i'}f(\vec k-\vec q_{i',j})\right]\right|^2\D\vec k
	\\&= |\bar w_i|^2\E_{\substack{\gamma^\top \vec \beta_{j}=\vec0\\\gamma^\top \vec \beta_{J}=\vec0}}\left[\int_{\F R^2} \vec k f(\vec k - \vec q_{i,j})f(\vec k - \vec q_{i,J})\D\vec k\right]
	\\&= |\bar w_i|^2\E_{\substack{\gamma^\top \vec \beta_{j}=\vec0\\\gamma^\top \vec \beta_{J}=\vec0}}\left[\int_{\F R^2} \left(\vec k + \tfrac{\vec q_{i,j}+\vec q_{i,J}}{2}\right) f\left(\vec k - \tfrac{\vec q_{i,j}-\vec q_{i,J}}{2}\right)f\left(\vec k + \tfrac{\vec q_{i,j}-\vec q_{i,J}}{2}\right)\D\vec k\right]
	\\&= |\bar w_i|^2\E_{\substack{\gamma^\top \vec \beta_{j}=\vec0\\\gamma^\top \vec \beta_{J}=\vec0}} \left[F_1(\vec q_{i,j}-\vec q_{i,J}) + \frac{\vec q_{i,j}+\vec q_{i,J}}{2}F_0(\vec q_{i,j}-\vec q_{i,J})\right]
	\\&= |\bar w_i|^2\E_{\substack{\gamma^\top \vec \beta_{j}=\vec0\\\gamma^\top \vec \beta_{J}=\vec0}} \left[\frac{\vec q_{i,j}+\vec q_{i,J}}{2}F_0(\vec q_{i,j}-\vec q_{i,J})\right]
	\end{align*}
	Finally, as $A_{j}+A_{J}$ is independent of $A_{j}-A_{J}$ we can translate this to $\vec q_{i,j}/\vec q_{i,J}$ and again apply Lemma~\ref{thm: product of independent}:
	\begin{equation}\label{eq: com square approximation}
	\E_{j,J} \left[\frac{\vec q_{i,j}+\vec q_{i,J}}{2}F_0(\vec q_{i,j}-\vec q_{i,J})\right] = \E_{j,J} \left[\frac{\vec q_{i,j}+\vec q_{i,J}}{2}\right]\E_{j,J}\left[F_0(\vec q_{i,j}-\vec q_{i,J})\right].\end{equation}
	Thus
	\begin{align*}
	\frac{\int_{|\vec k-\vec q_i|<\bar r} \vec k \bar D_\alpha(\vec k)\D\vec k}{\int_{|\vec k-\vec q_i|<\bar r} \bar D_\alpha(\vec k)\D\vec k} &= \frac{|\bar w_i|^2\E_{j,J}\left[\frac{\vec q_{i,j}+\vec q_{i,J}}{2}\right]\E_{j,J}\left[F_0(\vec q_{i,j}-\vec q_{i,J})\right]}{|\bar w_i|^2\E_{j,J}\left[F_0(\vec q_{i,j}-\vec q_{i,J})\right]}
	\\&= \frac12\E_{j,J}\left[\vec q_{i,j}+\vec q_{i,J}\right]
	= \E_{\gamma^\top \vec \beta_j=0}\vec q_{i,j} = \E_{\gamma^\top \vec \beta_j=0}\gamma^\top A_j^\top \vec p_i.
	\end{align*}
\end{proof}

\section{Precession angle estimation}\label{app: precession rule-of-thumb}
Figure~\ref{fig: precession angle} sketches the choice of precession angle for a deformed sample. There are two triangles of interest, one in the positive $z$ direction which accounts for the curvature of the sphere, and another in the negative direction to account for strain. The upper triangle is a right-angled triangle with one corner at the origin and the other at a point where the sphere of radius $P$ meets the Ewald sphere, say at point $(k,k_z(k))$. This gives the relationship
\begin{align*}
P^2 &= k^2 + k_z(k)^2
\\&= k^2 + \left[4\pi^2\lambda^{-2} - 4\pi\lambda^{-1}\sqrt{4\pi^2\lambda^{-2}-k^2} + |4\pi^2\lambda^{-2}-k^2|\right]
\\&=4\pi\lambda^{-1}\edit{k_z(k)}{\left[2\pi\lambda^{-1}-\sqrt{4\pi^2\lambda^{-2}-k^2}\right]}
\end{align*}
i.e. $k_z(k) = \frac{\lambda P^2}{4\pi}$. On the other hand, the strain moves the point $(P,0)$ a maximal distance of $\sigma P$ from its starting point and away from the Ewald sphere. Assuming the worst strain is a rotation, we get an isosceles triangle whose angle can be computed with the cosine rule:
$$ \cos(\theta) = \frac{2P^2-\sigma^2 P^2}{2P^2}\edit{ = 1 - \frac{\sigma^2}{2}}{}.$$
Combining\edit{}{ this with the sine rule}, the maximal angle is 
$$\alpha \edit{}{= \cos^{-1}\left(\frac{2P^2-\sigma^2 P^2}{2P^2}\right) + \sin^{-1}\left(\frac{k_z(k)}{P}\right)} = \cos^{-1}\left(1-\frac{\sigma^2}{2}\right) + \sin^{-1}\left(\frac{\lambda P}{4\pi}\right).$$
\begin{figure}[h]\begin{center}
	\begin{tikzpicture}[scale=5]
	\draw (0,-.7) -- (0,.6);
	\draw (0,0) -- (1,0);
	\draw (0,0) arc (270:320:1.4) node[right] {Ewald sphere};
	\filldraw (1,0) circle (1pt) node[right] {$P$};
	\draw[->, thick] (1,0) arc (0:20:1) node[midway, right] {curvature};
	\draw[->, thick] (1,0) arc (0:-35:1) node[midway, right] {strain};
	\draw[-, thick] (0,0) -- node[midway, above] {$P$} (.92,.35) node[above] {$\frac{\lambda P^2}{4\pi}$} -- (.92,0);
	\draw[-,thick] (0,0) -- node[midway, below] {$P$}(.81,-.57) -- node[midway, left] {$\sigma P$}(1,0);
	\end{tikzpicture}
\end{center}
\caption{Geometrical argument for choosing precession angle}\label{fig: precession angle}
\end{figure}

\end{appendix}

\SUPP{\newpage\pagenumbering{roman}\setcounter{figure}{0} \renewcommand\thefigure{S.\arabic{figure}}
\begin{center}{\textbf{\Huge Supplementary}}\end{center}
\section*{S.1.\  Continuous deformation phantom}\customlabel{supp: continuous deformation}{S.1}
Examples of strained discs for the layered phantoms are given in Figure~\ref{fig: examples} however the discs for the continuously deformed phantom looked qualitatively different. An example is given in Figure~\ref{fig: continuous deformation disc}.
\begin{figure}[h]\customlabel{fig: continuous deformation disc}{S.1}
	\centering
	\includegraphics[width=\textwidth]{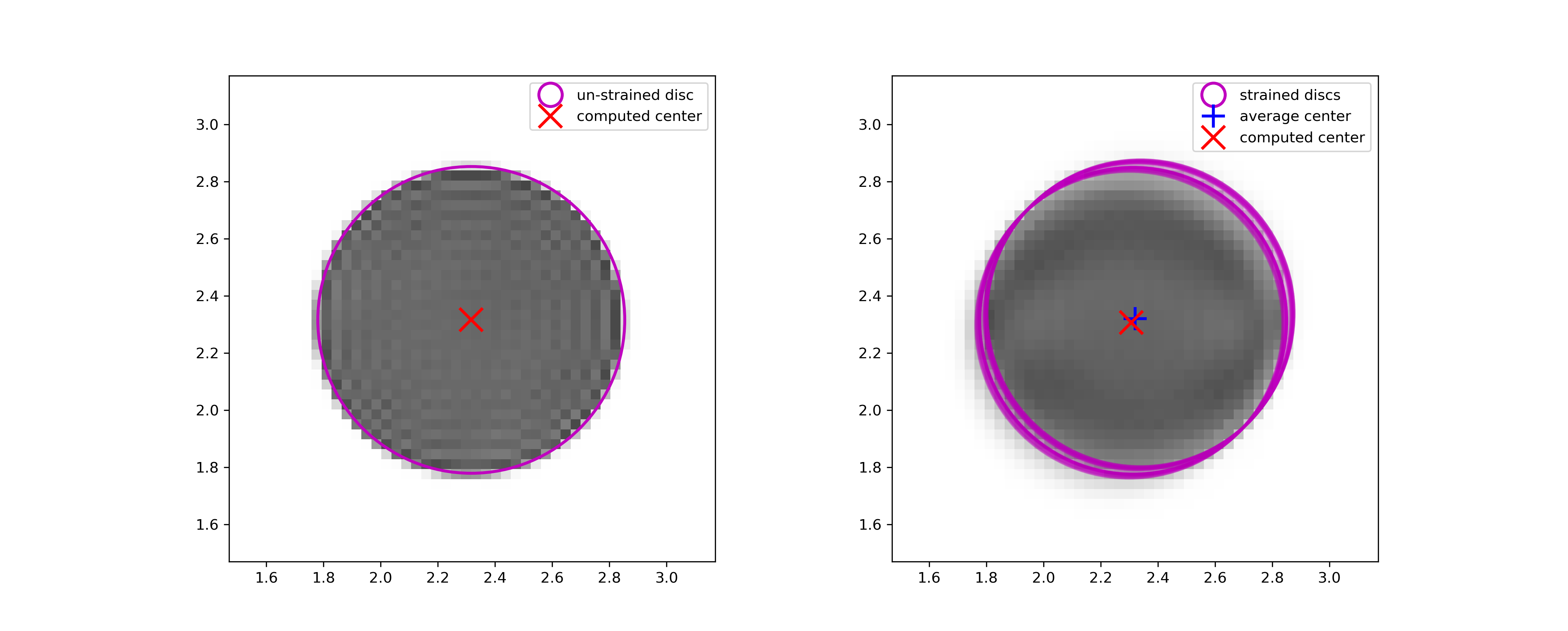}
	\caption{Example disc with worst centre of mass prediction error for the continuous deformation phantom. Left/right hand plots show un-strained/strained disc respectively.}
\end{figure}

\section*{S.2.\ Precession discretisation}\customlabel{supp: precession discretisation}{S.2}
For the dynamical simulations, multiple simulations were run with different numbers of points discretising precession, and the corresponding values in Table~\ref{tbl: quantitative results} are plotted in Figure~\ref{fig: precession points plot}. Dynamical results stated in Table~\ref{tbl: quantitative results} are with $2^8$ precession points. Because of the number of kinematical simulations, the full plots were not computed. Instead, as $\alpha=2^\circ$ was the slowest to converge for dynamical, the number of precession points was increased until the values stated in Table~\ref{tbl: quantitative results} for kinematical simulation at $\alpha=2^\circ$ had converged within the necessary two decimal place rounding error. Kinematical results stated in Table~\ref{tbl: quantitative results} are with $2^5$ precession points.
\begin{figure}[h]\customlabel{fig: precession points plot}{S.2}
	\centering
	\includegraphics[width=.95\textwidth,trim={10 5 10 5}, clip]{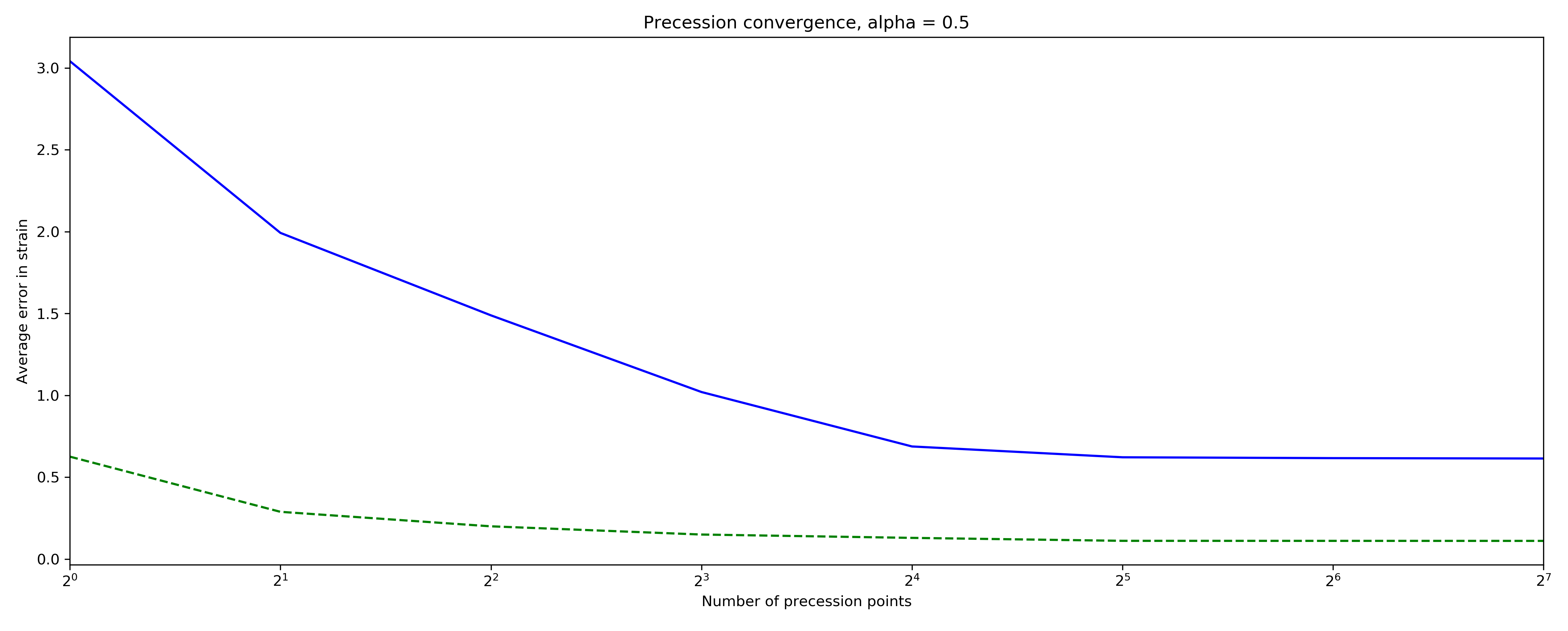}
	\includegraphics[width=.95\textwidth,trim={10 5 10 5}, clip]{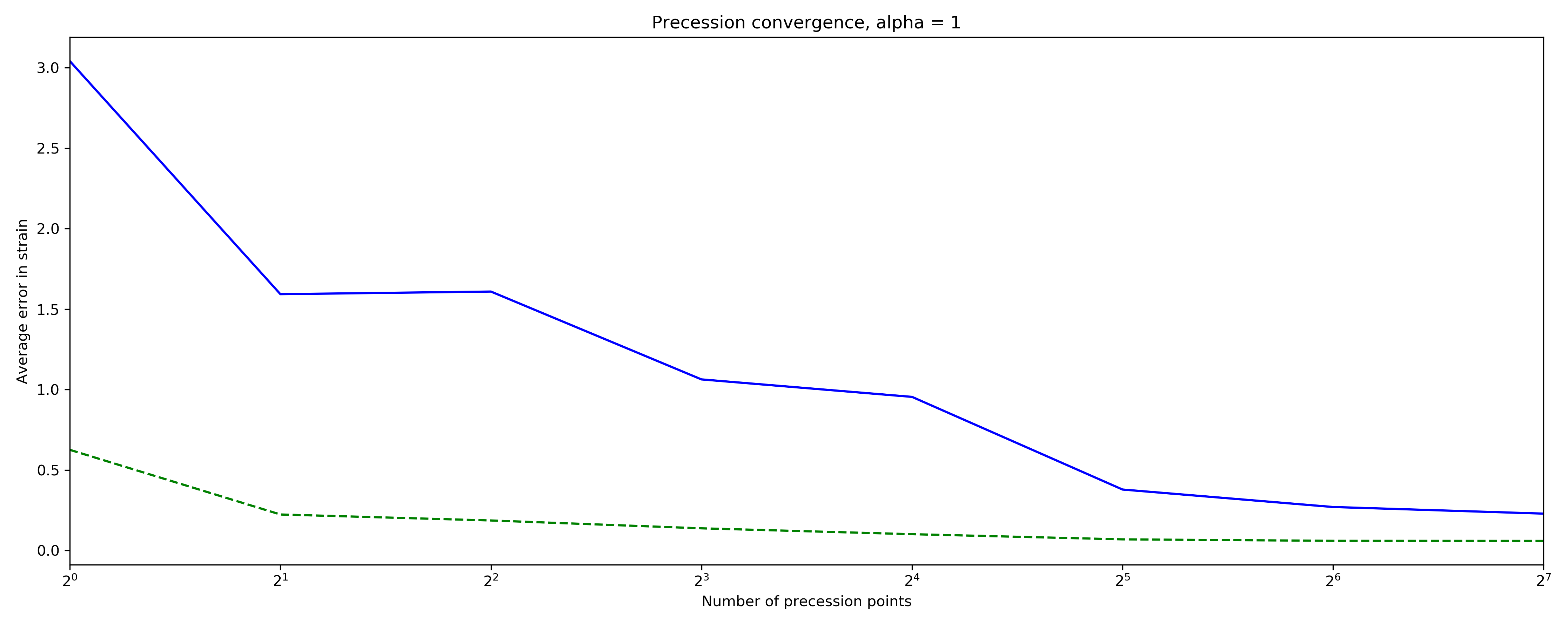}
	\includegraphics[width=.95\textwidth,trim={10 5 10 5}, clip]{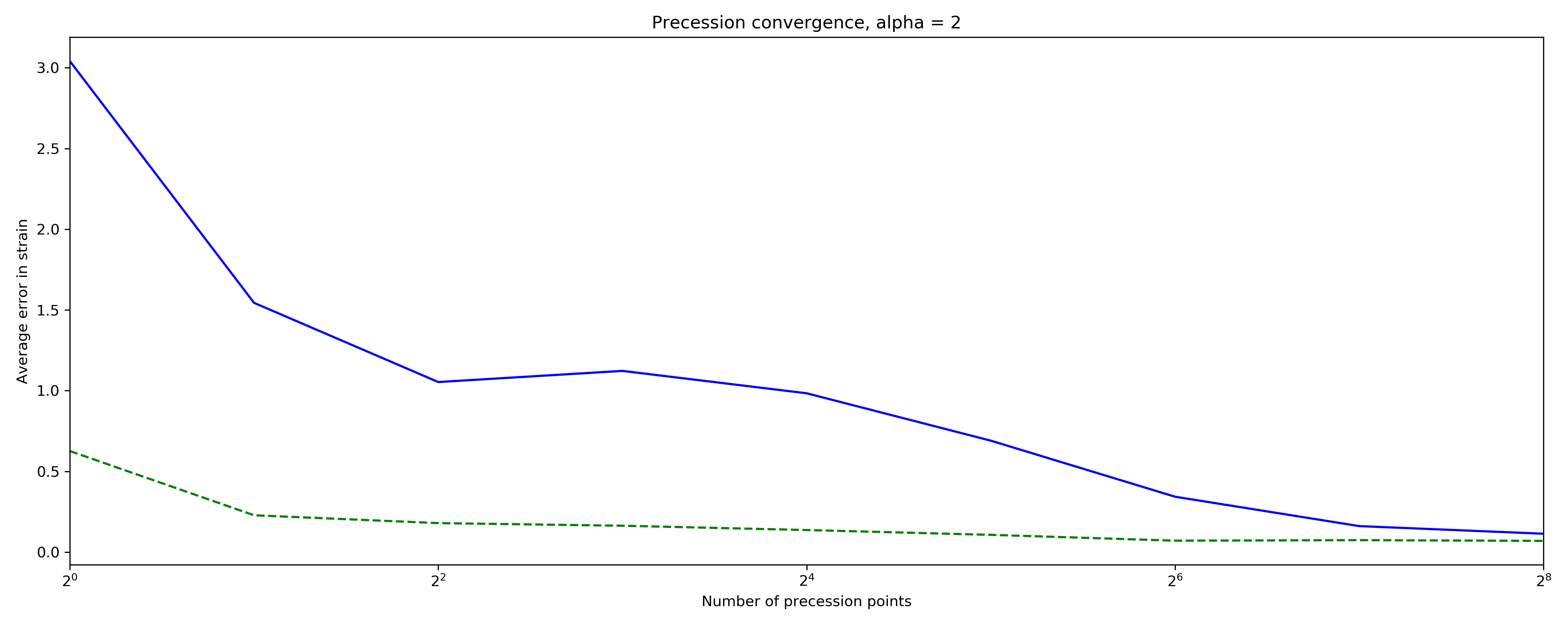}
	\caption{Convergence plot of the values in Table~\ref{tbl: quantitative results} for dynamical simulation for different numbers of precession points. The solid and dashed lines shows the convergence of centre of mass accuracy and registration method respectively.}
\end{figure}

\section*{S.3.\ Reconstruction error plots}\customlabel{supp: reconstruction error plots}{S.3}
In Section~\ref{sec: synth recon} a reconstructed gradient deformation tensor is computed and the 99\textsuperscript{th} percentile errors are reported. Figure~\ref{fig: phantom recon slices comparison} demonstrates the spread of errors in more detail. The first column represents the middle slice (a best-case analysis) which shows that errors can be up to a factor of four smaller than in the worst case (maximum shown in third column). The 99\textsuperscript{th} percentile gives an average worst-case. Qualitatively, this agrees much more with the middle slice indicating that the largest errors are achieved on the interface voxels between crystal and vacuum and on the smooth interior errors are much lower, up to a factor of three. 

\begin{figure}\customlabel{fig: phantom recon slices comparison}{S.3}\centering
	\begin{tikzpicture}
		\node at (0,0) {\includegraphics[height=.3\textheight, trim={20 10 10 5},clip]{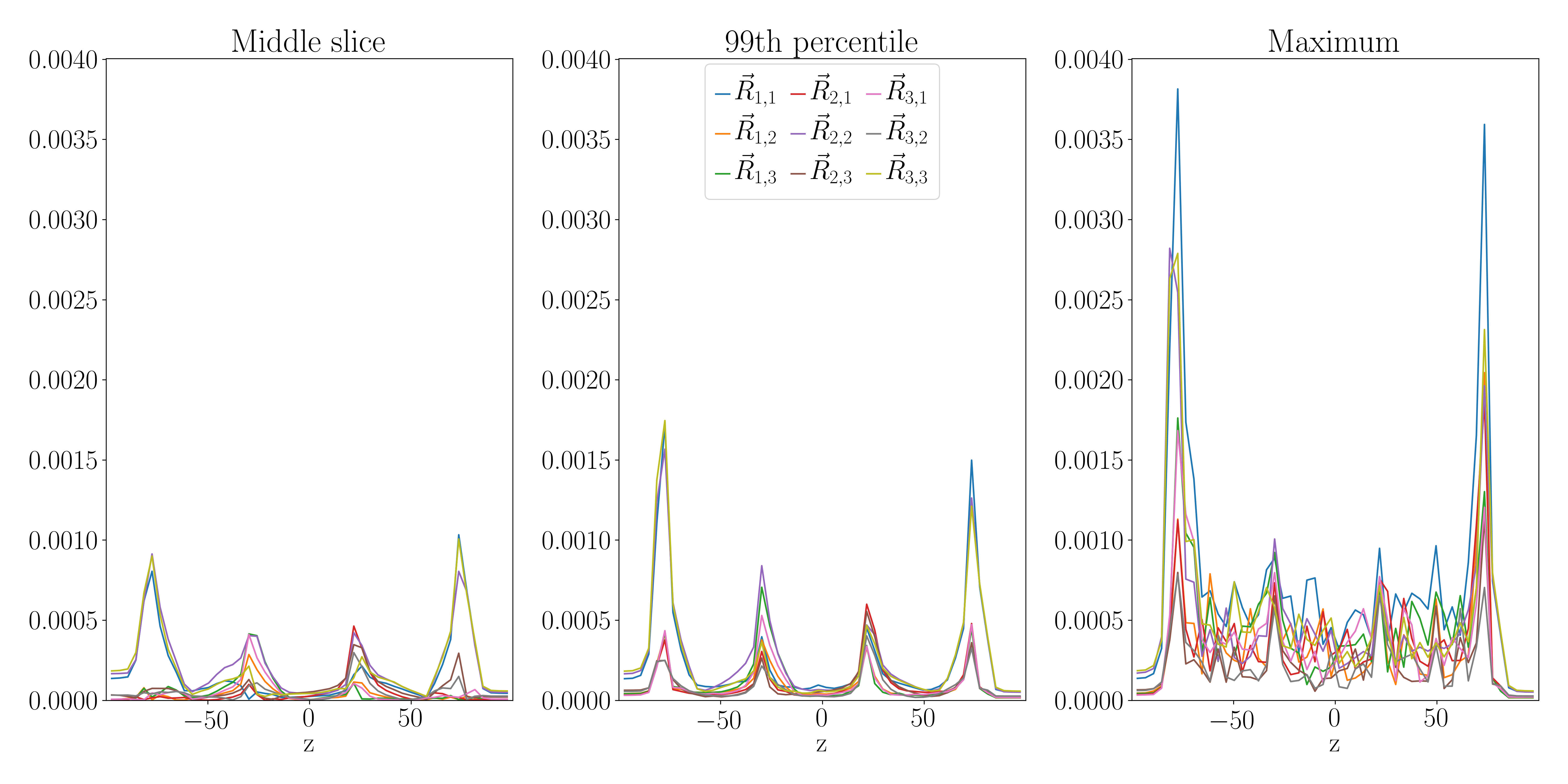}};
		\node at (-.5\textwidth,0) [rotate=90,text width=100pt, align=center] {Low noise};
		\node at (0,-.3\textheight) {\includegraphics[height=.3\textheight, trim={20 10 10 5},clip]{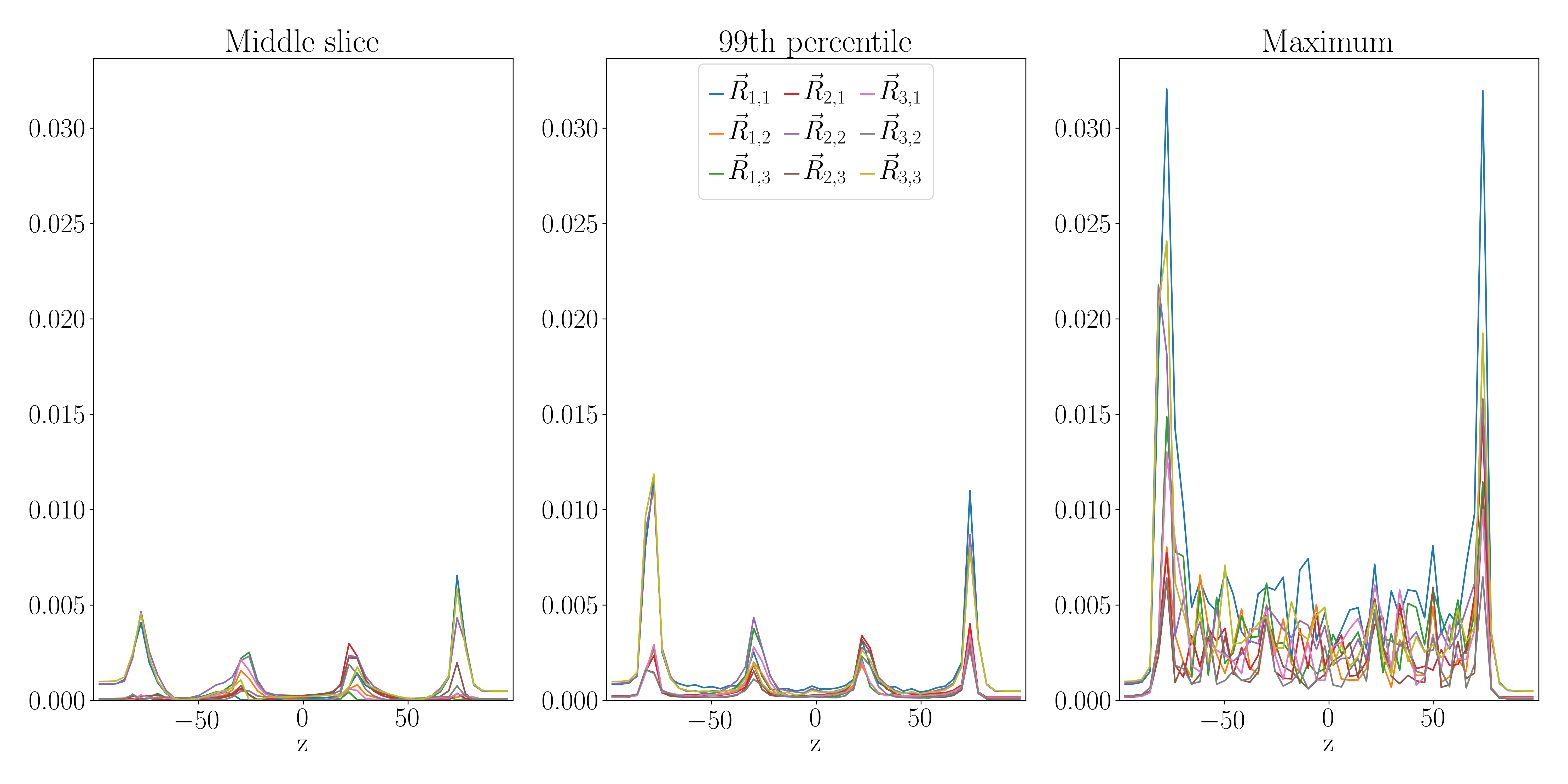}};
		\node at (-.5\textwidth,-.3\textheight) [rotate=90,text width=100pt, align=center] {High noise};
	\end{tikzpicture}
	\caption{1D projections of error of the gradient deformation tensor. For each 1D point shown, the corresponding 2D slice is projected to the reported error by the indicated method. The first pixel extracts the physically central pixel, the second computes the 99\textsuperscript{th} percentile, and the final computes the maximum over all pixels.}
\end{figure}}{}
\end{document}